\documentclass[12pt]{article}
\usepackage{subfigure}
\usepackage{tikz}
\usepackage{pstricks, amssymb, xr, multicol,textcomp,hyperref, amsmath} 
\definecolor{persianblue}{rgb}{0.11, 0.22, 0.73}
\hypersetup{colorlinks=true,
    linkcolor=persianblue,
    filecolor=black,      
    urlcolor=persianblue,
    citecolor=persianblue}
\usepackage{accents}
\usepackage{dirtytalk}
\usepackage{graphicx}
\def\hfl#1{\smash{\mathop{\hbox to 10mm{\rightarrowfill}}\limits^{\textstyle
#1}}}

\newtheorem{proposition}[equation]{Proposition} 
\newtheorem{corollary}[equation]{Corollary} 
\newtheorem{theorem}[equation]{Theorem} 
\newtheorem{exa}[equation]{Example} 
\newtheorem{ex}[equation]{Exercise} 
\newtheorem{s-ex}[equation]{Side-exercise} 
\newtheorem{exas}[equation]{Examples} 
\newtheorem{lemma}[equation]{Lemma} 
 
\newtheorem{remar}[equation]{Remark} 
\newtheorem{remars}[equation]{Remarks} 
\newtheorem{nota}[equation]{Notation} 
\newtheorem{sremar}[equation]{Side-remark} 
\newtheorem{definitio}[equation]{Definition}

\newenvironment{notation}{\begin{nota} \rm }{\end{nota}} 
\newenvironment{remark}{\begin{remar} \rm }{\end{remar}}

\newenvironment{examples}{\begin{exas} \rm }{\end{exas}} 
\newenvironment{example}{\begin{exa} \rm }{\end{exa}} 
\newenvironment{definition}{\begin{definitio} \rm }{\end{definitio}} 
\newcommand{\cconst}{c}
\newcommand{\petita}{a}
\newcommand{\petitb}{b}
\newcommand{\petitc}{c}
\newcommand{\granda}{A}
\newcommand{\grandb}{B}
\newcommand{\grandc}{C}

\newcommand{\talpha}{\tilde{\alpha}}
\newcommand{\alpham}{\alpha_-}
\newcommand{\calpham}{\check{\alpha}_-}
\newcommand{\tcalpham}{\tilde{\check{\alpha}}_-}
\newcommand{\talpham}{\tilde{\alpha}_-}
\newcommand{\betam}{\beta_-}
\newcommand{\tbeta}{\tilde{\beta}}
\newcommand{\tbetam}{\tilde{\beta}_-}
\newcommand{\cbetam}{\check{\beta}_-}
\newcommand{\gammam}{\gamma_-}
\newcommand{\ua}{u_{\alpha}}
\newcommand{\ub}{u_{\beta}}
\newcommand{\uc}{u_{\gamma}}

\newcommand{\amin}{a_-}
\newcommand{\bmin}{b_-}
\newcommand{\cmin}{c_-}
\newcommand{\apl}{a_+}
\newcommand{\bpl}{b_+}
\newcommand{\cpl}{c_+}

\newcommand{\tua}{\tilde{u}_{\alpha}}
\newcommand{\tub}{\tilde{u}_{\beta}}
\newcommand{\tuc}{\tilde{u}_{\gamma}}
\newcommand{\ddel}{\delta} 
\newcommand{\cdel}{\delta} 
\newcommand{\fdel}{\delta} 
\newcommand{\unit}{{\cal U}}

\newcommand{\CA}{{\cal A}}
\newcommand{\CB}{{\cal B}}
\newcommand{\CD}{{\cal D}}
\newcommand{\CE}{{\cal E}}
\newcommand{\CH}{{\cal H}}
\newcommand{\CL}{{\cal L}}
\newcommand{\CP}{{\cal P}}
\newcommand{\CR}{{\cal R}}
\newcommand{\Pcanp}{P} 
\newcommand{\CQ}{{\cal Q}}
\newcommand{\ssih}{{\mbox{\scriptsize Sh}}}
\newcommand{\sih}{{\mbox{Sh}}}
\newcommand{\fHone}{f}
\newcommand{\rats}{R} 
\newcommand{\ambe}{E}
\newcommand{\ambx}{X}
\newcommand{\extx}{X}
\newcommand{\ExtH}{E} 
\newcommand{\hb}{H}
\newcommand{\hbz}{H_0}
\newcommand{\mhbz}{\mathring{H}_0}
\newcommand{\Jknot}{J}
\newcommand{\Kknot}{K}
\newcommand{\Link}{L}

\newcommand{\alambicp}{\CA_{\ambe,\CP}}
\newcommand{\wdel}{w_{\delta}}
\newcommand{\ExtL}{X} 
\newcommand{\ZZ}{\mathbb{Z}} 
\newcommand{\RR}{\mathbb{R}} 
\newcommand{\QQ}{\mathbb{Q}} 
 
\newcommand{\NN}{\mathbb{N}} 
\newcommand{\bp}{\noindent {\sc Proof: }} 
 
\newcommand{\eop}{\nopagebreak \hspace*{\fill}{$\Box$} \medskip} 
\newcommand{\boxa}{\begin{tikzpicture}
\useasboundingbox (-.3,.1) rectangle (.3,.55);
\draw (-.2,0) rectangle (.2,.4) (0,.2) node{\scriptsize $\petita$};
\end{tikzpicture}}
\newcommand{\boxaband}{\begin{tikzpicture}
\useasboundingbox (-.3,0) rectangle (.3,.65);
\draw (-.2,.05) rectangle (.2,.35) (0,.2) node{\scriptsize $\petita$};
\draw (-.1,.05) -- (-.1,-.1) -- (.1,-.1) -- (.1,.05) (-.1,.35) -- (-.1,.5) -- (.1,.5) -- (.1,.35);
\end{tikzpicture}}
\newcommand{\smoothpctrig}{\begin{tikzpicture} \useasboundingbox (-.35,.1) rectangle (.35,.55);
\draw [out=90,in=-90] (.2,0) to (-.2,.4);
\draw [out=90,in=-90,draw=white,double=black,very thick] (-.2,0) to (.2,.4);
\draw [out=90,in=-90] (-.2,0) to (.2,.4);
\draw [draw=gray,->] (.1,.45) .. controls (0,.5) .. (-.1,.45);
\end{tikzpicture}}
\newcommand{\smoothpctrigtube}{\begin{tikzpicture} \useasboundingbox (-.35,.1) rectangle (.35,.55);
\draw [very thick, out=90,in=-90] (.2,0) to (-.2,.4);
\draw [out=90,in=-90,draw=white,double=black,very thick] (-.2,0) to (.2,.4);
\draw [very thick, out=90,in=-90] (-.2,0) to (.2,.4);
\draw (-.2,0) .. controls (-.2,-.04) and (-.1,-.09) .. (0,-.09) .. controls 
(.1,-.09)  and  (.2,-.04)  .. (.2,0);
\draw [dashed] (-.2,0) .. controls (-.2,.04) and (-.1,.09) .. (0,.09) .. controls 
(.1,.09)  and  (.2,.04)  .. (.2,0);
\begin{scope}[yshift=.4cm]
 \draw (-.2,0) .. controls (-.2,-.04) and (-.1,-.09) .. (0,-.09) .. controls 
(.1,-.09)  and  (.2,-.04)  .. (.2,0);
\draw [dashed] (-.2,0) .. controls (-.2,.04) and (-.1,.09) .. (0,.09) .. controls 
(.1,.09)  and  (.2,.04)  .. (.2,0);
\end{scope}
\draw [very thin] (-.2,0) -- (-.2,.4) (.2,0) -- (.2,.4);
\end{tikzpicture}}
\newcommand{\smoothpctrigtubeneg}{\begin{tikzpicture} \useasboundingbox (-.35,.1) rectangle (.35,.55);
\draw [very thick, out=90,in=-90] (-.2,0) to (.2,.4);
\draw [out=90,in=-90,draw=white,double=black,very thick] (.2,0) to (-.2,.4);
\draw [very thick, out=90,in=-90] (.2,0) to (-.2,.4);
\draw (-.2,0) .. controls (-.2,-.04) and (-.1,-.09) .. (0,-.09) .. controls 
(.1,-.09)  and  (.2,-.04)  .. (.2,0);
\draw [dashed] (-.2,0) .. controls (-.2,.04) and (-.1,.09) .. (0,.09) .. controls 
(.1,.09)  and  (.2,.04)  .. (.2,0);
\begin{scope}[yshift=.4cm]
 \draw (-.2,0) .. controls (-.2,-.04) and (-.1,-.09) .. (0,-.09) .. controls 
(.1,-.09)  and  (.2,-.04)  .. (.2,0);
\draw [dashed] (-.2,0) .. controls (-.2,.04) and (-.1,.09) .. (0,.09) .. controls 
(.1,.09)  and  (.2,.04)  .. (.2,0);
\end{scope}
\draw [very thin] (-.2,0) -- (-.2,.4) (.2,0) -- (.2,.4);
\end{tikzpicture}}
\newcommand{\tubea}{\begin{tikzpicture} \useasboundingbox (-.35,.1) rectangle (.35,.55);
\draw (0,.18) node{\scriptsize $\petita$} (-.2,0) .. controls (-.2,-.04) and (-.1,-.09) .. (0,-.09) .. controls 
(.1,-.09)  and  (.2,-.04)  .. (.2,0);
\draw [dashed] (-.2,0) .. controls (-.2,.04) and (-.1,.09) .. (0,.09) .. controls 
(.1,.09)  and  (.2,.04)  .. (.2,0);
\begin{scope}[yshift=.4cm]
 \draw (-.2,0) .. controls (-.2,-.04) and (-.1,-.09) .. (0,-.09) .. controls 
(.1,-.09)  and  (.2,-.04)  .. (.2,0);
\draw [dashed] (-.2,0) .. controls (-.2,.04) and (-.1,.09) .. (0,.09) .. controls 
(.1,.09)  and  (.2,.04)  .. (.2,0);
\end{scope}
\draw [very thin] (-.2,0) -- (-.2,.4) (.2,0) -- (.2,.4);
\end{tikzpicture}}
\newcommand{\smoothnc}{\begin{tikzpicture}
\useasboundingbox (-.35,.1) rectangle (.35,.55);
\draw [out=90,in=-90] (-.2,0) to (.2,.4);
\draw [out=90,in=-90,draw=white,double=black,very thick] (.2,0) to (-.2,.4);
\draw [out=90,in=-90] (.2,0) to (-.2,.4);
\end{tikzpicture}}

\begin{document} 
\title{On elementary invariants of genus one knots and Seifert surfaces
}

\author{Christine Lescop \thanks{Institut Fourier, CNRS, Univ. Grenoble Alpes}}
\maketitle
\begin{abstract}
This elementary article introduces easy-to-manage invariants of genus one knots in homology $3$-spheres. To prove their invariance, we investigate properties
of an invariant of $3$-dimensional genus two homology handlebodies called the \emph{Alexander form}.
The Alexander form of a $3$-manifold $E$ with boundary contains all Reidemeister torsions of link exteriors obtained by attaching $2$-handles along the boundary of $E$. It is a useful tool for studying Alexander polynomials and Reidemeister torsions. We extract invariants of genus one Seifert surfaces from the Alexander form of their exteriors.
\end{abstract}

\medskip 


\medskip 

\tableofcontents
\section{Introduction}

\paragraph{Preliminary conventions and definitions} In this article,
we consider manifolds of dimension at most three, and all manifolds are oriented. We orient boundaries of manifolds with the \emph{outward normal first} convention. Since any topological manifold of dimension at most three has a unique smooth structure up to diffeomorphism, we work with smooth (oriented) manifolds.
An \emph{isotopy} of such a manifold $M$ is a smooth map from $\Psi \colon [0,1] \times M$ to $M$ such that $\Psi(0,.)$ is the identity map and $\Psi(t,.)$ is a diffeomorphism for any $t \in [0,1]$.
Two submanifolds $P$ and $Q$ of $M$ are said to be \emph{isotopic} if there exists an isotopy of $M$ such that $\Psi\left( \{1\} \times P\right)=Q$.
We say that a function of submanifolds of a manifold $M$ is an \emph{invariant} if it maps any two isotopic submanifolds to the same value.
A \emph{knot} (resp. a \emph{link}) in a $3$-manifold is an isotopy class of (images of) embeddings of the circle (resp. a disjoint union of circles) into the manifold. In this article, a \emph{Seifert surface} is a connected compact (oriented) surface with one boundary component.
For $\Lambda$ = $\ZZ$, $\ZZ/2\ZZ$, or $\QQ$, a \emph{$\Lambda$-sphere} is a compact oriented $3$-manifold with the same homology with coefficients in $\Lambda$ as the standard $3$-sphere $S^3$. $\ZZ$-spheres (resp. $\QQ$-spheres) are also called \emph{integer (resp. rational) homology $3$-spheres}. We omit the $3$ and call them integer or rational homology spheres.
Unless otherwise mentioned, homology coefficients are in $\ZZ$. The cardinality or \emph{order} of the $H_1$ of a $\QQ$-sphere $\rats$ is denoted by $|H_1(\rats)|$.

\subsection{An overview}

Here is a short and vague introduction for readers familiar with Reidemeister torsions and Alexander polynomials. We provide complete definitions and background material in the following subsections.

We introduce an invariant $\wdel$ of genus one knots in homology $3$-spheres in Theorem~\ref{thmsimpleinvt}. To my knowledge, it is a new invariant. It is a simple combination of coefficients of Alexander polynomials of curves of a genus one Seifert surface.
Then we introduce a second independent simple invariant $w_{SL}$ of genus one knots  in $\ZZ$-spheres in Section~\ref{subSato}.
As we show in Theorem~\ref{thmw3decomp}, this second invariant $w_{SL}$ is a combination of $\wdel$ and a (degree $3$) knot invariant $w_3$, which appears in a surgery formula for a (degree $2$) invariant of $\QQ$-spheres.
To prove the invariance of $\wdel$, we show that $\wdel$ is a combination of $w_3$ and a fourth invariant $W_{\CD}$ of genus one knots, which is also new to my knowledge.

Let us outline the construction of $W_{\CD}$.
Let $K$ be a knot that bounds a genus one Seifert surface $\Sigma$ in a $\QQ$-sphere $\rats$. Assume that $[-1,1] \times \Sigma$ is embedded in $\rats$ so that $\Sigma=\{0\} \times \Sigma$. Let $\ExtH[K]$ be the manifold obtained from the \emph{exterior} $\ExtH=\rats \setminus \bigl(\left]-1,1\right[ \times \mathring{\Sigma}\bigr)$ of $\Sigma$ by adding a $2$-handle along $K$.
We extract a tautological invariant $W_s(\CD_2(\Sigma))$ of $\Sigma$ from a normalized Reidemeister torsion $\CD(\Sigma)$ of $\ExtH[K]$ in a canonical way.
Then we prove that $W_s(\CD_2(\Sigma))$ depends only on $K$, as outlined in Section~\ref{subintrosk}. Our fourth invariant $W_{\CD}(K)$ is $W_s(\CD_2(\Sigma))$.

To study these invariants and their relationships, we use and investigate \emph{Alexander forms} of $3$-manifolds. Alexander forms are generalizations of abelian Reidemeister torsions from link complements to more general $3$-manifolds.
We define them from the homology of the maximal free abelian covering spaces in Section~\ref{subdefAlef}. Alexander forms contain Reidemeister torsions (and Alexander polynomials) of link exteriors obtained by attaching $2$-handles.

A \emph{genus $g$ handlebody} is a manifold obtained from the unit ball $B^3$ of $\RR^3$ by attaching $g$ one-handles ($D^2 \times [-1,1]$ along $D^2 \times \partial [-1,1]$) to the boundary $\partial B^3$ of $B^3$.
For $\Lambda$ = $\ZZ$ or $\QQ$, a \emph{(genus $g$) $\Lambda$-handlebody} is a compact oriented $3$-manifold with the same homology with coefficients in $\Lambda$ as a (genus $g$) handlebody.
A (genus $g$) \emph{$\QQ$-handlebody} is also called a (genus $g$) \emph{rational homology handlebody}.

The exterior of a genus one Seifert surface in a $\QQ$-sphere is an example of a genus two rational homology handlebody.
We exhibit structural properties of the Alexander forms of genus two rational homology handlebodies in Theorem~\ref{thmAlexHg}, which is of independent interest.

\subsection{A first simple invariant of genus one knots}

The simplest (nontrivial) link invariant is the linking number.
It applies to $2$-component links $(J, K)$ in rational homology $3$-spheres.
In such manifolds, every null-homologous knot bounds a Seifert surface.\footnote{To be precise, we should take an embedding image representing the knot here, but it would make the text heavier.} The \emph{linking number} $lk(J, K)$ of two disjoint oriented knots $J$ and $K$ is the algebraic intersection of $J$ with a transverse Seifert surface $\Sigma_K$ of $K$ when $K$ is null-homologous. In general, $lk(J, K)$ is the algebraic intersection of $J$ with a rational chain bounded by $K$. 
The \emph{Alexander polynomial} $\Delta(K)$ is a simple invariant of knots in $\QQ$-spheres discovered in 1928 by Alexander.
It is the order of the first homology group of the maximal free abelian covering space of the knot complement (see Definition~\ref{defAlexonevorder}). Let us give an alternative definition of $\Delta(K)$, for a null-homologous knot $K$ in a $\QQ$-sphere $\rats$, using a Seifert surface $\Sigma_K$ of $K$, as in \cite[Proposition 2.3.13, page 27]{lespup}, for example. Embed $[-1,0] \times \Sigma_K$ in $\rats$ by an orientation-preserving embedding that naturally identifies $\{0\} \times \Sigma_K$ with $\Sigma_K$. Set $\Sigma_{K-} = \{-1\} \times \Sigma_K$. For a curve $v$, let $v_-$ denote the curve of $\{-1\} \times v$. 
The \emph{Seifert form} associated with $H_1(\Sigma_K)$ is the form 
$$\begin{array}{llll}V \colon &H_1(\Sigma_K) \otimes_s H_1(\Sigma_K)&\rightarrow 
&\ZZ\\
&[u] \otimes_s [v] &\mapsto& lk(u,v_-).
\end{array}$$
Then we have 
$$\Delta(K)=|H_1(\rats)|\mbox{det}\left(t^{1/2} V- t^{-1/2} (^TV)\right)$$
where $V$ also denotes a \emph{Seifert matrix} of the form $V$, and $^TV$ is its transpose.
We extract the (even simpler) invariant $$\lambda^{\prime}(K) = \frac{\Delta^{\prime\prime}(K)(1)}{2|H_1(\rats)|}$$ from $\Delta(K)$.
When $K$ has a genus one Seifert surface $\Sigma_K$, $\lambda^{\prime}(K)$ is simply the determinant of an associated Seifert matrix $V$. (We will give definitions of $\Delta$ and $\lambda^{\prime}$ for non-necessarily null-homologous knots in $\QQ$-spheres in Definitions~\ref{defAlexonevorder} and \ref{deflambdaprimeJknot}. See also Lemma~\ref{lemlambdaprimeJknot}.)

Let $\Sigma$ be a genus one Seifert surface. A \emph{symplectic basis} of $H_1(\Sigma)$ is a basis $(\alpha,\beta)$ of $H_1(\Sigma)$ such that the \emph{algebraic intersection number} $\langle \alpha,\beta\rangle_{\Sigma}$ is one. We say that an element $u$ of $H_1(\Sigma)$ is \emph{primitive} if there exists an element $v$ of $H_1(\Sigma)$ such that $\langle u,v\rangle_{\Sigma}=1$. Any primitive element of $H_1(\Sigma)$ can be represented by a non-separating simple closed curve of $\Sigma$.
We will prove the following easy lemma in Section~\ref{secpfisotopcurv}.

\begin{lemma}\label{lemhomcurvesisot}
 Two homologous non-separating simple closed curves of a genus one oriented surface $\Sigma$ with one boundary component are isotopic in $\Sigma$.
\end{lemma}

Thus, for such a genus one Seifert surface, any primitive element $u$ of $H_1(\Sigma)$ is represented by a connected curve of $\Sigma$, whose isotopy class in $\Sigma$ is determined by $u$. We also denote this curve associated to $u$ by $u$. We say that a Seifert surface $\Sigma$ in a rational homology sphere $\rats$ is  \emph{null-homologous} if the map induced by the inclusion from $H_1(\Sigma)$ to $H_1(\rats)$ maps $H_1(\Sigma)$ to $0$. 

We are now ready to define a simple invariant of genus one knots in integer homology $3$-spheres and state one of the main theorems of this article.

\begin{theorem}
\label{thmsimpleinvt}
Let $\Sigma=\Sigma_K$ be a
genus one Seifert surface of a knot $K$ in a $\QQ$-sphere $\rats$, let $(\alpha,\beta)$ be a symplectic basis of $H_1(\Sigma)$.
Let $\alpha$, $\beta$, and $\gamma=-\alpha-\beta$ denote both the corresponding primitive elements of $H_1(\Sigma)$ and the associated curves in $\Sigma$ and in the ambient manifold $\rats$.
Define $(\petita,\petitb,\petitc)$ so that
$$lk(\alpha,\alpham)=\frac{\petitb+\petitc}{2}\mbox{,}\;\;\;\;lk(\beta,\betam)=\frac{\petitc+\petita}{2}\mbox{,}\;\;\;\;\mbox{and}\;\;\;\;lk(\gamma,\gammam)=\frac{\petita+\petitb}{2}.$$
Set $$\lambda^{\prime}(\petita,\petitb,\petitc)= \frac{\petita\petitb+\petita\petitc+\petitb\petitc+1}{4}$$
and 
$$\wdel(\petita,\petitb,\petitc)=\frac{\petita+\petitb}2\frac{\petitb+\petitc}2\frac{\petitc+\petita}2-\lambda^{\prime}(\petita,\petitb,\petitc)(\petita+\petitb+\petitc).$$

Then $\lambda^{\prime}(K)= \lambda^{\prime}(\petita,\petitb,\petitc)$ and $$\wdel(\Sigma)=4\petita \lambda^{\prime}(\alpha) +4\petitb \lambda^{\prime}(\beta) +4\petitc \lambda^{\prime}(\gamma) +\wdel(\petita,\petitb,\petitc)$$ is an invariant of $\Sigma$.\footnote{It is independent of our choice of the symplectic basis
$(\alpha,\beta)$ of $H_1(\Sigma)$.}
Furthermore, if
$\Sigma$ is null-homologous, then $\wdel(\Sigma)$ depends only on $K=\partial \Sigma$. We denote it by $\wdel(K)$.
\end{theorem}

\begin{remark}
\label{rkdefabc}
Note that $(\beta,\gamma)$ and $(\gamma,\alpha)$ are symplectic bases of $H_1(\Sigma)$. The expression of $\wdel$ is invariant under a cyclic permutation of $(a,\alpha)$, $(b,\beta)$, and $(c,\gamma)$. We can equivalently define
 $\petita$ to be
 $$\petita = lk(\beta,\betam) + lk(\gamma,\gammam) - lk(\alpha,\alpham)$$
 or $$\petita=-lk(\beta,\gammam)-lk(\gamma,\betam),$$
where $lk(\beta,\gammam)-lk(\gamma,\betam)=\langle \beta,\gamma \rangle_{\Sigma}=1$.
Cyclic permutations lead to similar equalities for $\petitb$ and $\petitc$.
These equalities imply that $\petita$, $\petitb$, and $\petitc$ are odd integers when $\Sigma$ is null-homologous.
 \end{remark}
 
 \begin{remark}
\label{rkapgenusone}
The Alexander polynomial $\Delta(K)$ of a genus one knot is $\Delta(K) =1 +\lambda^{\prime}(K)(t -2+t^{-1})$.
 \end{remark}

\begin{figure}[h]
\begin{center}
\begin{tikzpicture} 
\draw [thick, rounded corners, fill=blue!20] (3.4,2.2) -- (3.4,3.8) -- (0,3.8) -- (0,2.2)
(0,2.2) -- (.6,2.2) 
(.6,2.2) -- (.6,3) -- (1.4,2.4) -- (1.4,2.2)
(1.4,2.2) -- (2,2.2)
(2,2.2) -- (2,2.4) -- (2.8,3) -- (2.8,2.2)
(2.8,2.2) -- (3.4,2.2);
\draw [rounded corners, fill=yellow]
(0,1.6) -- (0,0) -- (3.4,0) -- (3.4,1.6) (3.4,1.6) -- (2.8,1.6)
(2.8,1.6) -- (2.8,.8) -- (2,1.4) -- (2,1.6) (2,1.6) -- (1.4,1.6)
(1.4,1.6) -- (1.4,1.4) -- (.6,.8) -- (.6,1.6) (.6,1.6) -- (0,1.6);
\draw [rounded corners, thick]
(0,1.6) -- (0,0) -- (3.4,0) -- (3.4,1.6)
(2.8,1.6) -- (2.8,.8) -- (2,1.4) -- (2,1.6)
(1.4,1.6) -- (1.4,1.4) -- (.6,.8) -- (.6,1.6);
\draw [rounded corners,->] (3.1,3.5) node{\scriptsize $\gamma$} (2.92,3.5) -- (2.92,3.68) -- (.48,3.68) -- (.48,2.2) (2.92,2.2) -- (2.92,3.5);
\draw [dashed, rounded corners,->] (2.6,.3) node{\scriptsize $\gamma$} (2.6,.12) --  (3.28,.12) -- (3.28,1.6) (.12,1.6) -- (.12,.12) -- (2.6,.12);
\draw [rounded corners,->] (1.1,3.3) node{\scriptsize $\alpha$} (1.2,3.12) -- (1.52,3.12) -- (1.88,2.84) -- (1.88,2.2) (.12,2.2) -- (.12,3.12) -- (1.2,3.12) ;
\draw [dashed, rounded corners,-<] (1.12,.78)  node{\scriptsize $\alpha$}  (1.12,.98) -- (1.52,1.28) -- (1.52,1.6) (.48,1.6) -- (.48,.68) -- (.72,.68) -- (1.12,.98) ;
\draw [rounded corners,->] (2.3,3.3) node{\scriptsize $\beta$} (2.2,3.12) -- (3.28,3.12) -- (3.28,2.2) (1.52,2.2) -- (1.52,2.84) -- (1.88,3.12) -- (2.2,3.12);
\draw [dashed, rounded corners,->] (2.15,.82)  node{\scriptsize $\beta$}  (2.28,.98) -- (1.88,1.28) -- (1.88,1.6) (2.92,1.6) -- (2.92,.68) -- (2.68,.68) -- (2.28,.98);
\draw [thick,->] (.8,2.85) -- (1,2.7);
\draw [thick,->] (2.2,2.55) -- (2.4,2.7);
\draw [thick,->] (2,0) -- (1.5,0);
\draw [thick,->] (2,3.8) -- (1.7,3.8);
\draw [thick,->] (2.2, 1.25) -- (2.4,1.1);
\draw [thick,->] (1.7,3.75) node[above]{\scriptsize $K(\petita,\petitb,\petitc)$} (.8,.95) -- (1,1.1);
\draw (.3,1.9) node{\scriptsize $\petitb$} (-.1,1.6)  rectangle (.7,2.2);
\draw (1.7,1.9) node{\scriptsize $\petitc$} (2.1,1.6) rectangle (1.3,2.2);
\draw (3.1,1.9) node{\scriptsize $\petita$} (3.5,1.6) rectangle (2.7,2.2);
\fill[color=white] (.1 ,3.25) rectangle (.4 ,3.65) ;
       \draw (.25 ,3.45) node{\scriptsize $\Sigma$}; 
\fill[color=white] (1 ,.26) rectangle (2.4 ,.66) ;
       \draw (1.7 ,.44) node{\scriptsize $-\Sigma(\petita,\petitb,\petitc)$}; 

\begin{scope}[xshift=5cm]
\fill [rounded corners,blue!20] 
       (3.4,2.2) -- (3.4,3.8) -- (0,3.8) -- (0,2.11) -- (.3,1.85) -- (.6,2.11) -- (.6,3.3) -- (1.4,3.3) -- (1.4,2.11) -- (1.7,1.85) -- (2,2.11) -- (2,2.5) -- (2.8,2.5) -- (2.8,2.11) -- (3.1,1.85) -- (3.4,2.11);
       \fill[color=white] (1.5 ,3.1) rectangle (3.3 ,3.5) ;
       \draw (2.4,3.3) node{\scriptsize $\Sigma(-1,-1,1)$};
\fill [rounded corners, yellow] (0,1.6) -- (0,0) -- (2,0) -- (3.4,.0) -- (3.4,1.7) -- (3.1,1.95) -- (2.8,1.7) -- (2.8,1.3) -- (2,1.3) -- (2,1.7) -- (1.7,1.95) -- (1.4,1.7) -- (1.4,.5) -- (.6,.5) -- (.6,1.7) -- (.3,1.95) -- (0,1.7);
   \fill[color=white] (1.5 ,.3) rectangle (3.3 ,.7) ;
       \draw (2.4 ,.5) node{\scriptsize $-\Sigma(-1,-1,1)$}; 
\draw [rounded corners, thick] (3.4,2.2) -- (3.4,3.8) -- (0,3.8) -- (0,2.2) 
(0,1.6) -- (0,0) -- (2,0) -- (3.4,.0) -- (3.4,1.6)
(2.8,2.2) -- (2.8,2.5) -- (2,2.5) -- (2,2.2)
(2.8,1.6) -- (2.8,1.3) -- (2,1.3) -- (2,1.6)
(1.4,2.2) -- (1.4,3.3) -- (.6,3.3) -- (.6,2.2)
(1.4,1.6) -- (1.4,.5) -- (.6,.5) -- (.6,1.6);
\draw [thick,->] (1.7,3.75) node[above]{\scriptsize $K(-1,-1,1)$} (2,3.8) -- (1.7,3.8);
\draw [thick,->] (.8,3.3) -- (1,3.3);
\draw [thick,->] (1.5,0) -- (1,0);
\draw [thick,->] (.8,.5) -- (1,.5);
\draw [thick,out=90,in=-90] (2.8,1.6) to (3.4,2.2) ;
\draw [out=90,in=-90,draw=white,double=black,very thick] (3.4,1.6) to (2.8,2.2);
\draw [thick,out=90,in=-90] (3.4,1.6) to (2.8,2.2);
\draw [thick,out=90,in=-90] (2,1.6) to (1.4,2.2);
\draw [out=90,in=-90,draw=white,double=black,very thick] (1.4,1.6) to (2,2.2);
\draw [thick,out=90,in=-90] (1.4,1.6) to (2,2.2);
\draw [thick,out=90,in=-90] (0,1.6) to (.6,2.2);
\draw [out=90,in=-90,draw=white,double=black,very thick] (.6,1.6) to (0,2.2);
\draw [thick,out=90,in=-90] (.6,1.6) to (0,2.2);
\end{scope}
\begin{scope}[xshift=10cm]

\draw [rounded corners, thick] (3.4,2.2) -- (3.4,3.8) -- (0,3.8) -- (0,2.8) 
(0,1) -- (0,0) -- (2,0) -- (3.4,.0) -- (3.4,1.6)
(2.8,2.2) -- (2.8,2.5) -- (2,2.5) -- (2,2.2)
(2.8,1.6) -- (2.8,1.3) -- (2,1.3) -- (2,1.6)
(1.4,2.2) -- (1.4,3.3) -- (.6,3.3) -- (.6,2.8)
(1.4,1.6) -- (1.4,.5) -- (.6,.5) -- (.6,1);
\draw [thick,->] (1.7,3.75) node[above]{\scriptsize $K(-1,3,-1)$} (2,3.8) -- (1.7,3.8);
\draw [thick,->] (.8,3.3) -- (1,3.3);
\draw [thick,->] (1.5,0) -- (1,0);
\draw [thick,->] (.8,.5) -- (1,.5);
\draw [thick,out=90,in=-90] (2.8,1.6) to (3.4,2.2) ;
\draw [out=90,in=-90,draw=white,double=black,very thick] (3.4,1.6) to (2.8,2.2);
\draw [thick,out=90,in=-90] (3.4,1.6) to (2.8,2.2);
\draw [thick,out=90,in=-90] (1.4,1.6) to (2,2.2);
\draw [out=90,in=-90,draw=white,double=black,very thick] (2,1.6) to (1.4,2.2);
\draw [thick,out=90,in=-90] (2,1.6) to (1.4,2.2);
\draw [thick,out=90,in=-90] (.6,1) to (0,1.6);
\draw [out=90,in=-90,draw=white,double=black,very thick] (0,1) to (.6,1.6);
\draw [thick,out=90,in=-90] (0,1) to (.6,1.6);
\draw [thick,out=90,in=-90] (.6,1.6) to (0,2.2);
\draw [out=90,in=-90,draw=white,double=black,very thick] (0,1.6) to (.6,2.2);
\draw [thick,out=90,in=-90] (0,1.6) to (.6,2.2);
\draw [thick,out=90,in=-90] (.6,2.2) to (0,2.8);
\draw [out=90,in=-90,draw=white,double=black,very thick] (0,2.2) to (.6,2.8);
\draw [thick,out=90,in=-90] (0,2.2) to (.6,2.8);
\end{scope}\end{tikzpicture}
\caption{The surface $\Sigma(\petita,\petitb,\petitc)$ with its curves $\alpha$, $\beta$, and $\gamma$, $\Sigma(-1,-1,1)$, and the figure-eight knot $K(-1,3,-1)$}
\label{figSigmaabc}
\end{center}
\end{figure}
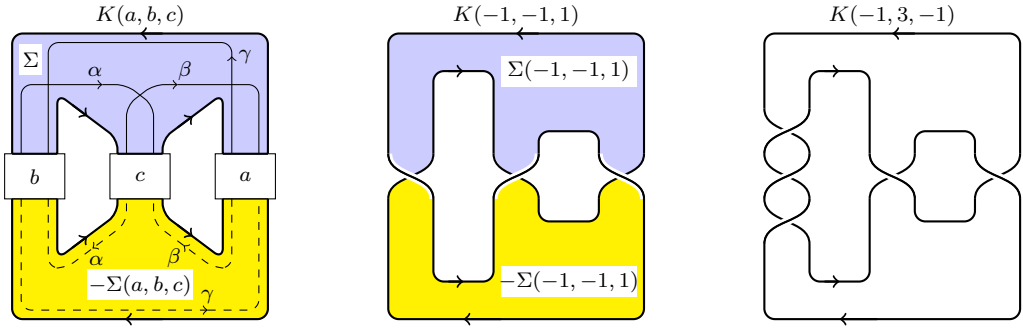

Consider the \emph{pretzel knots} $K(\petita,\petitb,\petitc)$ of Figure~\ref{figSigmaabc}. In this figure, for an integer $\petita$, the box \boxa denotes a two-strand braid 
with $|\petita|$ vertical juxtapositions of the motive \smoothpctrig, if $\petita$ is positive, and $|\petita|$ vertical juxtapositions of \smoothnc otherwise.
The surface 
$\Sigma(\petita,\petitb,\petitc)$ embedded in $(\RR^3\subset S^3)$ as in Figure~\ref{figSigmaabc} is a Seifert surface of $K(\petita,\petitb,\petitc)=\partial \Sigma(\petita,\petitb,\petitc)$.

\begin{examples}
The knot $K(-1,-1,1)$ is the trivial knot, and we have $\lambda^{\prime}(K(-1,1,-1))=0$.
The pictured curves $\alpha$, $\beta$, and $\gamma$ also represent trivial knots. In particular, we have $$\lambda^{\prime}(\alpha)= \lambda^{\prime}(\beta)= \lambda^{\prime}(\gamma)=0,$$
and $\wdel\left( K(\petita,\petitb,\petitc)\right) = \wdel\left(\petita,\petitb,\petitc\right)$.

For the left-handed trefoil $K(1,1,1)$, we have $\wdel(K(1,1,1))=1-3=-2$. For
the right-handed trefoil $K(-1,-1,-1)$, we have $\wdel(K(-1,-1,-1))=2$. (The \emph{signature} of $(V+^TV)$, which is another simple classical knot invariant, also distinguishes these knots.)

In \cite{OSknot}, Ozsváth and Szabó introduced a powerful \emph{categorification} $\widehat{HFK}$ of the Alexander polynomial, which detects the genus of knots. (The Alexander polynomial is a graded Euler characteristic of $\widehat{HFK}$.) Let $k$ be an integer greater than $1$.
According to a theorem of Hedden and Watson \cite[Theorem 1]{HeddenWatson}, the
pretzel knots of the family $\bigl(K(2n+1,2k+1,-2k-1)\bigr)_{n \in \NN}$
have the same Ozsváth--Szabó Heegaard--Floer homologies $\widehat{HFK}$ and $HFK$.\footnote{The pretzel knot $K(1,2k+1,-2k-1)$ is obtained by connecting the components of a trivial two-component link \begin{tikzpicture}
\useasboundingbox (-1.7,.1) rectangle (1.7,.9);
\draw [rounded corners] (-1.3,.7) -- (-1.3,.9) -- (1.3,.9) -- (1.3,.7)  (-.7,.7) -- (-.7,.8) -- (.7,.8) -- (.7,.7)
(-1.3,.3) -- (-1.3,.1) -- (1.3,.1) -- (1.3,.3)  (-.7,.3) -- (-.7,.2) -- (.7,.2) -- (.7,.3);
\fill [white] (-1.7,.3) rectangle (-.3,.7)
(1.7,.3) rectangle (.3,.7);
\draw (-1.7,.3) rectangle (-.3,.7)
(1.7,.3) rectangle (.3,.7)
(-1,.5)
node{\scriptsize $-2k-1$} (1,.5) node{\scriptsize $2k+1$};                                                              \end{tikzpicture}
 by a band $[0,1]^2$ that intersects the link along $[0,1]\times \partial [0,1]$ and replacing $[0,1]\times \partial [0,1]$ with $ (\partial [0,1])\times [0,1]$. It is nontrivial, since
$\lambda^{\prime}(K(1,2k+1,-2k-1))=-k-k^2$.}
We have $w_{\delta}\bigl(K(2n+1,2k+1,-2k-1)\bigr)= (k+k^2)(2n+1)$. So $w_{\delta}$ distinguishes these knots, which have identical Ozsváth--Szabó Heegaard--Floer homology, identical Alexander polynomial, and identical vanishing signature. (The Jones polynomial also distinguishes them.)
\end{examples}

The manifold obtained from a manifold $M$ by reversing its orientation is denoted by $(-M)$. 
\begin{proposition}
 Let $K$ be a knot that bounds a null-homologous genus one Seifert surface in a $\QQ$-sphere $\rats$. Then we have
  $$\wdel(K \subset \rats)=\wdel((-K) \subset \rats) \;\;\;\mbox{and} \;\;\;\wdel(K \subset (-\rats))=-\wdel(K \subset \rats).$$
\end{proposition}
\bp Changing the orientation of $K$ exchanges the roles of $\alpha$ and $\beta$. Since the expression of $\wdel$ is invariant under this exchange,
 the invariant $\wdel$ is independent of the orientation of $K$.\footnote{For a triple $K(\petita,\petitb,\petitc)$ of pairwise distinct odd integers greater than $2$, the knot $K(\petita,\petitb,\petitc)$ is \emph{non-invertible}, i.e., it is not isotopic to $-K(\petita,\petitb,\petitc)=K(\petitb,\petita,\petitc)$. In \cite{Trotternoninv}, Hale Trotter proved the existence of non-invertible knots by proving that these knots are non-invertible. Unfortunately, the invariant $\wdel$ cannot detect non-invertibility.}
 
  Reversing the orientation of the ambient rational homology sphere $\rats$ multiplies the linking numbers by $(-1)$. Therefore this operation, which changes $\rats$ to $(-\rats)$, multiplies $\petita$, $\petitb$, and $\petitc$ by $(-1)$.
\eop

The \emph{mirror image} of the isotopy class of a knot embedding $K \colon S^1 \hookrightarrow \RR^3$ is the class of the composition of $K$ by a reflection of $\RR^3$.
It is denoted by $\overline{K}$.
We have $\wdel(\overline{K})=-\wdel(K)$.
 Therefore, $\wdel$ vanishes for \emph{amphicheiral knots}, which are knots isotopic to their mirror images. The figure-eight knot $K(-1,3,-1)=K(1,-3,1)$
 is an example of an amphicheiral knot such that $\lambda^{\prime}(K(-1,3,-1))=-1$.

\subsection{A second simple invariant of genus one knots}
\label{subSato}
The intersection of two oriented surfaces in a $\QQ$-sphere $\rats$ is oriented so that the triple (positive normal to the first surface, positive normal to the second one, oriented tangent vector to the intersection) induces the orientation of $\rats$.

The definition of our second simple invariant of genus one knots is based on the following Sato--Levine invariant \cite{Sato}.

\begin{definition} \label{defzetatwo} Let $\Jknot$ and $\Kknot$
be two disjoint null-homologous knots in a rational homology sphere $\rats$, such that the linking number $lk(\Jknot,\Kknot)$ is zero. Let $\Sigma_{\Jknot}$ and $\Sigma_{\Kknot}$ be two transverse surfaces in $\rats$ respectively bounded by $\Jknot$ and $\Kknot$, such that $\Sigma_{\Jknot} \cap \Sigma_{\Kknot}$ is disjoint from $\Jknot$ and $\Kknot$. Let $\Sigma_{\Jknot,\parallel}$ be obtained from $\Sigma_{\Jknot}$ by a slight normal push.  Define the \emph{Sato--Levine invariant} $\lambda^{\prime}(\Jknot,\Kknot)$ of $(\Jknot,\Kknot)$ to be
$$\lambda^{\prime}(\Jknot,\Kknot)=-lk (\Sigma_{\Jknot} \cap \Sigma_{\Kknot},\Sigma_{\Jknot,\parallel} \cap \Sigma_{\Kknot}).$$  See Definition~\ref{deflambdaprimetwocomp} and Lemma~\ref{lemzetatwocomp} for more general definitions of $\lambda^{\prime}(\Jknot,\Kknot)$.  
\end{definition}

Let $\hb(\petita,\petitb,\petitc)$ be the genus two handlebody $[-1,0] \times \Sigma(\petita,\petitb,\petitc)$, where we identify $\Sigma(\petita,\petitb,\petitc)$ with $\{0\} \times \Sigma(\petita,\petitb,\petitc)$. 

The boundary $\partial \hb(\petita,\petitb,\petitc)$ of $\hb(\petita,\petitb,\petitc)$ is the union of $\Sigma(\petita,\petitb,\petitc)$ and $\left(-\Sigma(\petita,\petitb,\petitc)_-\right)$, where $$\Sigma(\petita,\petitb,\petitc)_- = \left(\{-1\} \times \Sigma(\petita,\petitb,\petitc)\right) \cup_{\{-1\} \times \partial \Sigma(\petita,\petitb,\petitc)} \left([-1,0] \times K(\petita,\petitb,\petitc)\right).$$ 
We think of $\partial \hb(\petita,\petitb,\petitc)$ as being obtained from the two copies $\Sigma(\petita,\petitb,\petitc)$ and $\Sigma(\petita,\petitb,\petitc)_-$ of $\Sigma(\petita,\petitb,\petitc)$ by inflating between these two glued copies whose boundaries are identified and sealed as in Figure~\ref{fighandlebo}. (In Figure~\ref{fighandlebo}, the tube \tubea represents $|\petita|$ vertical juxtapositions of the motive \smoothpctrigtube if $\petita$ is positive, and $|\petita|$ vertical juxtapositions of \smoothpctrigtubeneg otherwise.) In particular, $\hbz=\hb(\petita,\petitb,\petitc)$ does not depend on $(\petita,\petitb,\petitc)$, but the inclusion from $\{0\} \times \Sigma(\petita,\petitb,\petitc)$ to $\hbz$ does.

\begin{figure}[h]
\begin{center}
\begin{tikzpicture}
\draw [thick, rounded corners, fill=blue!20] (3.4,2) -- (3.4,3.8) -- (0,3.8) -- (0,2)
(0,2) -- (.6,2) 
(.6,2) -- (.6,3) -- (.8,3) -- (1.4,2.4) -- (1.4,2)
(1.4,2) -- (2,2)
(2,2) -- (2,2.4) -- (2.6,3) -- (2.8,3) -- (2.8,2)
(2.8,2) -- (3.4,2);
\begin{scope}[xshift=-1.3cm,yshift=2.9cm] 
\draw [-] (2,.9) .. controls (1.9,.9) and (1.75,.7) .. (1.75,.5);  
\draw (1.75,.5) .. controls (1.75,.3) and (1.9,.1) .. (2,.1);
\draw [dashed] (2,.9) .. controls (2.1,.9) and (2.25,.7) .. (2.25,.5) .. controls (2.25,.3) and (2.1,.1) .. (2,.1);
\end{scope}
\begin{scope}[xshift=.7cm,yshift=2.9cm] 
\draw [-] (2,.9) .. controls (1.9,.9) and (1.75,.7) .. (1.75,.5);  
\draw (1.75,.5) .. controls (1.75,.3) and (1.9,.1) .. (2,.1);
\draw [dashed] (2,.9) .. controls (2.1,.9) and (2.25,.7) .. (2.25,.5) .. controls (2.25,.3) and (2.1,.1) .. (2,.1);
\end{scope}
\draw [thick, rounded corners, fill=yellow]
(0,1.6) -- (0,0) -- (3.4,0) -- (3.4,1.6) (3.4,1.6) -- (2.8,1.6)
(2.8,1.6) -- (2.8,.8) -- (2.6,.8) -- (2,1.4) -- (2,1.6) (2,1.6) -- (1.4,1.6)
(1.4,1.6) -- (1.4,1.4)-- (.8,.8) -- (.6,.8) -- (.6,1.6) (.6,1.6) -- (0,1.6);
\begin{scope}[xshift=-1.3cm,yshift=-.1cm]
\draw [-] (2,.9) .. controls (1.9,.9) and (1.75,.7) .. (1.75,.5); 
\draw (1.75,.5) .. controls (1.75,.3) and (1.9,.1) .. (2,.1);
\draw [dashed] (2,.9) .. controls (2.1,.9) and (2.25,.7) .. (2.25,.5) .. controls (2.25,.3) and (2.1,.1) .. (2,.1);
\end{scope}
\begin{scope}[xshift=.7cm,yshift=-.1cm] 
\draw [-] (2,.9) .. controls (1.9,.9) and (1.75,.7) .. (1.75,.5); 
\draw (1.75,.5) .. controls (1.75,.3) and (1.9,.1) .. (2,.1);
\draw [dashed] (2,.9) .. controls (2.1,.9) and (2.25,.7) .. (2.25,.5) .. controls (2.25,.3) and (2.1,.1) .. (2,.1);
\end{scope}
\draw  (0,1.6) .. controls (0,1.52) and (.15,1.4) ..  (.3,1.4) .. controls (.45,1.4)  and  (.6,1.52)  .. (.6,1.6);
\fill [white] (0,2.2) -- (0,1.6) .. controls (0,1.52) and (.15,1.4) ..  (.3,1.4) .. controls (.45,1.4)  and  (.6,1.52)  .. (.6,1.6) -- (.6,2.2)  .. controls (.6,2.12)  and (.45,2) .. (.3,2) .. controls (.15,2) and (0,2.12) .. (0,2.2);
\draw [dashed] (0,1.6) .. controls (0,1.68) and (.15,1.8) .. (.3,1.8) .. controls 
(.45,1.8)  and  (.6,1.68)  .. (.6,1.6);
 \draw  (0,2.2) .. controls (0,2.12) and (.15,2) ..  (.3,2)  .. controls (.45,2)  and  (.6,2.12)  .. (.6,2.2);
\draw [dashed] (0,2.2) .. controls (0,2.28) and (.15,2.4) ..  (.3,2.4)  .. controls (.45,2.4)  and  (.6,2.28)  .. (.6,2.2);
\draw (.3,1.6) node{\scriptsize $\petitb$} (0,2.2) -- (0,1.6) (.6,1.6) -- (.6,2.2);
\begin{scope}[xshift=1.4cm]

\draw  (0,1.6) .. controls (0,1.52) and (.15,1.4) ..  (.3,1.4) .. controls (.45,1.4)  and  (.6,1.52)  .. (.6,1.6);
\fill [white] (0,2.2) -- (0,1.6) .. controls (0,1.52) and (.15,1.4) ..  (.3,1.4) .. controls (.45,1.4)  and  (.6,1.52)  .. (.6,1.6) -- (.6,2.2)  .. controls (.6,2.12)  and (.45,2) .. (.3,2) .. controls (.15,2) and (0,2.12) .. (0,2.2);
\draw [dashed] (0,1.6) .. controls (0,1.68) and (.15,1.8) .. (.3,1.8) .. controls 
(.45,1.8)  and  (.6,1.68)  .. (.6,1.6);
 \draw  (0,2.2) .. controls (0,2.12) and (.15,2) ..  (.3,2)  .. controls (.45,2)  and  (.6,2.12)  .. (.6,2.2);
\draw [dashed] (0,2.2) .. controls (0,2.28) and (.15,2.4) ..  (.3,2.4)  .. controls (.45,2.4)  and  (.6,2.28)  .. (.6,2.2);
\draw (.3,1.6) node{\scriptsize $\petitc$} (0,2.2) -- (0,1.6) (.6,1.6) -- (.6,2.2);
\end{scope}

\begin{scope}[xshift=2.8cm]
\draw  (0,1.6) .. controls (0,1.52) and (.15,1.4) ..  (.3,1.4) .. controls (.45,1.4)  and  (.6,1.52)  .. (.6,1.6);
\fill [white] (0,2.2) -- (0,1.6) .. controls (0,1.52) and (.15,1.4) ..  (.3,1.4) .. controls (.45,1.4)  and  (.6,1.52)  .. (.6,1.6) -- (.6,2.2)  .. controls (.6,2.12)  and (.45,2) .. (.3,2) .. controls (.15,2) and (0,2.12) .. (0,2.2);
\draw [dashed] (0,1.6) .. controls (0,1.68) and (.15,1.8) .. (.3,1.8) .. controls 
(.45,1.8)  and  (.6,1.68)  .. (.6,1.6);
 \draw  (0,2.2) .. controls (0,2.12) and (.15,2) ..  (.3,2)  .. controls (.45,2)  and  (.6,2.12)  .. (.6,2.2);
\draw [dashed] (0,2.2) .. controls (0,2.28) and (.15,2.4) ..  (.3,2.4)  .. controls (.45,2.4)  and  (.6,2.28)  .. (.6,2.2);
\draw (.3,1.6) node{\scriptsize $\petita$} (0,2.2) -- (0,1.6) (.6,1.6) -- (.6,2.2);
\end{scope}
\begin{scope}[xshift=4.5cm]
\draw [thick, rounded corners, fill=blue!20] (3.4,2) -- (3.4,3.8) -- (0,3.8) -- (0,2)
(0,2) -- (.6,2) 
(.6,2) -- (.6,3) -- (.8,3) -- (1.4,2.4) -- (1.4,2)
(1.4,2) -- (2,2)
(2,2) -- (2,2.4) -- (2.6,3) -- (2.8,3) -- (2.8,2)
(2.8,2) -- (3.4,2); 
\draw [thick, rounded corners, fill=yellow]
(0,1.6) -- (0,0) -- (3.4,0) -- (3.4,1.6) (3.4,1.6) -- (2.8,1.6)
(2.8,1.6) -- (2.8,.8) -- (2.6,.8) -- (2,1.4) -- (2,1.6) (2,1.6) -- (1.4,1.6)
(1.4,1.6) -- (1.4,1.4)-- (.8,.8) -- (.6,.8) -- (.6,1.6) (.6,1.6) -- (0,1.6);
\begin{scope}[xshift=-1.3cm,yshift=-.1cm] 
\draw [-] (2,.9) .. controls (1.9,.9) and (1.75,.7) .. (1.75,.5); 
\draw (1.75,.5) .. controls (1.75,.3) and (1.9,.1) .. (2,.1);
\draw [dashed] (2,.9) .. controls (2.1,.9) and (2.25,.7) .. (2.25,.5) .. controls (2.25,.3) and (2.1,.1) .. (2,.1);
\end{scope}
\begin{scope}[xshift=.7cm,yshift=-.1cm]
\draw [-] (2,.9) .. controls (1.9,.9) and (1.75,.7) .. (1.75,.5); 
\draw (1.75,.5) .. controls (1.75,.3) and (1.9,.1) .. (2,.1);
\draw [dashed] (2,.9) .. controls (2.1,.9) and (2.25,.7) .. (2.25,.5) .. controls (2.25,.3) and (2.1,.1) .. (2,.1);
\end{scope}
\draw  (0,1.6) .. controls (0,1.52) and (.15,1.4) ..  (.3,1.4) .. controls (.45,1.4)  and  (.6,1.52)  .. (.6,1.6);
\fill [white] (0,2.2) -- (0,1.6) .. controls (0,1.52) and (.15,1.4) ..  (.3,1.4) .. controls (.45,1.4)  and  (.6,1.52)  .. (.6,1.6) -- (.6,2.2)  .. controls (.6,2.12)  and (.45,2) .. (.3,2) .. controls (.15,2) and (0,2.12) .. (0,2.2);
\draw [dashed] (0,1.6) .. controls (0,1.68) and (.15,1.8) .. (.3,1.8) .. controls 
(.45,1.8)  and  (.6,1.68)  .. (.6,1.6);
 \draw  (0,2.2) .. controls (0,2.12) and (.15,2) ..  (.3,2)  .. controls (.45,2)  and  (.6,2.12)  .. (.6,2.2);
\draw [dashed] (0,2.2) .. controls (0,2.28) and (.15,2.4) ..  (.3,2.4)  .. controls (.45,2.4)  and  (.6,2.28)  .. (.6,2.2);
\draw (.3,1.6) node{\scriptsize $\petitb$} (0,2.2) -- (0,1.6) (.6,1.6) -- (.6,2.2);
\begin{scope}[xshift=1.4cm]
\draw  (0,1.6) .. controls (0,1.52) and (.15,1.4) ..  (.3,1.4) .. controls (.45,1.4)  and  (.6,1.52)  .. (.6,1.6);
\fill [white] (0,2.2) -- (0,1.6) .. controls (0,1.52) and (.15,1.4) ..  (.3,1.4) .. controls (.45,1.4)  and  (.6,1.52)  .. (.6,1.6) -- (.6,2.2)  .. controls (.6,2.12)  and (.45,2) .. (.3,2) .. controls (.15,2) and (0,2.12) .. (0,2.2);
\draw [dashed] (0,1.6) .. controls (0,1.68) and (.15,1.8) .. (.3,1.8) .. controls 
(.45,1.8)  and  (.6,1.68)  .. (.6,1.6);
 \draw  (0,2.2) .. controls (0,2.12) and (.15,2) ..  (.3,2)  .. controls (.45,2)  and  (.6,2.12)  .. (.6,2.2);
\draw [dashed] (0,2.2) .. controls (0,2.28) and (.15,2.4) ..  (.3,2.4)  .. controls (.45,2.4)  and  (.6,2.28)  .. (.6,2.2);
\draw (.3,1.6) node{\scriptsize $\petitc$} (0,2.2) -- (0,1.6) (.6,1.6) -- (.6,2.2);
\end{scope}

\begin{scope}[xshift=2.8cm]
\draw  (0,1.6) .. controls (0,1.52) and (.15,1.4) ..  (.3,1.4) .. controls (.45,1.4)  and  (.6,1.52)  .. (.6,1.6);
\fill [white] (0,2.2) -- (0,1.6) .. controls (0,1.52) and (.15,1.4) ..  (.3,1.4) .. controls (.45,1.4)  and  (.6,1.52)  .. (.6,1.6) -- (.6,2.2)  .. controls (.6,2.12)  and (.45,2) .. (.3,2) .. controls (.15,2) and (0,2.12) .. (0,2.2);
\draw [dashed] (0,1.6) .. controls (0,1.68) and (.15,1.8) .. (.3,1.8) .. controls 
(.45,1.8)  and  (.6,1.68)  .. (.6,1.6);
 \draw  (0,2.2) .. controls (0,2.12) and (.15,2) ..  (.3,2)  .. controls (.45,2)  and  (.6,2.12)  .. (.6,2.2);
\draw [dashed] (0,2.2) .. controls (0,2.28) and (.15,2.4) ..  (.3,2.4)  .. controls (.45,2.4)  and  (.6,2.28)  .. (.6,2.2);
\draw (.3,1.6) node{\scriptsize $\petita$} (0,2.2) -- (0,1.6) (.6,1.6) -- (.6,2.2);
\end{scope}
\draw [rounded corners,->]  (2.2,3.65) -- (.15,3.65) -- (.15,.15) -- (3.25,.15) -- (3.25,3.65) -- (2.2,3.65);
\draw [rounded corners,->]  (1.15,2.95) -- (1.55,2.55) -- (1.55,1.25) -- (.95,.65) -- (.45,.65) -- (.45,3.15) -- (.95,3.15) -- (1.15,2.95);
\draw [rounded corners,->]  (2.05,2.75) -- (2.45,3.15) -- (2.95,3.15) -- (2.95,.65) -- (2.45,.65) -- (1.85,1.25) -- (1.85,2.55) -- (2.05,2.75);
\draw (2.2,3.7) node[below]{\scriptsize $\grandc$};
\draw (1.05,3.05) node[right]{\scriptsize $\granda$};
\draw (2.1,2.85) node[left]{\scriptsize $\grandb$};
\end{scope}

\end{tikzpicture}
\caption{The handlebody $\hbz=\hb(\petita,\petitb,\petitc)$ and the curves $\granda$, $\grandb$, and $\grandc$}
\label{fighandlebo}
\end{center}
\end{figure}

Let $\rats$ denote a $\QQ$-sphere. 
Consider the three curves $\granda$, $\grandb$, and $\grandc$ of the right-hand side of Figure~\ref{fighandlebo}. Let  $\phi \colon \hbz \hookrightarrow \rats$ be an embedding that maps $\granda$, $\grandb$, and $\grandc$ to
null-homologous curves in the \emph{exterior} $\ExtH= \rats \setminus \phi\bigl(\mhbz\bigr)$ of $\phi(\hbz)$. 
Any null-homologous genus one Seifert $\Sigma$ in $\rats$ may be written as $\phi\bigl(\Sigma(\petita,\petitb,\petitc)\bigr)$ for such an embedding and for three odd integers $\petita$, $\petitb$, and $\petitc$ (which can be chosen as in Theorem~\ref{thmsimpleinvt} or Remark~\ref{rkdefabc}) so that this hypothesis is fulfilled.
Note that $lk\bigl(\phi(\granda),\phi(\grandb)\bigr)=0$ and that $$\lambda^{\prime}\bigl(\phi(\granda),\phi(\grandb)\bigr)=\lambda^{\prime}\bigl(\phi(\grandb),\phi(\grandc)\bigr)=\lambda^{\prime}\bigl(\phi(\grandc),\phi(\granda)\bigr).$$ 

\begin{theorem}
\label{thmSatoinvariant} Let  $\phi \colon \hbz \hookrightarrow \rats$ be an embedding of $\hbz$ into a $\QQ$-sphere $\rats$ that maps $\granda$ and $\grandb$ to 
null-homologous curves in the exterior of $\phi(\hbz)$.
With the notation of Theorem~\ref{thmsimpleinvt}, set
 $$w_{SL}\Bigl(\phi\bigl( \Sigma(\petita,\petitb,\petitc)\bigr)\Bigr) = \lambda^{\prime}\bigl(\phi(\granda),\phi(\grandb)\bigr) - \frac1{12} \wdel(\petita,\petitb,\petitc).$$ 
 Then $w_{SL}\bigl(\phi\bigl( \Sigma(\petita,\petitb,\petitc)\bigr)\bigr)$ is a topological invariant of the genus (at most) one knot $\partial\phi( \Sigma(\petita,\petitb,\petitc)) \subset \rats$. We denote it by $w_{SL}\bigl(\partial\phi\bigl(\Sigma(\petita,\petitb,\petitc)\bigr)\bigr)$.
\end{theorem}

The invariant $w_{SL}$ is a fortiori an invariant of null-homologous genus one Seifert surfaces in $\QQ$-spheres. It is generalized to all genus one Seifert surfaces in $\QQ$-spheres in Theorem~\ref{thmwsl}.

\begin{examples}
\label{exawsldelta} Let $\phi_0$ denote the trivial embedding of Figure~\ref{fighandlebo} of $\hbz$ into $\RR^3 \subset S^3$.
 For the knots $K(\petita,\petitb,\petitc)=\phi_0\bigl(K(\petita,\petitb,\petitc)\bigr)$, the Sato--Levine invariant $\lambda^{\prime}\bigl(\phi_0(\granda),\phi_0(\grandb)\bigr)$ is zero since $\phi_0(\granda)$ bounds a disk disjoint from $\phi(\grandb)$. Thus $w_{SL}\bigl(K(\petita,\petitb,\petitc)\bigr)= - \frac1{12} \wdel(\petita,\petitb,\petitc)$.
 
 Let $B_a$ be a ball
$B_a$ of $\RR^3$ in a neighborhood of the box \boxa of Figure~\ref{figSigmaabc} that intersects $\phi_0\bigl(\Sigma(\petita,\petitb,\petitc)\bigr)$ as the twisted band $\phi_0\bigl([0,1]^2\bigr) =\boxaband$
so that $\left(\partial B_a\right) \cap \phi_0([0,1]^2) = \phi_0([0,1] \times \partial [0,1])$.
 Change the trivial embedding $\phi_0$ of $\Sigma(\petita,\petitb,\petitc)$ to an embedding $\phi_J$ such that
 \begin{itemize}
  \item $\phi_0$ and $\phi_J$ coincide outside the ball
$B_a$,
\item $\phi_J$ maps the band $[0,1]^2$ to the ball $B_a$
so that $\phi_J([0,1] \times \partial [0,1])=\phi_0([0,1] \times \partial [0,1])$,
\item the self-linking number of $\phi_J(\beta)$ is still $\frac{c+a}{2}$, and 
\item the knot $\phi_J(\beta)$ is a knot $J$.
 \end{itemize}
Figure~\ref{figperturbeight} shows an example of a knot $\phi_J(K(-1,3,-1))$ obtained in this way, where $J$ is the figure-eight knot.
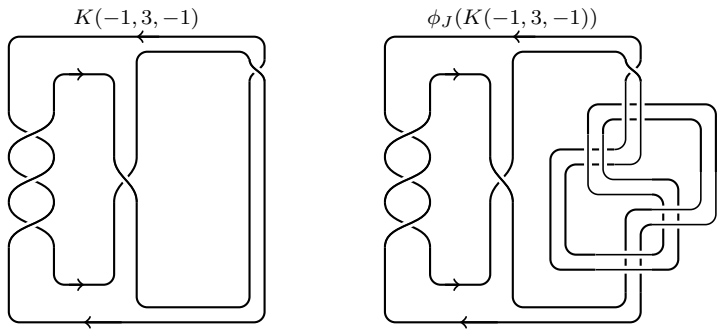
\begin{figure}[h]
\begin{center}
\begin{tikzpicture}
\begin{scope}[xshift=-5cm]
\draw [rounded corners, thick] (3.4,3.5) -- (3.4,3.8) -- (0,3.8) -- (0,2.8) 
(0,1) -- (0,0) -- (3.4,0) -- (3.4,3.2)
(3.2,3.5) -- (3.2,3.6) -- (1.7,3.6) -- (1.7,2.2)
(3.2,3.2) -- (3.2,.2) -- (1.7,.2) -- (1.7,1.6)
(1.4,2.2) -- (1.4,3.3) -- (.6,3.3) -- (.6,2.8)
(1.4,1.6) -- (1.4,.5) -- (.6,.5) -- (.6,1);
\draw [thick,->] (1.7,3.75) node[above]{\scriptsize $K(-1,3,-1)$} (2,3.8) -- (1.7,3.8);
\draw [thick,->] (.8,3.3) -- (1,3.3);
\draw [thick,->] (1.5,0) -- (1,0);
\draw [thick,->] (.8,.5) -- (1,.5);

\draw [thick,out=90,in=-90] (3.2,3.2) to (3.4,3.5) ;
\draw [out=90,in=-90,draw=white,double=black,very thick] (3.4,3.2) to (3.2,3.5);
\draw [thick,out=90,in=-90] (3.4,3.2) to (3.2,3.5);

\draw [thick,out=90,in=-90] (1.4,1.6) to (1.7,2.2);
\draw [out=90,in=-90,draw=white,double=black,very thick] (1.7,1.6) to (1.4,2.2);
\draw [thick,out=90,in=-90] (1.7,1.6) to (1.4,2.2);
\draw [thick,out=90,in=-90] (.6,1) to (0,1.6);
\draw [out=90,in=-90,draw=white,double=black,very thick] (0,1) to (.6,1.6);
\draw [thick,out=90,in=-90] (0,1) to (.6,1.6);
\draw [thick,out=90,in=-90] (.6,1.6) to (0,2.2);
\draw [out=90,in=-90,draw=white,double=black,very thick] (0,1.6) to (.6,2.2);
\draw [thick,out=90,in=-90] (0,1.6) to (.6,2.2);
\draw [thick,out=90,in=-90] (.6,2.2) to (0,2.8);
\draw [out=90,in=-90,draw=white,double=black,very thick] (0,2.2) to (.6,2.8);
\draw [thick,out=90,in=-90] (0,2.2) to (.6,2.8);
\end{scope}

\draw [rounded corners, thick] (3.4,3.5) -- (3.4,3.8) -- (0,3.8) -- (0,2.8) 
(0,1) -- (0,0) -- (3.4,0) -- (3.4,.4)
(3.2,3.5) -- (3.2,3.6) -- (1.7,3.6) -- (1.7,2.2)
(3.2,.4) -- (3.2,.2) -- (1.7,.2) -- (1.7,1.6)
(1.4,2.2) -- (1.4,3.3) -- (.6,3.3) -- (.6,2.8)
(1.4,1.6) -- (1.4,.5) -- (.6,.5) -- (.6,1);
\draw [thick,->] (1.7,3.75) node[above]{\scriptsize $\phi_J(K(-1,3,-1))$} (2,3.8) -- (1.7,3.8);
\draw [thick,->] (.8,3.3) -- (1,3.3);
\draw [thick,->] (1.5,0) -- (1,0);
\draw [thick,->] (.8,.5) -- (1,.5);

\draw [thick,out=90,in=-90] (3.2,3.2) to (3.4,3.5) ;
\draw [out=90,in=-90,draw=white,double=black,very thick] (3.4,3.2) to (3.2,3.5);
\draw [thick,out=90,in=-90] (3.4,3.2) to (3.2,3.5);

\draw [thick,out=90,in=-90] (1.4,1.6) to (1.7,2.2);
\draw [out=90,in=-90,draw=white,double=black,very thick] (1.7,1.6) to (1.4,2.2);
\draw [thick,out=90,in=-90] (1.7,1.6) to (1.4,2.2);

\draw [thick,out=90,in=-90] (.6,1) to (0,1.6);
\draw [out=90,in=-90,draw=white,double=black,very thick] (0,1) to (.6,1.6);
\draw [thick,out=90,in=-90] (0,1) to (.6,1.6);
\draw [thick,out=90,in=-90] (.6,1.6) to (0,2.2);
\draw [out=90,in=-90,draw=white,double=black,very thick] (0,1.6) to (.6,2.2);
\draw [thick,out=90,in=-90] (0,1.6) to (.6,2.2);
\draw [thick,out=90,in=-90] (.6,2.2) to (0,2.8);
\draw [out=90,in=-90,draw=white,double=black,very thick] (0,2.2) to (.6,2.8);
\draw [thick,out=90,in=-90] (0,2.2) to (.6,2.8);

\draw [rounded corners, thick] (4.2,2) -- (4.2,2.7) -- (2.9,2.7) -- (2.9,2.5)
 (4.4,2) -- (4.4,2.9) -- (2.7,2.9) -- (2.7,2.5)
 (3.05,2.3) -- (2.2,2.3) -- (2.2,.7) -- (2.8,.7)
(3.05,2.1) -- (2.4,2.1) -- (2.4,.9) -- (2.8,.9)
(3.2,.4) -- (3.2,1.5) -- (3.55,1.5) 
(3.4,.4) -- (3.4,1.3) -- (3.55,1.3) 
(3.55,.7) -- (3.9,.7) -- (3.9,1.9) -- (3.55,1.9)
(3.55,.9) -- (3.7,.9) -- (3.7,1.7) -- (3.55,1.7);

\draw [rounded corners,draw=white,double=black,very thick] 
(3.2,3.2) -- (3.2,2.3) -- (3.05,2.3)
(3.4,3.2) -- (3.4,2.1) -- (3.05,2.1)
(2.7,2.5) -- (2.7,1.7) -- (3.55,1.7)
(2.9,2.5) -- (2.9,1.9) -- (3.55,1.9)
(4.2,2) -- (4.2,1.5) -- (3.55,1.5)
(4.4,2) -- (4.4,1.3) -- (3.55,1.3)
(2.8,.7) -- (3.55,.7)
(2.8,.9) -- (3.55,.9);
\end{tikzpicture}
\caption{A modification of the figure-eight knot $K(-1,3,-1)$}
\label{figperturbeight}
\end{center}
\end{figure}

 This operation changes $K(\petita,\petitb,\petitc)$ to a knot
 $\phi_J\bigl(K(\petita,\petitb,\petitc)\bigr)$ with the same Alexander polynomial as $K(\petita,\petitb,\petitc)$ such that
  $$\wdel\Bigl(\phi_J\bigl(K(\petita,\petitb,\petitc)\bigr)\Bigr) = \wdel\bigl(K(\petita,\petitb,\petitc)\bigr) + 4(\petitb+\petitc) \lambda^{\prime}(J).$$
In the meantime, we have $$w_{SL}\Bigl(\phi_J\bigl(K(\petita,\petitb,\petitc)\bigr)\Bigr) = w_{SL}\bigl(K(\petita,\petitb,\petitc)\bigr)$$ since $\phi_J(\granda)$ still bounds a disk disjoint from $\phi_J(\grandb)$.
  In particular, when $4(b+c) \lambda^{\prime}(J) \neq 0$, which is easily realizable, $\wdel$ distinguishes $\phi_J\bigl(K(\petita,\petitb,\petitc)\bigr)$ from $K(\petita,\petitb,\petitc)$ while $w_{SL}$ does not. (In the case of Figure~\ref{figperturbeight}, $(\petitb+\petitc)=2$ and $\lambda^{\prime}(J)=-1$.) Thus $\wdel$ and $w_{SL}$ are independent knot invariants.
\end{examples}

\begin{remark}
To compute $\lambda^{\prime}\bigl(\phi(\granda),\phi(\grandb)\bigr)$, it suffices to
\begin{itemize}
\item choose a surface $\Sigma_{\granda}$ bounded by $\phi(\granda)$ in the exterior of $\phi(\grandb)$,\footnote{The \emph{exterior} of a knot $\Jknot$ in $\rats$ is the complement of an open tubular neighborhood of $\Jknot$ in $\rats$.}
\item pick a curve $c$ of $H_1(\Sigma_{\granda})$, such that, $\langle x,c \rangle_{\Sigma_{\granda}}=lk\bigl(x,\phi(\grandb)\bigr)$ for any curve $x$ of $\Sigma_{\granda}$, and
\item compute $\lambda^{\prime}\bigl(\phi(\granda),\phi(\grandb)\bigr)=lk(c,c_{\parallel})$, where $c_{\parallel}$ is a curve parallel to $c$ on $\Sigma_{\granda}$.
\end{itemize}
Fix the curve $\phi(\granda)$ of $\partial \ExtH$. The kernel of the map induced by the inclusion from $H_1(\partial \ExtH)$ to $H_1(\ExtH)$ is $\ZZ[\phi(\granda)] \oplus \ZZ[\phi(\grandb)]$. So, the homology class of $\phi(\grandb)$ in $\partial \ExtH \setminus \phi(\granda)$ is well determined up to multiplication by $\pm 1$ and addition of curves parallel to $\phi(\granda)$.
In particular, $\lambda^{\prime}(\phi(\granda),\phi(\grandb))$ is an invariant of the curve $\phi(\granda)$ of $\partial \ExtH$. (It is similarly an invariant of the curve $\phi(\grandb)$.) The homology class of $\phi(\granda)$ on $\partial \ExtH$ is characterized by the properties that $\phi(\granda)$ bounds in $\ExtH$ and that $\phi(\granda)$ is homologous to $\phi(\alpha)$ in $\phi(\hbz)$. Performing a Dehn twist about $\alpha$ on $\Sigma(\petita,\petitb,\petitc)$
does not change the homology class of $\phi(\granda)$ on $\partial \ExtH$.
But it changes $\wdel(\petita,\petitb,\petitc)$ as soon as $\petitb(\petitb^2-1) \neq 0$, as Lemma~\ref{lemvarwdelDehn} shows. Therefore, $\lambda^{\prime}\bigl(\phi(\granda),\phi(\grandb)\bigr)$ is not an invariant of the homology class of $\phi(\granda)$
in $H_1(\partial \ExtH)$. Theorem~\ref{thmSatoinvariant} can be used to reduce the computation of its variation to the computation of the variation of $\wdel(\petita,\petitb,\petitc)$.
\end{remark}

Theorem~\ref{thmSatoinvariant} is a direct consequence
of Theorem~\ref{thmsimpleinvt} and of the following proposition \cite[Proposition 7.3]{lesformagt}.

\begin{proposition}
\label{propgenusone}
There exists an invariant $w_3$ of knots in rational homology $3$-spheres that satisfies the following property.
For any embedding $\phi$ of $\hbz$ in a rational homology sphere $\rats$ such that
$\phi(\granda)$ and $\phi(\grandb)$ are null-homologous in the exterior $\rats \setminus (\phi\bigl(\mhbz\bigr))$ of $\phi(H)$, we have
$$w_3\Bigl(\phi(K(\petita,\petitb,\petitc))\Bigr)=\frac{3}{2}\lambda^{\prime}(\phi(\granda),
\phi(\grandb))-\frac14 \wdel(\petita,\petitb,\petitc) - \frac{\petita}{2} \lambda^{\prime}(\phi(\granda)) - 
\frac{\petitb}{2}  
\lambda^{\prime}(\phi(\grandb))
- \frac{\petitc}{2} \lambda^{\prime}(\phi(\grandc)).$$ 
\end{proposition}

For an integer $q$, let $\rats(K;1/q)$ denote the manifold obtained by $1/q$-surgery on a null-homologous knot $K$ in a $\QQ$-sphere $\rats$.
In \cite[Theorem 7.1]{lesformagt}, I proved a surgery formula 
$$\lambda_2(\rats(K;1/q))=\lambda^{\prime\prime}_2(K)q^2 +w_3(K)q+\lambda_2(\rats)$$
for an invariant $\lambda_2$ of rational homology $3$-spheres.
The invariant $w_3$ that governs the linear term of this formula satisfies the property of Proposition~\ref{propgenusone}, according to \cite[Proposition 7.3]{lesformagt}. It also satisfies
$w_3\bigl(\Kknot \subset (-\rats)\bigr)=-w_3\bigl(\Kknot \subset \rats\bigr).$
When the ambient rational homology $3$-sphere $\rats$ is $S^3$, the invariant
$w_3$ is the degree $3$ knot invariant that changes sign under mirror image and maps the chord diagram with three diameters to $(-1)$. It is a coefficient of the Jones polynomial in a suitable normalization. 

For knots that bound null-homologous Seifert surfaces, Theorem~\ref{thmsimpleinvt} and Theorem~\ref{thmSatoinvariant} stated above can be summarized as follows.

\begin{theorem} \label{thmw3decomp} If $K$ is a knot in a 
$\QQ$-sphere $\rats$ such that $K$ bounds a null-homologous Seifert surface in $\rats$, then the invariant $w_3(K)$ is the sum of the following independent knot invariants described in Theorems~\ref{thmsimpleinvt} and \ref{thmSatoinvariant}.
 $$w_3(K) = \frac{3}{2}w_{SL}(K) - \frac{1}{8} \wdel(K).$$ 
\end{theorem}

This article and the discovery of Theorem~\ref{thmsimpleinvt} grew up from an attempt to understand surgery formulas for finite type invariants of $\QQ$-spheres and finite type invariants of knots, such as $w_3$, from finite type invariants of curves on Seifert surfaces. (For a null-homologous knot $K$ in a $\QQ$-sphere $\rats$, the \say{degree one} Casson invariant $\lambda_C$ satisfies the surgery formula $\lambda_C(\rats(K;1/q))=\lambda^{\prime}(K)q+\lambda_C(\rats)$, and the knot invariant $\lambda^{\prime\prime}_2(K)$ of the quadratic term of the surgery formula above for the \say{degree two} invariant $\lambda_2$ is a polynomial in the entries of a Seifert matrix. See \cite[Theorem 7.1]{lesformagt}.)

\subsection{Definition of Alexander forms and Reidemeister torsions}
\label{subdefAlef}

We will prove Theorems~\ref{thmw3decomp} and \ref{thmsimpleinvt} by studying a Reidemeister torsion associated with a genus one Seifert surface.
In this subsection, we define the needed Reidemeister torsions in a way appropriate for our study from the Alexander forms described below.

Alexander forms were introduced in \cite[Section 3]{lesinv}.\footnote{Alexander forms were called Alexander functions in \cite[Section 3]{lesinv}.} They help study Alexander polynomials with several variables, and the coefficients $\zeta$ of \cite{lespup}, in particular.
Below, we recall their definition and some of their properties. We develop more examples and more properties in Section~\ref{secAlef}.

Throughout this section, $\ambx$ denotes a connected compact oriented 3-manifold 
with a nonempty boundary, equipped with a basepoint $\star$.
Define the \emph{genus} $g(\ambx)$ of $\ambx$ to be $$ g(\ambx) = 1 - \chi(\ambx) \; \bigl(= 1 -\frac12 \chi(\partial \ambx)\bigr).$$
Set $g=g(\ambx)$ and assume $g >0$.

Let $\tilde{\ambx}$ denote the maximal free abelian covering space of $\ambx$, and let $p_{\ambx}\colon \tilde{\ambx}\to \ambx$ be the associated
covering map.
Define the \emph{group ring} $\Lambda_{\ambx}$ of the covering group ${H_1(\ambx)}/{\mbox{\scriptsize Torsion}}$ of $\tilde{\ambx}$ to be the $\ZZ$-module
$$\Lambda_{\ambx} = \ZZ\left[\frac{H_1(\ambx)}{{\rm Torsion}(H_1(\ambx))}\right] = \bigoplus_{x \in \frac{H_1(\ambx)}{\rm Torsion}} \ZZ\, \exp(x), $$ equipped with the
 $\ZZ$-bilinear multiplication law that maps $\bigl(\exp(x),\exp(y)\bigr)$ to $\exp(x+y)$ for any $x, y \in {H_1(\ambx)}/{\rm Torsion}$.
We use the notation $\exp(x)$ or $e^{x}$ to denote an element $x$ of ${H_1(\ambx)}/{\rm Torsion}$ viewed as a generator of $\Lambda_{\ambx}$ to remind
this multiplication law.
Similarly define the group rings $\Lambda_{\ambx}^{\QQ}=\ZZ\left[H_1(\ambx;\QQ)\right]$ and $\Lambda_{\ambx}^{\frac12}=\ZZ\bigl[\frac12{H_1(\ambx)}/{{\rm Torsion}(H_1(\ambx))}\bigr]$, in which $\Lambda_{\ambx}$ embeds. These group rings are equipped with the linear involution $$\sum_x\lambda_x \exp(x) \mapsto \left( \overline{\sum_x\lambda_x\exp(x)}=\sum_x\lambda_x\exp(-x)\right).$$

The units of $\Lambda_{\ambx}$, are its elements of the form $\exp\left(x \in {H_1(\ambx)}/{\mbox{Torsion}}\right)$ called the \emph{positive units},
and its elements of the form $ - \exp\left(x \in {H_1(\ambx)}/{\mbox{Torsion}}\right)$.

The Alexander form of $\ambx$ is the following invariant of the $\Lambda_{\ambx}$-module
$$ \CH_{\ambx}=H_1\bigl(\tilde{\ambx}, p_{\ambx}^{-1}(\star)\bigr).$$
\begin{definition} 
\label{defalexform}
Let $\CP$ be a presentation of $\CH_{\ambx}$ over $\Lambda_{\ambx}$ with $(r+g)$ generators $\gamma_1,\dots,\gamma_{r+g}$, and $r$ relators $\rho_1, \dots,\rho_r$, which are $\Lambda_{\ambx}$-linear combinations of the $\gamma_i$.\footnote{Such a presentation exists. It may be obtained from a cellular decomposition of $\ambx$ with one $0$-cell, $(r+g)$ $1$-cells, $r$ $2$-cells and no other cells.} Define the element $\hat{\gamma}$ of $\bigwedge^{r+g} \left( \bigoplus_{i=1}^{r+g}\Lambda_{\ambx} \gamma_i \right) $ to be $\hat{\gamma} = \gamma_1 \wedge \dots \wedge \gamma_{r+g}$. Similarly define $\hat{\rho}=\rho_1 \wedge \dots \wedge \rho_r \in \bigwedge^{r} \left( \bigoplus_{i=1}^{r+g}\Lambda_{\ambx} \gamma_i \right)$.
The \emph{Alexander  form} $\CA_{\ambx,\CP}$ of $(\ambx,\CP)$ is the $\Lambda_{\ambx}$-linear form 
$$ \CA_{\ambx,\CP}\colon \bigwedge^{g} \CH_{\ambx} \longrightarrow \Lambda_{\ambx}  $$ that maps an element
$\hat{u} = u_1 \wedge \dots \wedge u_g$ of $\bigwedge^g \CH_{\ambx}$ to the
element $\CA_{\ambx,\CP}(\hat{u})$ of $\Lambda_{\ambx}$ defined by the equality
\begin{equation}
\label{edaa}
\hat{\rho} \wedge \hat{u} = \CA_{\ambx,\CP}(\hat{u}) \hat{\gamma} 
\end{equation}
where the $u_i$ are represented as combinations of the $\gamma_j$, and the exterior products in Equation~\ref{edaa} belong to
$\bigwedge^{r+g} \left( \bigoplus_{i=1}^{r+g}\Lambda_{\ambx} \gamma_i \right) $. 

As proved in \cite[Section 3.1]{lesinv}, if $\CP$ and $\CQ$ are two presentations of $\CH_{\ambx}$ as above, then there exists a unit $\pm \exp(u)$ of $\Lambda_{\ambx}$ such that $\CA_{\ambx,\CQ}=\pm \exp(u)\CA_{\ambx,\CP}$. We will often fix $\CP$, omit $\CP$ from the notation, keep in mind that $\CA_{\ambx} = \CA_{\ambx,\CP}$ is defined up to a multiplication by a unit of $\Lambda_{\ambx}$, and call $\CA_{\ambx}$ the \emph{Alexander form} of $\ambx$.
\end{definition}

The Alexander form $\CA$ is a topological invariant of the pairs $(\ambx,\star)$ up to homotopy equivalence of pairs (connected compact oriented 3-manifold with boundary and with non-negative genus, basepoint). Homotopy equivalences
provide identifications between the involved $\Lambda$s and $\CH$s.
Changing the location of the basepoint $\star$ of $\ambx$ does not affect the homotopy equivalence class of $(\ambx,\star)$.
Thus $\CA$ is an invariant of $\ambx$. But we need a reference basepoint, which we may choose anywhere we want in $\ambx$.
So, we fix a basepoint $\star$ of $\ambx$ and a preferred lift $\star_0$ of $\star$ in $\tilde{\ambx}$.

When $c$ is an oriented curve in $\ambx$, $c$ also denotes the homology class
it carries and the associated variable in $\exp(c) \in \Lambda_{\ambx}$. If furthermore $c$ is based at $\star$, then
$c$ also denotes its own class in $\pi_1(\ambx,\star)$. We denote the class in $\CH_{\ambx}$ of the preferred lift of $c$ starting at 
$\star_0$ in $\tilde{\ambx}$ by $\tilde{c}$.

\begin{example} \label{exaAlextriv}
Let $\ExtH_0$ be the closure of the complement of the handlebody $\hbz$ of Figure~\ref{fighandlebobasedu} trivially embedded in $S^3$, as the figure suggests. We equip $\partial \hbz$ with the based curves $\ua$, $\ub$, and $\uc$ of the figure. Their product $\ua\ub\uc$ is trivial in $\pi_1(\partial \hbz,\star)$.

\begin{figure}[h]
\begin{center}
\begin{tikzpicture}

 \begin{scope}[xscale=4]

\draw [thick, rounded corners, fill=blue!20] (3.4,2) -- (3.4,3.8) -- (0,3.8) -- (0,2)
(0,2) -- (.6,2) 
(.6,2) -- (.6,3) -- (.8,3) -- (1.4,2.4) -- (1.4,2)
(1.4,2) -- (2,2)
(2,2) -- (2,2.4) -- (2.6,3) -- (2.8,3) -- (2.8,2)
(2.8,2) -- (3.4,2); 
\draw [rounded corners, draw=blue!20, fill=blue!20]
(0,2.2) -- (0,0) -- (3.4,0) -- (3.4,2.2) (3.4,2.2) -- (2.8,2.2)
(2.8,2.2) -- (2.8,.8) -- (2.6,.8) -- (2,1.4) -- (2,2.2) (2,2.2) -- (1.4,2.2)
(1.4,2.2) -- (1.4,1.4)-- (.8,.8) -- (.6,.8) -- (.6,2.2) (.6,2.2) -- (0,2.2);
\draw [rounded corners,thick]
(0,1.6) -- (0,0) -- (3.4,0) -- (3.4,1.6)
(2.8,1.6) -- (2.8,.8) -- (2.6,.8) -- (2,1.4) -- (2,1.6)
(1.4,1.6) -- (1.4,1.4) -- (.8,.8) -- (.6,.8) -- (.6,1.8);
\begin{scope}[xshift=-1.3cm,yshift=-.1cm] 
\draw [-] (2,.9) .. controls (1.96,.9) and (1.92,.7) .. (1.92,.5);
\draw (1.92,.5) .. controls (1.92,.3) and (1.96,.1) .. (2,.1);
\draw [dashed] (2,.9) .. controls (2.04,.9) and (2.08,.7) .. (2.08,.5) .. controls (2.08,.3) and (2.04,.1) .. (2,.1);
\end{scope}
\begin{scope}[xshift=.7cm,yshift=-.1cm]
\draw [-] (2,.9) .. controls (1.96,.9) and (1.92,.7) .. (1.92,.5); 
\draw (1.92,.5) .. controls (1.92,.3) and (1.96,.1) .. (2,.1);
\draw [dashed] (2,.9) .. controls (2.04,.9) and (2.08,.7) .. (2.08,.5) .. controls (2.08,.3) and (2.04,.1) .. (2,.1);
\end{scope}
 \draw [-<] (0,2.2) .. controls (0,2.12) and (.15,2) ..  (.3,2);
 \draw (.24,1.92) node[above]{\scriptsize $\ub$} (.3,2)  .. controls (.45,2)  and  (.6,2.12)  .. (.6,2.2);
\draw [dashed,->] (0,2.2) .. controls (0,2.28) and (.15,2.4) ..  (.3,2.4);
\draw [dashed] (.3,2.4)  .. controls (.45,2.4)  and  (.6,2.28)  .. (.6,2.2);
\draw (0,2.2) -- (0,1.6) (.6,1.6) -- (.6,2.2);
\begin{scope}[xshift=1.4cm]
 \draw [-<] (0,2.2) .. controls (0,2.12) and (.15,2) ..  (.3,2);
 \draw (.24,1.92) node[above]{\scriptsize $\uc$} (.3,2)  .. controls (.45,2)  and  (.6,2.12)  .. (.6,2.2);
\draw [dashed,->] (0,2.2) .. controls (0,2.28) and (.15,2.4) ..  (.3,2.4);
\draw [dashed] (.3,2.4)  .. controls (.45,2.4)  and  (.6,2.28)  .. (.6,2.2);
\draw (0,2.2) -- (0,1.6) (.6,1.6) -- (.6,2.2);
\end{scope}
\begin{scope}[xshift=2.8cm]
 \draw [-<] (0,2.2) .. controls (0,2.12) and (.15,2) ..  (.3,2);
 \draw (.36,1.92) node[above]{\scriptsize $\ua$} (.3,2)  .. controls (.45,2)  and  (.6,2.12)  .. (.6,2.2);
\draw [dashed,->] (0,2.2) .. controls (0,2.28) and (.15,2.4) ..  (.3,2.4);
\draw [dashed] (.3,2.4)  .. controls (.45,2.4)  and  (.6,2.28)  .. (.6,2.2);
\draw (0,2.2) -- (0,1.6) (.6,1.6) -- (.6,2.2);
\end{scope}
\draw [rounded corners,->]  (2.2,3.65) -- (.15,3.65) -- (.15,.15) -- (3.25,.15) -- (3.25,3.65) -- (2.2,3.65);
\draw [rounded corners,->]  (1.15,2.95) -- (1.55,2.55) -- (1.55,1.25) -- (.95,.65) -- (.45,.65) -- (.45,3.15) -- (.95,3.15) -- (1.15,2.95);
\draw [rounded corners,->]  (2.25,2.95) -- (2.45,3.15) -- (2.95,3.15) -- (2.95,.65) -- (2.45,.65) -- (1.85,1.25) -- (1.85,2.55) -- (2.25,2.95);
\draw (2.3,3.72) node[below]{\scriptsize $\grandc$};
\draw (1.05,3.15) node[right]{\scriptsize $\granda$};
\draw (2.25,3.05) node[left]{\scriptsize $\grandb$};
\draw (1.7,3.65) -- (1.7,3.45) -- (2.05,2.75) (1.7,3.45) -- (1.35,2.75);
\draw (1.7,3.45) -- (1.7,2);
\draw [rounded corners] (.3,2) -- (.4,3.25) -- (1.7,3.45) -- (3,3.25) -- (3.1,2);
\end{scope}

\end{tikzpicture}
\caption{The based curves $\ua$, $\ub$, $\uc$, $\granda$, $\grandb$, and $\grandc$ on $\partial \hbz$}
\label{fighandlebobasedu}
\end{center}
\end{figure}

Since $H_1(\ExtH_0)=\ZZ \ua \oplus \ZZ \ub$,
we have $\Lambda_{\ExtH_0}=\ZZ[e^{\ua}, e^{-\ua},e^{\ub}, e^{-\ub}].$
We can also use the following symmetric expression
$$\Lambda_{\ExtH_0}=\frac{\ZZ[e^{\ua}, e^{-\ua},e^{\ub}, e^{-\ub},e^{\uc}, e^{-\uc}]}{(e^{\ua}e^{\ub}e^{\uc}=1)}.$$
Since $\ExtH_0$ retracts on the wedge of the two circles $\ua$ and $\ub$,
$\tilde{\ExtH_0}$ is diffeomorphic to a regular neighborhood in $\RR^3$
of the grid $\bigl(\RR \times \ZZ \times\{0\}\bigr) \cup \bigl(\ZZ \times \RR \times\{0\}\bigr)$ of $\RR^3$ in Figure~\ref{figabeliancover}. In this grid, the preimage of $\star$ is $\ZZ \times \ZZ \times\{0\}$.

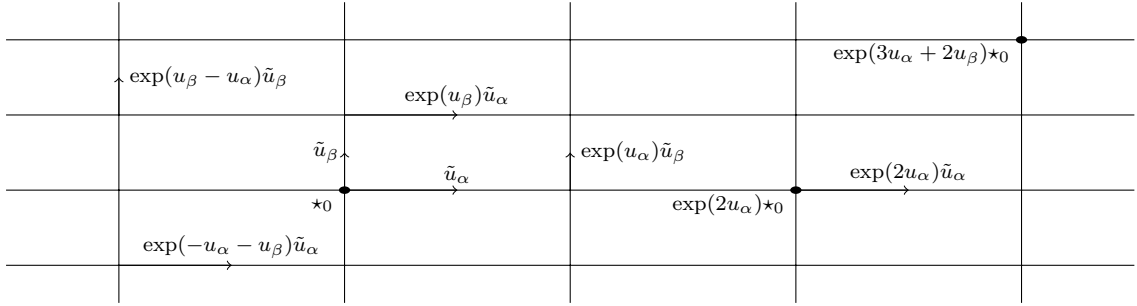
\begin{figure}[h]
\begin{center}
\begin{tikzpicture}[xscale=1.5] 
\draw (-3,-1) -- (7,-1) (-3,0) -- (7,0) (-3,1) -- (7,1) (-3,2) -- (7,2)
(-2,-1.5) -- (-2,2.5) (0,-1.5) -- (0,2.5) (2,-1.5) -- (2,2.5) (4,-1.5) -- (4,2.5) (6,-1.5) -- (6,2.5) (0,-.2) node[left]{\scriptsize $\star_0$}   (6,1.8) node[left]{\scriptsize $\exp(3\ua+2\ub)\star_0$} (4,-.2) node[left]{\scriptsize $\exp(2\ua)\star_0$};
\draw[->] (1,-.05) node[above]{\scriptsize $\tua$} (0,0) -- (1,0);
\draw[->] (5,-.05) node[above]{\scriptsize $\exp(2\ua)\tua$} (4,0) -- (5,0);
\draw[->] (1,.95) node[above]{\scriptsize $\exp(\ub)\tua$} (0,1) -- (1,1);
\fill (0,0) circle (.05) (4,0) circle (.05) (6,2) circle (.05);
\draw[->]  (0.05,.5) node[left]{\scriptsize $\tub$} (0,0) -- (0,.5);
\draw[->]  (2,.5) node[right]{\scriptsize $\exp(\ua)\tub$} (2,0) -- (2,.5);
\draw[->] (-1,-1.05) node[above]{\scriptsize $\exp(-\ua-\ub)\tua$} (-2,-1) -- (-1,-1);
\draw[->]  (-2,1.5) node[right]{\scriptsize $\exp(\ub-\ua)\tub$} (-2,1) -- (-2,1.5);
\end{tikzpicture}
\caption{The one-skeleton of $\tilde{\ExtH_0}$}
\label{figabeliancover}
\end{center}
\end{figure}
In particular, we have
$\CH_{\ExtH_0}=\Lambda_{\ExtH_0}\tua \oplus \Lambda_{\ExtH_0}\tub.$
The associated presentation $\CP_0$ of $\CH_{\ExtH_0}$ with two generators and no relator provides the Alexander form 
$\CA_{\ExtH_0,\CP_0}$ such that $\CA_{\ExtH_0,\CP_0}(\tua \wedge \tub)=1$. Since 
$\granda$ bounds a disk in $\ExtH_0$, the linear form $\CA_{\ExtH_0,\CP_0}(\tilde{\granda} \wedge .)$ vanishes.

Note that the following \emph{Fox calculus rules} apply in $\CH_{\ExtH}$.
For any two based curves $\gamma$ and $\cdel$ of $\ambe$, we have
$\widetilde{\gamma\cdel}=\tilde{\gamma}+\exp(\gamma)\tilde{\cdel}$ and
$\widetilde{\gamma^{-1}}=-\exp(-\gamma)\tilde{\gamma}$. See Figure~\ref{figabeliancover}.

As a consequence, we get the following lemma.
\begin{lemma} \label{lemuaubuc}
 Let $\ambe$ be a compact $3$-manifold whose boundary is the genus $2$ surface $(-\partial \hbz)$ of Figure~\ref{fighandlebobasedu}. The classes of the based curves $\ua$, $\ub$, and $\uc$ of Figure~\ref{fighandlebobasedu}
 satisfy the following relations.
$$\ua\ub\uc=1\;\;\mbox{in}\;\; \pi_1(\ambe),$$
$$\tua + \exp(\ua) \tub + \exp (-\uc)\tuc=0\;\;\mbox{in}\;\;\CH_{\ambe},$$
$$\tub \wedge \tuc=\exp(\uc)\tua \wedge \tub\;\;\mbox{and}\;\;\tuc \wedge\tua=\exp(-\ub) \tua \wedge \tub \;\;\mbox{in}\;\; \CH_{\ambe} \wedge \CH_{\ambe}.$$
\end{lemma}
\eop

In particular, back to Example~\ref{exaAlextriv}, we have
$\CA_{\ExtH_0,\CP_0}(\tub \wedge \tuc)=\exp(\uc)$ and $\CA_{\ExtH_0,\CP_0}(\tuc \wedge \tua)=\exp(-\ub)$.
\end{example}

\begin{notation}
The choice of the preferred lift $\star_0$ of $\star$ in $\tilde{\ambx}$ provides a natural isomorphism from $H_0(p_{\ambx}^{-1}(\star)) = \Lambda_{\ambx}[\star_0]$ to $\Lambda_{\ambx}$. The boundary map
from $\CH_{\ambx}$ to $H_0(p_{\ambx}^{-1}(\star))$ composed with this isomorphism is denoted by $\partial$. The map $\partial$ maps the class $\tilde{x}$ of the lift of a based curve $x$ that starts at $\star_0$ to $\partial(\tilde{x})=\exp(x) -1$.
\end{notation}

When $\ambx$ is a link exterior, its genus $g(\ambx)$ is $1$. In this case, the \emph{Reidemeister torsion} of $\ambx$ is the ratio of the following proposition, which is \cite[Property 1, p. 640]{lesinv}. We use this proposition as a definition for the Reidemeister torsion.

The sign $\doteq$ stands for \say{equal up to multiplication by a positive unit (of $\Lambda_{\ambx}$ or $\Lambda_{\ambx}^{\frac12}$)}.

\begin{proposition}
\label{proprtort} Let $\ambx$ be a link exterior. Let $\tau(\ambx)$ be its Reidemeister torsion described in \cite{turaevReid}. For a given normalization $\CA_{\ambx,\CP}$ of $\CA_{\ambx}$, the linear form $\CA_{\ambx,\CP}$ is proportional to $\partial$.
For any element $u$ of $\CH_{\ambx}$, we have
$$\CA_{\ambx}(u) \doteq \pm \tau(\ambx)\partial(u).$$
\end{proposition}

Let $c$ be a simple closed curve in the boundary $\partial \ambx$ of a $3$-manifold $\ambx$. Let $[0,1] \times c$ be a tubular neighborhood of $c$ in $\partial \ambx$.
Then the manifold obtained by \emph{attaching a $2$-handle} along $c$
is obtained by gluing a solid cylinder $[0,1] \times D^2$ along $[0,1] \times c$ so that $[0,1] \times \partial D^2$ is identified to $[0,1] \times c$, and by smoothing the angles. We denote it by $\ambx[c]$.

Alexander forms are useful to study Reidemeister torsions since $\CA_{\ambx}$ recovers the Reidemeister torsions of the link complements obtained from $\ambx$ by attaching $2$-handles along $\partial \ambx$ as above.
For example, if $c$ is based and if $\ambx[c]$ is a link complement, then the genus of $\ambx$ is $2$ and, up to multiplication by a unit, $\CA_{\ambx[c]}(.)$ is obtained from $\pm\CA_{\ambx}(\tilde{c}\wedge .)$ by composition by the quotient map $\Lambda_{\ambx} \to \Lambda_{\ambx[c]}$ induced by the inclusion.

\subsection{Invariants of genus one surfaces in rational homology spheres}
\label{subintrocd}

Assume that the manifold $\ambx$ of the previous subsection is the exterior of a genus one surface $\Sigma=\phi\bigl(\Sigma(\petita,\petitb,\petitc)\bigr)$ for an embedding 
$\phi \colon \hbz \hookrightarrow \rats$ of $\hbz$ into a $\QQ$-sphere $\rats$.
Let $\ExtH[K]$ be the $3$-manifold obtained from this exterior $\ExtH= \rats \setminus \phi\bigl(\mhbz\bigr)$ by attaching a $2$-handle along $(\partial \Sigma = K)$. 
Then $\ExtH$, $\ExtH[K]$, and the Reidemeister torsion 
of $\ExtH[K]$ are invariants of $\Sigma$.
We use the Alexander form of $\ExtH$ to normalize the Reidemeister torsion of $\ExtH[K]$ as follows.

\begin{notation}
Let $\varepsilon$ denote the $\ZZ$-linear \emph{augmentation morphism}
$$\begin{array}{llll}\varepsilon \colon & \Lambda_{\ambe} &\longrightarrow &\ZZ\\
&\exp\left(x \in {H_1(\ambe)}/{\mbox{Torsion}}\right)& \mapsto & 1. \end{array}$$
\end{notation}

An \emph{orientation} of a rational vector space is an equivalence class of bases under the equivalence relation that identifies two bases related by a linear transformation with a positive determinant.
For example, symplectic bases orient $H_1(\Sigma;\QQ)$. The natural isomorphism
$$\begin{array}{llll}\fHone \colon &H_1(\Sigma;\QQ) &\rightarrow &H_1(\ExtH;\QQ)\\
&z &\mapsto& z-z_-
\end{array}$$
carries the orientation of $H_1(\Sigma;\QQ)$ to an orientation of  $H_1(\ExtH;\QQ)$ induced by the orientation of $\Sigma$.

Let $(\alpha,\beta)$ be a symplectic basis of $\Sigma$. View $\partial \ExtH$ as $\Sigma \cup_K(-\Sigma_-)$.
Assume that the basepoint of $\ExtH$ is on $K \subset \partial \ExtH$ and that $\alpha$ and $\beta$ are attached to this basepoint of $\ExtH$. So $\alpha$ and $\beta$ are based curves, and their copies $\alpham$ and $\betam$ in $\Sigma_-$ are similarly based. The basis $(\alpha-\alpham,\beta-\betam)$ of $H_1(\ExtH;\QQ)$ defines the above orientation of $H_1(\ExtH;\QQ)$ and allows us to fix the sign of $\CA_{\ExtH}$ so that $\varepsilon\bigl(\CA_{\ExtH}(\talpha-\talpham\wedge\tbeta-\tbetam)\bigr)=|H_1(\rats)|$. So $\CA_{\ExtH}$ is defined up to multiplication by a positive unit.

Recall that $K$ is oriented as the boundary of $\Sigma$.
We fix the sign of the Reidemeister torsion of $\ExtH[K]$ so that 
$$\CA_{\ExtH}(\tilde{K}\wedge u) \doteq \tau(\ExtH[K])\partial(u)$$
for any element $u$ of $\CH_{\ExtH}=\CH_{\ExtH[K]}$.
Finally, thanks to the well-known symmetry \cite{turaevReid} of the Reidemeister torsion up to units, we can multiply a representative of $\tau(\ExtH[K])$ by a unit of $\Lambda_{\ExtH}^{\frac12}$ to obtain a well-defined element $\CD(\Sigma)$ of $\Lambda_{\ExtH}^{\frac12}$
 such that $\overline{\CD(\Sigma)}= \CD(\Sigma)$.

Note that we used twice the orientation of $\Sigma$ to determine the sign of 
$\CD(\Sigma)$. So, we have $\CD(\Sigma)=\CD(-\Sigma)$. We have proved the following proposition.

\begin{proposition}\label{propfHone}
For any oriented genus one surface $\Sigma$ in a rational homology sphere,
the Reidemeister torsion $\CD(\Sigma)$ normalized as above is an invariant of (the isotopy class of the image of the embedding of) $\Sigma$, with values in the group ring $\Lambda_{\ExtH}^{\frac12}$. 

Let
$\fHone_{\ast} \colon \ZZ[H_1(\Sigma;\QQ)] \to \Lambda_{\ExtH[K]}^{\QQ}$ be the isomorphism induced by the isomorphism $\fHone \colon H_1(\Sigma;\QQ) \rightarrow H_1(\ExtH;\QQ)$ such that $f(z)=z-z_-$.
Then $\fHone_{\ast}^{-1}(\CD(\Sigma))$ is an invariant of $\Sigma$ valued in $\ZZ[H_1(\Sigma;\QQ)]$.
\end{proposition}

We are going to extract numerical invariants from the invariant $\CD(\Sigma)$ in Theorem~\ref{thmalexcanphi} and Proposition~\ref{propsergenone}.

\begin{definition}
 Let $V$ be a rational vector space. 
 The \emph{symmetric algebra} $S(V)$ of $V$ is the quotient of the tensor algebra of $V$ by the two-sided ideal generated by the elements of the form $u \otimes v - v \otimes u$. 
It is a graded algebra $S(V)=\oplus_{n \in \NN} S_n(V)$ whose degree $n$ part $S_n(V)$ is the \emph{symmetrized $n^{\mbox{\scriptsize th}}$ tensor power} $\otimes_s^n V$ of $V$. We set $\QQ[[V]]=\prod_{n\in \NN}S_n(V)$.
\end{definition}

In particular, $S_2(V)$ is the symmetrized tensor product $V \otimes_s V$. (It is the quotient of $ V \otimes V$ by the vector space generated by the differences $(u\otimes v-v\otimes u)$ for pairs $(u,v) \in V^2$.) We often denote products $u \otimes_s v$ of $V \otimes_s V$ by $uv$, and forget the signs $\otimes_s$ in $S_n(V)$ in the same way.

A morphism $\phi \colon V \to W$ of rational vector spaces induces a natural algebra morphism $\phi_{\ast} \colon \QQ[[V]] \to \QQ[[W]]$. When $V=\bigoplus_{i=1}^k\QQ v_i$, the algebra $\QQ[[V]]$ is isomorphic to the algebra $\QQ[[v_1,\dots, v_k]]$ of formal power series in the commuting variables $v_i$, for $i=1, \dots, k$. 
Note the following lemma.

\begin{lemma}
\label{lemexpans} For any finite-dimensional rational vector space $V$, the linear map
$$ \begin{array}{lll}\ZZ[V] & \to & \QQ[[V]]\\
 \exp\left(v\right) & \mapsto & \sum_{n \in \NN}\frac{1}{n !}v^n,
             \end{array}$$
is an injective algebra morphism.
\end{lemma}
\bp It is easy to see that the above linear map is an algebra morphism. When $V=\QQ$, it is injective because the vectors
$(a_1,a_2,\dots, a_n)$, $(a_1^2,a_2^2,\dots, a_n^2)$, \dots, $(a_1^n,a_2^n,\dots, a_n^n)$ are independent in $\QQ^n$ for $n$ pairwise distinct $a_i$.
Let $u_1$, \dots, $u_n$ be pairwise distinct vectors in $V$. Then there exists a linear form $\psi$ on $V$ such that $\psi(u_1)$, \dots, $\psi(u_n)$ are pairwise distinct. (Choose $\psi$ outside the finitely many hyperplanes of $\mbox{Hom}(V;\QQ)$ of linear forms that send $(u_j-u_i)$ to zero for pairs $\{u_i,u_j\}$ of distinct elements.) Then the previous case implies that $\bigl(\exp(u_1), \dots, \exp(u_n)\bigr)$ is a free system in $\QQ[[V]]$.
\eop

Lemma~\ref{lemexpans} allows us to consider $\Lambda_{\ExtL}^{\QQ}$ as a subring of $\QQ[[H_1(\ExtL;\QQ)]]$ by writing $\exp(v)$ as $\sum_{n \in \NN}\frac{1}{n !}v^n$.
The elements of $\Lambda_{\ExtL}^{\QQ}$ that are (infinite) sums of tensors of even (resp. odd) degree in $\QQ[[V]]$ are called \emph{even} (resp. \emph{odd}).

With these conventions, $\CD(\Sigma)$ is an even element of
$\QQ[[H_1(\ExtH;\QQ) ]]$, 
whose degree $0$ part is $|H_1(\rats)|\lambda^{\prime}(K)$, according to Lemma~\ref{lemCDdegzero}.
We determine the whole series $\CD(\Sigma)$ in Section~\ref{subdescDsigma}. See Proposition~\ref{propgenDSigma}.
Its degree $2$ part $\CD_2(\Sigma)$ belongs to $H_1(\ExtH;\QQ) \otimes_s H_1(\ExtH;\QQ)$. 

\begin{notation} \label{notsymSeif} Let $V_s$ denote 
the symmetrized \emph{Seifert form} 
$$\begin{array}{llll}V_s \colon &H_1(\Sigma;\QQ) \otimes_s H_1(\Sigma;\QQ)&\rightarrow 
&\QQ\\
&u \otimes_s v &\mapsto& lk(u,v_-)+lk(u_-,v).
\end{array}$$
\end{notation}
If $\phi$ maps the curves $\granda$ and $\grandb$ of Figure~\ref{fighandlebo} to $0$ in $H_1(\ExtH;\QQ)$, then $V_s(\alpha \otimes_s \alpha)=\petitb+\petitc$, $V_s(\beta \otimes_s \beta)=\petitc+\petita$, and $V_s(\alpha \otimes_s \beta)= -\petitc$.
The form $V_s$ allows us to 
extract the numerical invariant $V_s\left(\fHone^{-1}_{\ast}\CD_2(\Sigma) \right)$ of $\Sigma$.

\begin{theorem}
\label{thmalexcanphi} 
Let $\Sigma$ be a genus one Seifert surface in a $\QQ$-sphere $\rats$.
Set $W_s=W_s(\Sigma)=V_s \circ (\fHone^{-1} \otimes_s \fHone^{-1})$.
Then $V_s\left(\fHone^{-1}_{\ast}\CD_2(\Sigma) \right)=W_s\left(\CD_2(\Sigma) \right)$ is a topological invariant of $\Sigma$.
Let $w_{SL}(\Sigma)$ be the invariant defined in Theorem~\ref{thmSatoinvariant} when $\Sigma$ is null-homologous, and in Theorem~\ref{thmwsl}, in general.
We have
$$W_s\left(\CD_2(\Sigma) \right) = |H_1(\rats)| \left( \left(4\lambda^{\prime}(\partial \Sigma)-1\right)w_{SL}(\Sigma) +\frac14 \lambda^{\prime}(\partial \Sigma)\wdel(\Sigma) \right).$$
\end{theorem}

The stated invariance of $W_s\left(\CD_2(\Sigma)\right)$ is tautological. We compute $W_s\left(\CD_2(\Sigma) \right)$ in Section~\ref{secAlexEKlow}, where Theorem~\ref{thmalexcanphi} is proved.

If $K$ bounds a null-homologous genus one Seifert in a $\QQ$-sphere, then $\lambda^{\prime}(K) \in \ZZ$. In particular, we have $7\lambda^{\prime}(K) \neq 1$. Therefore, Theorem~\ref{thmalexcanphi} and Proposition~\ref{propgenusone} have the following immediate corollary.

\begin{corollary}
\label{corwww}
 For any null-homologous 
 genus one Seifert surface $\Sigma$ in a $\QQ$-sphere $\rats$, we have
 $$w_{SL}(\Sigma)=\frac{1}{7\lambda^{\prime}(\partial \Sigma)-1}\left( 2\lambda^{\prime}(\partial \Sigma)w_3(\partial \Sigma) + \frac{W_s\left(\CD_2(\Sigma) \right)}{|H_1(\rats)|}\right)$$
 and 
 $$\wdel(\Sigma)=\frac{4}{7\lambda^{\prime}(\partial \Sigma)-1}\left( \left(2-8\lambda^{\prime}(\partial \Sigma)\right)w_3(\partial \Sigma) + 3 \frac{W_s\left(\CD_2(\Sigma) \right)}{|H_1(\rats)|}\right).$$
\end{corollary}

\subsection{Sketch of the proof of the invariance \texorpdfstring{of $\wdel$}{}}
\label{subintrosk}

According to Corollary~\ref{corwww},
to prove Theorem~\ref{thmsimpleinvt} for null-homologous genus one Seifert surfaces, it suffices to prove that 
$W_s\left(\CD_2(\Sigma) \right)$ does not depend on the surface $\Sigma$ with boundary $K$, as stated in Theorem~\ref{thmcorwratgen} below. 
In Sections~\ref{secinvcob} and \ref{secKcob}, we will prove the following two lemmas, which imply Theorem~\ref{thmcorwratgen}.

Say that two oriented surfaces $\Sigma$ and $\Sigma^{\prime}$ with boundary $K$
in a $\QQ$-sphere $\rats$ are \emph{$K$-cobordant} if there exists a diffeomorphism $\psi$ of $\rats$ isotopic to the identity of $\rats$ such that $\psi(K)=K$ and $\psi(\Sigma) \cap \Sigma^{\prime}=K$. Since the latter condition implies $\Sigma \cap \psi^{-1}(\Sigma^{\prime})=K$, this definition is symmetric with respect to the exchange of $\Sigma$ and $\Sigma^{\prime}$.

\begin{lemma}
\label{leminvKcob}
Let $\Sigma$ and $\Sigma^{\prime}$ be two genus one Seifert surfaces of a knot $K$ in a rational homology $3$-sphere $\rats$. If $\Sigma$ and $\Sigma^{\prime}$ are $K$-cobordant, then 
$$W_s(\Sigma^{\prime})\left(\CD_2(\Sigma^{\prime}) \right)= W_s(\Sigma)\left(\CD_2(\Sigma) \right).$$
\end{lemma}

\begin{lemma}
\label{lemcob}
Let $\rats$ be a 
$\QQ$-sphere. For any two genus one Seifert surfaces $\Sigma$ and $\Sigma^{\prime}$ of $K$ in $\rats$, there exists a sequence $(\Sigma_i)_{i=1, \dots, k}$ of genus one Seifert surfaces of $K$ such that $\Sigma=\Sigma_1$, $\Sigma^{\prime}=\Sigma_k$, and, for any $i=1, \dots, k-1$, $\Sigma_i$ and $\Sigma_{i+1}$ are $K$-cobordant.
\end{lemma}

\begin{theorem}
\label{thmcorwratgen}
Let $K$ be a knot in a
 $\QQ$-sphere. If $K$ bounds a genus one Seifert surface $\Sigma$, then the invariant $W_s\left(\CD_2(\Sigma) \right)$ of Theorem~\ref{thmalexcanphi} does not depend on such a genus one Seifert surface $\Sigma$. It is an invariant of the knot $K$.
\end{theorem}
\eop

\subsection{Alexander forms of genus two homology handlebodies}

Our proof of Theorem~\ref{thmalexcanphi} relies on a study of Alexander forms of genus two rational homology handlebodies. This study produces the following structure theorem (\ref{thmAlexHg}), which is of independent interest. 
Theorem~\ref{thmAlexHg} involves the normalized multivariable Alexander polynomial $\Delta$ of two-component links. The polynomial of such a link is obtained from the Reidemeister torsion of its exterior by a multiplication by a unit and a change of variables. We can use Theorem~\ref{thmAlexHg} as a definition for the specific needed Alexander polynomial. See Definition~\ref{defAPmv} for a general definition of the normalized multivariable Alexander polynomial. 

In rings $\QQ[[V]]$ of formal power series, for a positive integer $n$, we use the notation $O(n)$ to denote a series of tensors (or monomials) of degree at least $n$.
Set $\sih(u)=\exp\left(\frac{u}2\right) - \exp\left(-\frac{u}2\right)(=2 \sinh(u)).$

\begin{theorem}
\label{thmAlexHg}
Let $\ExtH$ be a genus two rational homology handlebody.
Identify $\partial \ExtH$ with the surface $(-\partial \hbz)$ of Figures~\ref{fighandlebo} and \ref{fighandlebobasedu} in such a way that $H_1(\ExtH;\QQ)=\QQ \ua \oplus \QQ \ub$. Let $\rats$ be the $\QQ$-sphere obtained by attaching $2$-handles along $\ua$ and $\ub$, and by gluing a $3$-ball to the remaining boundary, which is a $2$-dimensional sphere.\footnote{We could state this theorem without the $\QQ$-sphere $\rats$. We have $|H_1(\rats)|=|{H_1(\ExtH)}/{\ZZ \ua \oplus \ZZ \ub}|$ and
$\granda=\ell \ua$ in $H_1(\ExtH;\QQ)/\QQ[\ub]$. The coefficients $\lambda^{\prime}(\granda)$, $\lambda^{\prime}(\grandb)$, and $\lambda^{\prime}(\grandc)$ could be respectively defined from the Reidemeister torsions of $\ExtH[\ua]$, $\ExtH[\ub]$, and $\ExtH[\uc]$ valued in quotients of $\ZZ[H_1(\ExtH;\QQ)]$. See Lemma~\ref{lemlambdaprimeJknot}.}
Set $$
\ddel_{\Delta,\ExtH}(\alpha,\beta)=-\lambda^{\prime}(\granda)\ub\uc -  
\lambda^{\prime}(\grandb) \uc\ua - \lambda^{\prime}(\grandc)\ua\ub$$ and $\ell=lk(\granda,\grandb)$.
Let $\CA_{\ExtH,\CP}$ be a representative of the Alexander form of $\ExtH$.
There exist a unique unit $\unit$ of $\Lambda_{\ExtH}^{\QQ}\subset {\QQ[[\ua,\ub,\uc]]}/(\ua+\ub+\uc=0)$ and a unique rational number $\lambda^{\prime}(\ExtH;\ua,\ub,\uc)$ such that $\CA_{\ExtH}=\unit\CA_{\ExtH,\CP}$ satisfies the following equality.\footnote{Lemma~\ref{lemnormAlexH} constrains the unit $\unit$ to live in a smaller ring.}
$$\CA_{\ExtH}(\tua \wedge \tub)=|H_1(\rats)|\left(1 +\ddel_{\Delta,\ExtH}(\alpha,\beta) + \frac12 \lambda^{\prime}(\ExtH;\ua,\ub,\uc) \ua\ub\uc\right) +O(4).$$
This equality defines $\lambda^{\prime}(\ExtH;\ua,\ub,\uc)$.\footnote{See Corollary~\ref{corlambdaprimethree} for an alternative definition of $\lambda^{\prime}(\ExtH;\ua,\ub,\uc)$.}
If $\ell=0$, then we have $\lambda^{\prime}(\ExtH;\ua,\ub,\uc)=\lambda^{\prime}(\granda,\grandb)$.

With this normalization of $\CA_{\ExtH}$, the multivariable Alexander polynomial $\Delta(\granda,\grandb)$ of $(\granda,\grandb)$ satisfies the following equations.
$$\sih(\uc)\Delta(\granda,\grandb) = \exp\left(-\frac{\ell}{2} \uc \right)\overline{\CA_{\ExtH}(\tua \wedge \tub )} - \exp\left(\frac{\ell}{2} \uc \right)\CA_{\ExtH}(\tua \wedge \tub )$$ 
and
$$\Delta(\granda,\grandb)=\exp\left(\frac{\ell-1}{2}\uc \right)\CA_{\ExtH}(\tub \wedge \tilde{\granda})=\exp\left(\frac{\ell-1}{2}\uc -\ua - \grandb\right)\CA_{\ExtH}(\tua \wedge \tilde{\grandb}). $$
\end{theorem}

Theorem~\ref{thmAlexHg} shows the structure of the restriction of $\CA_{\ExtH}$ to the exterior products of classes of $\partial \ExtH$, as stated in Corollary~\ref{coreasyone} below.
In Section~\ref{subsecpfstructA}, we prove Theorem~\ref{thmAlexHg} and give other corollaries of Theorem~\ref{thmAlexHg}. See Corollaries~\ref{coreasytwo} and \ref{corlambdaprimethree}.

Note the following easy lemma.
\begin{lemma}\label{lemeqABCH} We have
$\granda\ub\granda^{-1}\ub^{-1}=\grandb^{-1}\ua^{-1}\grandb\ua$ in $\pi_1(\partial \ExtH)$, and
$$\partial  \tilde{\granda} \tub - \partial \tub \tilde{\granda} =\exp\left(-\ua - \grandb\right) \left(\partial  \tilde{\grandb} \tua - \partial \tua \tilde{\grandb} \right)$$ in $\CH_{\ExtH}$.
In particular, if $\granda$ and $\grandb$ vanish in $H_1(\ExtH;\QQ)$, then we have
$\partial \tub \tilde{\granda}=\exp(-\ua)\partial \tua \tilde{\grandb}$ in $\CH_{\ExtH}$.
\end{lemma}
\eop

Lemma~\ref{lemeqABCH} implies that $\CA_{\ExtH}(\tilde{\granda} \wedge \tilde{\grandb})=0$
when  $\granda$ and $\grandb$ are rationally null-homologous in $\ExtH$.
In general, it shows how
the restriction of $\CA_{\ExtH}$ to exterior products of classes of curves of $\partial \ExtH$  is determined by $\CA_{\ExtH}(\tua \wedge \tub)$, $\CA_{\ExtH}(\tua \wedge \tilde{\granda})$, and $\CA_{\ExtH}(\tub \wedge \tilde{\granda})$.
The following lemma shows how $\CA_{\ExtH}(\tua \wedge \tilde{\granda})$ is determined by $\CA_{\ExtH}(\tub \wedge \tilde{\granda})$ and $\CA_{\ExtH}(\tua \wedge \tub)$.

\begin{lemma}
\label{lemthreedel}
 For any three elements
$x$, $y$, and $z$ of $\CH_{\ExtH}$, we have
$$\partial x \CA_{\ExtH}(y \wedge z) +\partial y \CA_{\ExtH}(z \wedge x)+ \partial z \CA_{\ExtH}(x \wedge y)=0.$$
In particular, if $z$ is trivial in $H_1(\ExtH;\QQ)$, then we have
$\partial x \CA_{\ExtH}(y \wedge z)=\partial y \CA_{\ExtH}(x \wedge z)$.
\end{lemma}
\bp This is a particular case of Property~2 of \cite[Section 3.2]{lesinv}.
\eop

So, Theorem~\ref{thmAlexHg} has the following easy corollary.
\begin{corollary} \label{coreasyone}
Under the hypotheses of Theorem~\ref{thmAlexHg}, we have
$$\sih(\uc)\CA_{\ExtH}(\tub \wedge \tilde{\granda})=\exp\left(\frac{1-2\ell}{2} \uc \right)\overline{\CA_{\ExtH}(\tua \wedge \tub )} - \exp\left(\frac{1}{2} \uc \right)\CA_{\ExtH}(\tua \wedge \tub ),$$
and the restriction of $\CA_{\ExtH}$ to the exterior products of classes of based curves of $\partial \ExtH$ is determined by $\ell$ and the series $\CA_{\ExtH}(\tua\wedge\tub)$, which is determined by any determinant $\CA_{\ExtH,\CP}(\tua\wedge\tub)$.
\end{corollary}
\eop

In Section~\ref{subdelta}, we prove that $\wdel$ is an invariant of genus one Seifert surfaces $\Sigma$ in $\QQ$-spheres, as announced in Theorem~\ref{thmsimpleinvt}. We furthermore prove that $\wdel$ lifts to an invariant $\ddel_{\Sigma}$ valued in $H_1(\Sigma;\QQ) \otimes_s H_1(\Sigma;\QQ)$, via the form $V_s$ of Notation~\ref{notsymSeif},
as stated in the following theorem.

\begin{theorem}
\label{thmsecondinvariantbis}
Let $\Sigma$ be a genus one Seifert surface in a $\QQ$-sphere $\rats$, and let $(\alpha,\beta)$ be a symplectic basis of $H_1(\Sigma)$. Set $\gamma=-\alpha-\beta$. Represent $\alpha$, $\beta$, and $\gamma$ by simple curves of $\Sigma$.
Define $\petita=-lk(\beta,\gammam)- lk(\betam,\gamma)$,
 $\petitb=-lk(\alpha,\gammam)- lk(\alpham,\gamma)$, and
 $\petitc=-lk(\alpha, \betam)- lk(\alpham, \beta)$.
Define the following elements of $H_1(\Sigma;\QQ) \otimes_s H_1(\Sigma;\QQ)$.
$$\begin{array}{ll}\ddel_{2,\Sigma}(\alpha,\beta)=&-\frac{1}{24} (\petita^2+2\petita\petitb+2\petita\petitc+3)\alpha\otimes_s \alpha -\frac{1}{24}(\petitb^2+2\petitb\petitc+2\petitb\petita+3)\beta \otimes_s \beta \\&-\frac{1}{24}(\petitc^2+2\petitc\petita+2\petitc\petitb+3)\gamma \otimes_s \gamma
\end{array}$$
and $$\ddel_{\Delta,\Sigma}(\alpha,\beta)=\fHone_{\ast}^{-1}(\ddel_{\Delta,\ExtH}(\alpha,\beta))
=-\lambda^{\prime}(\alpha) \beta \otimes_s \gamma - \lambda^{\prime}(\beta) \gamma \otimes_s \alpha - \lambda^{\prime}(\gamma) \alpha \otimes_s \beta.$$
Then the element $$\ddel_{\Sigma}=\ddel_{\Sigma}(\alpha,\beta) = \ddel_{2,\Sigma}(\alpha,\beta) + 4\ddel_{\Delta,\Sigma}(\alpha,\beta)$$ of $H_1(\Sigma;\QQ) \otimes_s H_1(\Sigma;\QQ)$
is an invariant of $\Sigma$ such that $\wdel(\Sigma)=V_s(\ddel_{\Sigma})$.
 \end{theorem}

We will prove the following theorem in the end of Section~\ref{subdelta}, with the help of Theorem~\ref{thmsecondinvariantbis}. It extends the invariant $w_{SL}$ of Theorem~\ref{thmSatoinvariant} to all genus one Seifert surfaces of knots in $\QQ$-spheres.

\begin{theorem}
\label{thmwsl}
 Let $\Sigma$ be a genus one Seifert surface in a $\QQ$-sphere $\rats$. Identify $\Sigma$ to some $\phi\left( \Sigma(\check{a}, \check{b}, \check{c}) \subset \hbz\right)$ for an embedding $\phi \colon \hbz \hookrightarrow \rats$ and a triple $(\check{a}, \check{b}, \check{c})$ of odd integers. Let $\ExtH= \rats \setminus \phi\bigl(\mhbz\bigr)$ be the exterior of $\Sigma$. Let $\alpha=\phi(\alpha(\check{a}, \check{b}, \check{c}))$ and $\beta=\phi(\beta(\check{a}, \check{b}, \check{c}))$ be the corresponding curves of $\Sigma$, as in Figure~\ref{figSigmaabc}.
 Define $\petita=lk(\beta,\alpham + \betam)+ lk(\betam,\alpha + \beta)$,
 $\petitb=lk(\alpha,\alpham + \betam)+ lk(\alpham,\alpha + \beta)$, and
 $\petitc=-lk(\alpha, \betam)- lk(\alpham, \beta)$. Recall the definitions of $\lambda^{\prime}(\ExtH;\ua,\ub,\uc)$ from Theorem~\ref{thmAlexHg} and of $\wdel(\petita,\petitb,\petitc)$ from Theorem~\ref{thmsimpleinvt}.
 Then $$w_{SL}(\Sigma)=\lambda^{\prime}(\ExtH;\ua,\ub,\uc)- \frac1{12} \wdel(\petita,\petitb,\petitc)$$ is an invariant of $\Sigma$.
\end{theorem}

\subsection{More numerical invariants of genus one surfaces}

The torsion $\fHone_{\ast}^{-1}\left(\CD(\Sigma)\right)$ is obviously an invariant of the surface $\Sigma$. So is its degree $(2k)$ part $\fHone_{\ast}^{-1}\left(\CD_{2k}(\Sigma)\right)$, which is valued in the symmetric tensor product
$\bigotimes_s^{2k} H_1(\Sigma;\QQ)$.
We can endow this symmetric product with the canonical linear form 
$V_{s,2k} \colon \bigotimes_s^{2k} H_1(\Sigma;\QQ) \to \QQ$
that maps $(u_1\otimes u_2\otimes\dots \otimes u_{2k})$ to the sum over the partitions of $\{1,2,\dots, 2k\}$ into $k$ pairs $\{n_{2i-1}, n_{2i}\}$ of the products $\prod_{i=1}^{k}V_s(u_{n_{2i-1}} \otimes u_{n_{2i}})$. (There are $\frac{(2k)!}{k!2^k}$ such partitions, and $V_{s,2}=V_s$.) We get the following tautological result.

\begin{proposition}
\label{propsergenone}
For any genus one Seifert surface in a rational homology sphere, 
 $\fHone_{\ast}^{-1}(\CD(\Sigma))$ is an invariant of $\Sigma$ valued in $\ZZ[H_1(\Sigma;\QQ)]$. The extracted series
  $$\biggl(V_{s,2k}\Bigl(\fHone_{\ast}^{-1}\bigl(\CD_{2k}(\Sigma)\bigr)\Bigr)\biggr)_{k \in \NN}$$ is a series of numerical invariants of $\Sigma$.
\end{proposition}

\subsection{Questions}

How does $\wdel$ compare to other knot invariants? Is $\wdel$ a function of the Jones polynomial for knots in $S^3$?
Or does $\wdel$ distinguish knots that are not distinguished by the Jones polynomial?
Do the products by $(7\lambda^{\prime}(\partial \Sigma)-1)$ of the formulas of Corollary~\ref{corwww} hold when $\Sigma$ is not null-homologous?
What are the properties of the series of invariants of genus one Seifert surfaces of Proposition~\ref{propsergenone}?

\section{More on Alexander forms}
\label{secAlef}
\label{AA}\setcounter{equation}{0}

Proposition~\ref{proprtort} allowed us to define Reidemeister torsions from Alexander forms.
In this section, we also use Alexander forms to define normalized Alexander polynomials.

\subsection{Alexander polynomials of knots from Alexander forms}

\begin{definition} \label{defAlexonevorder}
Let $\Jknot$ be a knot in a $\QQ$-sphere $\rats$. Let $O(\Jknot)$ be the minimal positive integer such that $O(\Jknot)\Jknot$ is null-homologous. (When $\Jknot$ is null-homologous, we have $O(\Jknot)=1$.)
Let $X(\Jknot)$ be the exterior of $\Jknot$. Let $m$ be a based meridian of $\Jknot$. Set \begin{equation*}t^{\frac{1}{2O(\Jknot)}}=\exp\left(\frac{1}{2O(\Jknot)}m\right) \in \Lambda_{X(\Jknot)}^{\frac12}.\end{equation*}
The \emph{normalized Alexander polynomial} $\Delta(\Jknot)$ of $\Jknot$ is the Laurent polynomial in $\ZZ\bigl[t^{{1}/{(2O(\Jknot))}},t^{-{1}/{(2O(\Jknot))}}\bigr]$ determined by the following conditions
$$\Delta(\Jknot) \doteq \pm \frac{\exp\left(\frac{1}{O(\Jknot)} m\right)-1}{\exp(m)-1}\CA_{X(\Jknot)}(\tilde{m}),$$ 
$$\Delta(\Jknot)\left(t^{{1}/{(2O(\Jknot))}}=1\right)=\frac{|H_1(\rats)|}{O(\Jknot)}\;\;\mbox{and}\;\; \Delta(\Jknot)(t^{-1/(2O(\Jknot))})=\Delta(t^{1/(2O(\Jknot))}).$$
See \cite[2.1.2 page 22 and 2.2.1 page 23]{lespup}.
\end{definition}

\begin{remark}
\label{rkDeltaonecomp}
Let $g$ be a based curve that represents a generator of $H_1(X(\Jknot))/\mbox{Torsion}$. 
Thanks to the Poincar\'e duality, the algebraic intersection of such a curve with an integral chain whose boundary is $O(\Jknot)\Jknot$ is $\pm 1$. So, we can assume that $m=O(\Jknot)g$ in $H_1(X(\Jknot))/\mbox{Torsion}$. As in Proposition~\ref{proprtort}, we have $\CA_{X(\Jknot)}(\tilde{m}) (\exp(g)-1)=\CA_{X(\Jknot)}(\tilde{g}) (\exp(m)-1)$.
So, we have $\Delta(\Jknot) \doteq \pm \CA_{X(\Jknot)}(\tilde{g})$. Since 
$\CH_{X(\Jknot)}=H_1(\widetilde{X(\Jknot)}) \oplus \Lambda_{X(\Jknot)} [\tilde{g}]$, the polynomial
$\CA_{X(\Jknot)}(\tilde{g})$ is the \emph{order} of $H_1(\widetilde{X(\Jknot)})$, i.e., it is the determinant of a square matrix of presentation of $H_1(\widetilde{X(\Jknot)})$ up to multiplication by a unit.
We also have that 
$$\bigl|H_1(\rats)\bigr|=\Bigl|\varepsilon\bigl(\CA_{X(\Jknot)}(\tilde{m})\bigr)\Bigr|=O(\Jknot)\Bigl|\varepsilon\bigl(\CA_{X(\Jknot)}(\tilde{g})\bigr)\Bigr|=O(\Jknot)\Bigl|\mbox{Torsion}\bigl(H_1(\rats \setminus \Jknot)\bigr)\Bigr|.$$
\end{remark}

 \begin{example} 
 \label{exampleDeltaGranda}
The exterior $X(\granda)$ of the knot $\phi(\granda)$ of the right-hand side of Figure~\ref{fighandlebo} is obtained from $\ExtH= \rats \setminus \phi\bigl(\mhbz\bigr)$ by attaching a $2$-handle along the curve $\ua$ of Figure~\ref{fighandlebobasedu}. So, we have $X(\granda)=\ExtH[\ua]$ and  $H_1(X(\granda);\QQ)={H_1(\ExtH;\QQ)}/{\QQ [\ua]}$. 
Let $\alambicp$ be a representative of $\CA_{\ExtH}$.
The Alexander form $\CA_{X(\granda)}$ is obtained 
 from $\alambicp(\tua \wedge .)$ by mapping $\ua$ to zero (and therefore $\exp(\ua)$ to $1$). We
 write
 $$\CA_{X(\granda)}(u)\doteq \pm \mbox{ev}(\ua=0)\bigl(\alambicp(\tua \wedge u)\bigr).$$
 Since $[-\ub]$ is a meridian of $\granda$ (pushed inside $\phi(\mhbz)$), 
 we have
 $$\Delta\bigl(\phi(\granda)\bigr)(t)\doteq \mbox{sign}\Bigl(\varepsilon \bigl(\alambicp(\tua\wedge\tub)\bigr)\Bigr)\mbox{ev}(\ua=0)\left( \frac{\sih\left(\frac{1}{O(\phi(\granda))}\ub\right)}{\sih(\ub)} \alambicp(\tua \wedge \tub)\right)_{t=\exp(\uc=-\ub)}.$$
 Once Theorem~\ref{thmAlexHg} is proved, we will know that the above expression is even with the normalization of $\CA_{\ExtH}$ of Theorem~\ref{thmAlexHg} instead of $\alambicp$, and we will replace $(\alambicp,\doteq)$ by $(\CA_{\ExtH},=)$.
 \end{example}
 
 \begin{definition} \label{deflambdaprimeJknot}
Let $\Jknot$ be a knot in a $\QQ$-sphere $\rats$. We define
$\lambda^{\prime}(\Jknot)$ from the Alexander polynomial $\Delta(\Jknot)$ of Definition~\ref{defAlexonevorder} to be
$$\lambda^{\prime}(\Jknot) = \frac{O(\Jknot)\Delta^{\prime\prime}(\Jknot)(1)}{2|H_1(\rats)|} + \frac1{24}\left(1-\frac1{O(\Jknot)^2}\right),$$
where the second derivative of $\Delta$ is taken with respect to the variable $t=\exp(m)$.
\end{definition}

Lemma~\ref{lemlambdaprimeJknot} below justifies the above definition. It provides an alternative definition for $\lambda^{\prime}(\Jknot)$.

\begin{lemma}
\label{lemlambdaprimeJknot} Let $\Jknot$ be a knot in a $\QQ$-sphere $\rats$. Let $m$ be a based meridian of $\Jknot$.
Let $\CA_{X(\Jknot),\CP}$ be a representative of the Alexander form of $X(\Jknot)$. There exists a unit $\unit$ of $\Lambda_{\extx(\Jknot)}^{\frac12} \subset \QQ[[m]]$ such that 
$$\Delta(\Jknot)=\frac{\sih\left(\frac{1}{O(\Jknot)} m \right)}{\sih(m)} \unit\CA_{X(\Jknot),\CP}(\tilde{m})$$
and 
 $$\unit\CA_{X(\Jknot),\CP}(\tilde{m})=|H_1(\rats)|\left(1 + \lambda^{\prime}(\Jknot)m^2 \right) +O(4).$$
 If $\rats$ is a $\frac{\ZZ}{2\ZZ}$--sphere or if $\Jknot$ is null-homologous, then $\unit$ and $\Delta(\Jknot)$ are in $\Lambda_{\extx(\Jknot)}$.
\end{lemma}
\bp According to Remark~\ref{rkDeltaonecomp}, there exists a based curve
$g$ such that $m=O(\Jknot)g$ in $H_1(X(\Jknot))/\mbox{Torsion}$. Then we have ${H_1(X(\Jknot))}/{\mbox{Torsion}}=\ZZ g$.
Without loss of generality, we choose $\CA_{X(\Jknot),\CP}$ so that $\varepsilon\left(\CA_{X(\Jknot),\CP}(\tilde{g})\right) >0$.
By Definition~\ref{defAlexonevorder}, we have $$\Delta(\Jknot) \doteq \pm \frac{\partial{\tilde{g}}}{\partial{\tilde{m}}}\CA_{X(\Jknot),\CP}(\tilde{m}).$$
So $\Delta(\Jknot) \doteq \pm \CA_{X(\Jknot),\CP}(\tilde{g})$, and
there exists $k_1\in \ZZ$ such that 
$\Delta(\Jknot)=\exp\left(\frac12 k_1 g\right)\CA_{X(\Jknot),\CP}(\tilde{g})$.
Since $\CA_{X(\Jknot),\CP}(\tilde{g}) \in 
\Lambda_{\extx(\Jknot)}$,
we have $\CA_{X(\Jknot),\CP}(\tilde{g})={|H_1(\rats)|}/{O(\Jknot)} + kg +O(2) $ in $\QQ[[g]]$ for some $k$ in $\ZZ$.
Since $\Delta(\Jknot)$ is even, the degree one part of $\exp(k_1 g/2)\CA_{X(\Jknot),\CP}(\tilde{g})$ is zero. 
So, we have $k_1= -2\frac{kO(\Jknot)}{|H_1(\rats)|} \in \ZZ$, and $k_1$ is even if $\frac{|H_1(\rats)|}{O(\Jknot)}$ is odd, which holds as soon as $\rats$ is a $\frac{\ZZ}{2\ZZ}$-sphere.
In this case, $\Delta(\Jknot)$ belongs to $\Lambda_{\extx(\Jknot)}$.
Set $$\unit=\exp\left(-\frac{O(\Jknot)-1}2g\right)\exp\left(\frac12 k_1 g\right) \;\; \mbox{and}\;\;\CA_{X(\Jknot)}= \unit\CA_{X(\Jknot),\CP}.$$ Then $$\CA_{X(\Jknot)}(\tilde{m})=\exp\left(-\frac{O(\Jknot)-1}2g\right)\frac{\partial{\tilde{m}}}{\partial{\tilde{g}}}\exp\left(\frac12 k_1 g\right)\CA_{X(\Jknot),\CP}(\tilde{g}) =\frac{\sih(O(\Jknot) g)}{\sih(g)} \Delta(\Jknot)$$ is even, and we have
$$\CA_{X(\Jknot)}(\tilde{m})=|H_1(\rats)|(1 + \cconst m^2) + O(4).$$

Since $\sih(x)=\exp\left(\frac{x}2\right) - \exp\left(-\frac{x}2\right)=x\left(1+\frac{x^2}{24}+O(4) \right),$
we get
$$\Delta(\Jknot)
=\frac{\sih\left(\frac1{O(\Jknot)} m\right)}{\sih(m)}\CA_{X(\Jknot)}(\tilde{m}) =\frac{|H_1(\rats)|}{O(\Jknot)}\left(1 - \frac1{24}\left(1-\frac1{O(\Jknot)^2}\right)m^2 \right)\left( 1 + \cconst m^2\right)+O(4).$$
 On the other hand, we have
$$\Delta(\Jknot)=\frac{|H_1(\rats)|}{O(\Jknot)} +\frac{\Delta^{\prime\prime}(\Jknot)(1)}{2}m^2 +O(4).$$
So, we get \begin{equation*}
\cconst= \frac1{24}\left(1-\frac1{O(\Jknot)^2}\right) +  \frac{O(\Jknot)\Delta^{\prime\prime}(\Jknot)(1)}{2|H_1(\rats)|}=\lambda^{\prime}(\Jknot).\end{equation*}

Since $(O(\Jknot)-1)$ is even when $\rats$ is a $\frac{\ZZ}{2\ZZ}$--sphere, $\unit$ belongs to $\Lambda_{\extx(\Jknot)}$ in this case. When $\Jknot$ is null-homologous, $O(\Jknot)$ equals one, and the expression  $\Delta(\Jknot)=|H_1(\rats)|\mbox{det}\left(t^{1/2} V- t^{-1/2} (^TV)\right) $ shows that $\Delta(\Jknot)$ is in $\Lambda_{\extx(\Jknot)}$. So, we have $k_1\in 2\ZZ$ and $\unit \in \Lambda_{\extx(\Jknot)}$, too.
\eop

\subsection{Normalized Alexander polynomials of links}

\begin{definition}
 \label{defAPmv} Let $k$ be an integer such that $k \geq 2$. Set $\underline{k}=\{1,\dots,k\}$.
Let $\Link=\left(K_i\right)_{i  \in \underline{k}}$ be a link of a $\QQ$-sphere $\rats$ with $k$ components. For $i \in \underline{k}=\{1,\dots,k\}$, let $m_i$ be a meridian of the component $K_i$ of $\Link$, and set $O_i=O(K_i)$.
The Reidemeister torsion $\tau(\extx(\Link))$ of the exterior $\extx(\Link)$ of $\Link$ is an element of
$$\ZZ\left[ e^{ \frac{1}{O_1}m_1}, e^{- \frac{1}{O_1}m_1},e^{ \frac{1}{O_2}m_2}, e^{- \frac{1}{O_2}m_2},\dots, e^{ \frac{1}{O_k}m_k}, e^{- \frac{1}{O_k}m_k}\right].\footnote{Here, we see $H_1(\extx(\Link))/{\mbox{Torsion}}$ as a submodule of $H_1(\extx(\Link);\QQ)=\oplus_{i=1}^k\QQ m_i$. Let $N(K_i)$ denote a tubular neighborhood of $K_i$. Let $\Sigma_i$ be a surface embedded in the exterior $\rats \setminus \mathring{N}(K_i)$ of $K_i$
whose boundary $\partial \Sigma_i$ is on $\partial N(K_i)$ and is homologous to $O_iK_i$ in $N(K_i)$. Every element $x$ of $H_1(\extx(\Link);\QQ)=\oplus_{i=1}^k\QQ m_i$ may be written as $x=\sum_{i=1}^k ({\langle x, \Sigma_i\rangle}/{O_i}) m_i$.}$$
It is defined up to multiplication by a unit $\exp\bigl(x \in H_1\bigl(\extx(\Link)\bigr)/\mbox{Torsion}\bigr)$. The multivariable Alexander polynomial $\Delta(\Link)$ of $\Link$ is an element of $\ZZ\bigl[t_1^{{1}/{(2O_1)}},t_1^{-{1}/{(2O_1)}},\dots, t_k^{{1}/{(2O_k)}},  t_k^{-{1}/{(2O_k)}}\bigr]$. It is obtained from $\tau(\extx(\Link))$, by a multiplication by a unit of $\Lambda^{\frac12}_{\extx(\Link)}$, 
and by the change of variables \begin{equation*}t_i^{\frac{1}{2O_i}}=\exp\left(\frac1{2O_i} m_i\right).\end{equation*} It satisfies the following property $$\Delta(\Link)\bigl(t_1^{\frac{1}{2O_1}}, t_2^{\frac{1}{2O_2}}, \dots, t_k^{\frac{1}{2O_k}}\Bigr) = \pm \Delta(\Link)\Bigl(t_1^{-\frac{1}{2O_1}}, t_2^{-\frac{1}{2O_2}}, \dots, t_k^{-\frac{1}{2O_k}}\Bigr),$$
which determines it up to sign.
To fix the sign of $\Delta$, as Turaev did in \cite{turaevReid}, we use Alexander forms as follows.

Let $B$ be a $3$-ball in $\extx(\Link)$ that intersects the boundary $\partial N(K_i)$ of the tubular neighborhood $N(K_i)$ of its component $K_i$ along a disk $D_i$, for each component $K_i$ of $\Link$.
So $B \cup \cup_{i=1}^k N(K_i)$ is a handlebody in $\rats$, and $\extx(\Link)$ is obtained from the exterior $\ExtH$ of this handlebody by attaching the $2$-cells $D_1, \dots, D_{k-1}$ to $\partial \ExtH$ and by thickening them. Let $\ell_i$ be a longitude of $K_i$, we have $\langle m_i,\ell_i \rangle_{\partial N(K_i)} =1$. Attach $m_i$ and $\ell_i$ to a basepoint in $\partial \ambe \cap B$ by a path on $\partial \ambe \cap (\partial N(K_i) \cup B)$ so that the boundary $\cdel_i$ of $\pm D_i$ reads $m_i\ell_im_i^{-1}\ell_i^{-1} $ up to conjugation in $\pi_1(\partial \ambe)$. 
Then $\tilde{\cdel}_i=(1 -\exp(\ell_i))\tilde{m}_i +(\exp(m_i)-1)\tilde{\ell}_i$ in $\CH_{\ambe}$, up to multiplication by a positive unit of $\Lambda_{\ambe}$. 
According to \cite[Property 8, p. 647]{lesinv},
for any $j$ and $r$ in $\underline{k}$, the \emph{multivariable Alexander polynomial $\Delta(\Link)$ of $\Link$}  satisfies
$$ \Delta(\Link) \doteq\mbox{sign}\Bigl(\varepsilon \bigl(\CA_{\ExtH}(\hat{m})\bigr)\Bigr) \frac{\CA_{\ExtH}\left(\hat{\cdel}\left(\frac{\tilde{m}_j}{\tilde{\cdel}_r}\right)\right)}{\partial \tilde{m}_j},$$
 where $\hat{m} =\tilde{m}_1 \wedge \dots \wedge \tilde{m}_k$, $\hat{\cdel} =\tilde{\cdel}_1 \wedge \dots \wedge \tilde{\cdel}_k$,
$\hat{\cdel}\left(\frac{\tilde{m}_j}{\tilde{\cdel}_r}\right)= \tilde{\cdel}_1 \wedge \dots \wedge \tilde{\cdel}_{r-1} \wedge \tilde{m}_j \wedge \tilde{\cdel}_{r+1} \wedge \dots \wedge \tilde{\cdel}_k,$
and we keep the change of variables $t_i^{{1}/{(2O_i)}}=\exp\bigl(\frac1{2O_i} m_i\bigr)$ understood.
See also \cite[2.2.1 page 23]{lespup}.
 \end{definition}

\begin{example} \label{exatwocomp}
  For a general two-component link $\Link=(K_1,K_2)$ as above, we have
  $$ \begin{array}{ll}\Delta(\Link)\left(t_1=\exp(m_1), t_2
 =\exp(m_k)\right)
 &\doteq\mbox{sign}\Bigl(\varepsilon \bigl(\CA_{\ExtH}(\tilde{m}_1 \wedge \tilde{m}_2)\bigr)\Bigr) \frac{\CA_{\ExtH}\left((1 -\exp(\ell_1))\tilde{m}_1 +(\exp(m_1)-1)\tilde{\ell}_1 \wedge \tilde{m}_1\right)}{\partial \tilde{m}_1}\\
 &\doteq\mbox{sign}\Bigl(\varepsilon \bigl(\CA_{\ExtH}(\tilde{m}_1 \wedge \tilde{m}_2)\bigr)\Bigr) \CA_{\ExtH}\left(\tilde{\ell}_1 \wedge \tilde{m}_1\right)\\
 &\doteq\mbox{sign}\Bigl(\varepsilon \bigl(\CA_{\ExtH}(\tilde{m}_1 \wedge \tilde{m}_2)\bigr)\Bigr) \CA_{\ExtH}\left(\tilde{m}_2 \wedge \tilde{\ell}_2\right).
 \end{array}$$
 Since $\ell_1=lk(K_1,K_2)m_2$ in $H_1(\extx(\Link))/\mbox{Torsion}$, we have
 $\varepsilon(\Delta(\Link))=-lk(K_1,K_2)|H_1(\rats)|$.
 Observe that $|\varepsilon(\Delta(\Link))|$ is the order of the quotient of $H_1(\extx(\Link))$ by the image of the $H_1$ of one of the boundary components.
In particular, for our genus one knot $K$ and the two-component link whose
exterior is the manifold $\ExtH[K]$ obtained from $\ExtH$ by adding a $2$-handle along $K$, $|\varepsilon(\CD(\Sigma))|$ is $$\biggl|\varepsilon\Bigl(\CA_{\ExtH}\bigl(\tilde{\alpha} \wedge \tilde{\beta}\bigr)\Bigr)\biggr|=\biggl|\varepsilon\Bigl(\CA_{\ExtH}\bigl(\talpham \wedge \tbetam\bigr)\Bigr)\biggr|=\Bigl|H_1(\rats)\Bigr|\Bigl|\mbox{det}(V)\Bigr|=\Bigl|H_1(\rats)\Bigr|\Bigl|\lambda^{\prime}\bigl(\phi(K)\bigr)\Bigr|$$
where $V$ is the Seifert matrix of $\Sigma$.
 \end{example}
 
 \begin{remark}
  The condition that $\cdel_i$ is conjugate to $ m_i\ell_im_i^{-1}\ell_i^{-1} $ is satisfied when $m_i$ and $\ell_i$ are represented by two simple closed curves that intersect at one transverse point and when they are attached to the basepoint by the same path from the intersection point to the basepoint of $\partial \ExtH$.
Figure~\ref{figsurfmili} shows the expression of $\cdel_i$ in the $\pi_1$ of the regular neighborhood of $m_i \cup \ell_i$ with respect to the shown basepoints and paths.

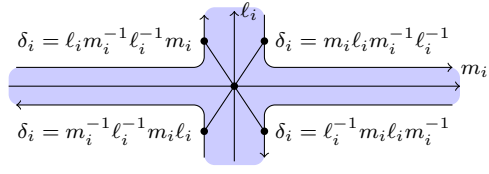
\begin{figure}[h]
\begin{center}
\begin{tikzpicture}
 \useasboundingbox (-4,-1.1) rectangle (4,1.1);
\fill[rounded corners,blue!20] (-3,0) -- (-3,.25) -- (-.4,.25) -- (-.4,1.05) -- (.4,1.05) -- (.4,.25) --  (3,.25)  --  (3,-.25) -- (.4,-.25) -- (.4,-1.05) -- (-.4,-1.05) -- (-.4,-.25) -- (-3,-.25) --  (-3,0);
\draw[->] (3.2,.2) node{\scriptsize $m_i$} (-3,0) -- (3,0);
\draw[->] (.2,1) node{\scriptsize $\ell_i$} (0,-1) -- (0,1);
\draw [rounded corners,->] (-2.9,.25) -- (-.4,.25) -- (-.4,.95);
\draw (-.4,.6) node[left]{\scriptsize $\cdel_i=\ell_im_i^{-1}\ell_i^{-1}m_i$} (-.4,.6) -- (0,0);
\draw [rounded corners,<-] (2.9,.25) -- (.4,.25) -- (.4,.95);
\draw (.4,.6) node[right]{\scriptsize $\cdel_i=m_i\ell_im_i^{-1}\ell_i^{-1}$} (.4,.6) -- (0,0);
\draw [rounded corners,<-] (-2.9,-.25) -- (-.4,-.25) -- (-.4,-.95);
\draw (-.4,-.6) node[left]{\scriptsize $\cdel_i=m_i^{-1}\ell_i^{-1}m_i\ell_i$} (-.4,-.6) -- (0,0);
\draw [rounded corners,->]  (2.9,-.25) -- (.4,-.25) -- (.4,-.95);
\draw (.4,-.6) node[right]{\scriptsize $\cdel_i=\ell_i^{-1}m_i\ell_im_i^{-1}$} (.4,-.6) -- (0,0);
\fill (0,0) circle (.05) (-.4,.6) circle (.05) (-.4,-.6) circle (.05) (.4,-.6) circle (.05) (.4,.6) circle (.05);
\end{tikzpicture}
\caption{A piece of the genus one surface obtained by plumbing tubular neighborhoods of $m_i$ and $\ell_i$}
\label{figsurfmili}
\end{center}
\end{figure}
Be aware that $\cdel_i$ might no longer be conjugate to $m_i\ell_im_i^{-1}\ell_i^{-1} $ if one of the curves is attached to the basepoint by another path.
For example, replace $\ell_i$ by $\check{\ell_i}$ so that $\check{\ell_i}=\cdel_i\ell_i\cdel_i^{-1}$ in $\pi_1(\ambe)$, 
then $\check{\cdel_i}=m_i\check{\ell_i}m_i^{-1}\check{\ell_i}^{-1}$ and
$$\tilde{\check{\cdel_i}}=\tilde{\cdel_i} + \bigl(\exp(m_i) -1\bigr)\bigl(1-\exp(\ell_i)\bigr)\tilde{\cdel_i}.$$
So $\tilde{\check{\cdel_i}}$ and $\tilde{\cdel_i}$ are not obtained from one another by a multiplication by a positive unit of $\Lambda_{\ambe}$ when $\ell_i$ is nonzero in $H_1\bigl(\extx(\Link);\QQ\bigr)$, and $\cdel_i$ and $\check{\cdel_i}$ are not conjugate in this case.
 \end{remark}

 \begin{example} 
 \label{exampleDeltaGrandaB}
For the links $(\phi(\granda),\phi(\grandb))$, $(\phi(\alpha),\phi(\beta_-))$, and $(\phi(\alpha_-),\phi(\beta))$, the meridian of $\phi(\granda)$, $\phi(\alpha)$, or $\phi(\alpha_-)$ is $\ub^{-1}$, and the meridian of $\phi(\grandb)$, $\phi(\beta_-)$, or $\phi(\beta)$ is $\ua$.
If $\varepsilon \left(\CA_{\ExtH}(\ub^{-1}\wedge\ua)\right) >0$, then we have
$$\begin{array}{ll}\Delta\bigl(\phi(\granda),\phi(\grandb)\bigr)\bigl(t_{\granda}=\exp(m_{\granda}=-\ub), t_{\grandb}=\exp(m_{\grandb}=\ua)\bigr)&\doteq \CA_{\ExtH}(\tua\wedge \tilde{\grandb}),\\
\Delta\bigl(\phi(\alpha),\phi(\beta_-)\bigr) &\doteq \CA_{\ExtH}\left(\tilde{\alpha}\wedge\widetilde{\ub^{-1}}\right),\\
\Delta\bigl(\phi(\alpha_-),\phi(\beta)\bigr) &\doteq \CA_{\ExtH}(\tua\wedge\tbeta).  
  \end{array}
$$
 \end{example}

\begin{definition} \label{defzetatwocomp} Let $(K_1,K_2)$ be a two-component link in a $\QQ$-sphere $\rats$.
 The coefficient $\zeta(K_1,K_2)$ is the coefficient of $m_1m_2$ in the multivariable Alexander polynomial $\Delta$ of Definition~\ref{defAPmv} viewed as a series in $\QQ[[m_1,m_2]]$. \end{definition}
 
 The following lemma will allow us to generalize the description of $\lambda^{\prime}(K_1,K_2)$ of Definition~\ref{defzetatwo} in Definition~\ref{deflambdaprimetwocomp}.
 
 \begin{lemma} \label{lemzetatwocomp}
  Let $(K_1,K_2)$ be a $2$-component link in a $\QQ$-sphere $\rats$ such that
  $lk(K_1,K_2)=0$.
  For $i=1,2$, let $N(K_i)$ denote a tubular neighborhood of $K_i$, and let $\Sigma_i$ be a surface embedded in $\rats \setminus \bigl(N(K_{3-i}) \sqcup \mathring{N}(K_i)\bigr)$
such that the boundary $\partial \Sigma_i$ of $\Sigma_i$ consists of $r_i$ parallel copies of a primitive curve $\lambda_i$ on $\partial N(K_i)$, with $r_i>0$, and $\partial \Sigma_i$ is homologous to $O_iK_i$ in $N(K_i)$, where $O_i=O(K_i)$. Assume that $\Sigma_1$ and $\Sigma_2$ are transverse. 
  Let $\Sigma_{1,\parallel}$ be obtained from $\Sigma_{1}$ by a slight normal push. Then the Alexander polynomial $\Delta(K_1,K_2)$ of $(K_1,K_2)$ satisfies the following equation in $\QQ[[m_1,m_2]]$
$$\Delta(K_1,K_2)\left(\exp\left(\frac{1}{2O_1}m_1\right),\exp\left(\frac{1}{2O_2}m_{2}\right)\right)=\zeta(K_1,K_2)m_1m_2 + O(4),$$
where $$\zeta(K_1,K_2)=-\frac{|H_1(\rats)|}{r_1r_2O_1O_2}lk (\Sigma_{1} \cap \Sigma_{2},\Sigma_{1,\parallel} \cap \Sigma_{2}).$$
 \end{lemma}
 \bp The degree $0$ part of $\Delta(K_1,K_2)$ is $-lk(K_1,K_2)|H_1(\rats)|$, as in Example~\ref{exatwocomp}. According to the surgery formula \cite[1.4.8 page 12,  T2 page 13]{lespup}, the invariant $\lambda$ of the manifold $\chi=\rats\bigl((K_1;0),(K_2;0)\bigr)$ obtained from $\rats$ by $0$-surgery on $(K_1,K_2)$ is \begin{equation*}\lambda(\chi) =\zeta(K_1,K_2).\end{equation*}
 Define the quadratic form $q \colon  H_2(\chi) \wedge H_2(\chi) \to \QQ$ that maps 
 the exterior product of the homology classes of two transverse surfaces $F_1$ and $F_2$ in $\chi$ to $lk({F}_1 \cap {F}_2, {F}_{1,\parallel} \cap {F}_2)$.
 Let $([{A}_1],[{A}_2])$ be a basis of $H_2(\chi)=\ZZ [{A}_1] \oplus \ZZ [{A}_2]$. 
 Then \cite[T5.2, page 13]{lespup} implies \begin{equation*}
   \lambda(\chi) =-|\mbox{Torsion}(H_1(\chi))|q([{A}_1] \wedge [{A}_2]).                                        
                                           \end{equation*}
                                           
For $i \in \{1,2\}$, let $T_i$ denote the solid torus in $\chi$
with meridian $\lambda_i$ added during the $0$--surgery along $K_i$. 
Let $\extx=\extx(K_1,K_2)$ be the exterior of $(K_1,K_2)$.
 There is a canonical isomorphism from $H_2(\chi)$ to $H_2(\extx,\partial \extx)$. It maps the class of a surface $F$ in $\chi$ transverse to $T_i$ to the class $[\check{F}=F\cap \extx]$.
We can assume that the above basis $([{A}_1],[{A}_2])$ of $H_2(\chi)$ is represented by surfaces ${A}_1$ and ${A}_2$ intersecting the $T_i$ as meridian disks of the $T_i$ and we do.
For $i \in \{1,2\}$, set $\check{A}_i ={A}_i \cap \extx$ and let $\hat{\Sigma}_i$ be the closed surface of $\chi$ obtained from $\Sigma_i$ by gluing meridian disks of $T_i$ along $\partial \Sigma_i$. 
Set   \begin{equation*}r=\left| \frac{\ZZ [\lambda_1] \oplus \ZZ [\lambda_2]}{\ZZ [ \partial \check{A}_1] \oplus \ZZ [ \partial \check{A}_2]}\right|. \end{equation*} We have 
\begin{equation*}
q([{A}_1] \wedge [{A}_2])=\frac{r^2}{r^2_1r^2_2}q([\hat{\Sigma}_1] \wedge [\hat{\Sigma}_2]). \end{equation*}                                           
For $i \in \{1,2\}$, let $\mu_i$ be a curve on the boundary of the solid torus $T_i$ such that $\langle \mu_i,\lambda_i \rangle_{\partial T_i}=1$.
 The Mayer--Vietoris sequence associated with the decomposition of $\chi$ as the union of $\extx$ and $T_1 \cup T_2$ leads
 to an exact sequence 
 $$\frac{\ZZ [\lambda_1] \oplus \ZZ [\lambda_2]}{\ZZ [ \partial \check{A}_1] \oplus \ZZ [ \partial \check{A}_2]}  \oplus \ZZ \mu_1 \oplus \ZZ \mu_2 \hookrightarrow H_1(\extx) \oplus H_1(T_1) \oplus H_1(T_2)\to H_1(\chi) \to 0,$$ 
 which implies $\left|\mbox{Torsion}\bigl(H_1(\extx)\bigr)\right|=\left|\mbox{Torsion}\bigl(H_1(\chi)\bigr)\right|r$. We obtain $$\zeta(K_1,K_2)=-\frac{\left|\mbox{Torsion}\bigl(H_1(\extx)\bigr)\right|}{r}\frac{r^2}{r^2_1r^2_2} q([\hat{\Sigma}_1] \wedge [\hat{\Sigma}_2]).$$
Since $H_1(\rats)=\frac{H_1(\extx)}{\ZZ m_1 \oplus \ZZ m_2}$, we have
 $|H_1(\rats)|=\Bigl|\det\bigl(\langle m_i,\check{A}_j \rangle_{\rats} \bigr)\Bigr| \left|\mbox{Torsion}\bigl(H_1(\extx)\bigr)\right|$.
 Thus $$|H_1(\rats)|
 =\frac{r\bigl|\langle m_1,\Sigma_1 \rangle_{\rats}\langle m_2,\Sigma_2 \rangle_{\rats}\bigr|}{r_1r_2}\Bigl|\mbox{Torsion}\bigl(H_1(\extx)\bigr)\Bigr|= \frac{r}{r_1r_2}O_1O_2\left|\mbox{Torsion}\bigl(H_1(\extx)\bigr)\right|.$$
 Therefore, we obtain the desired formula for $\zeta(K_1,K_2)$.
 \eop
 
 \begin{definition} \label{deflambdaprimetwocomp} Under the assumptions of Definition~\ref{defzetatwocomp}, we set $\lambda^{\prime}(K_1,K_2)= \frac{\zeta(K_1,K_2)}{|H_1(\rats)|}$. \end{definition}

 \subsection{A few lemmas about the structure of Alexander forms}

In this subsection, $\ExtH$ is a rational homology genus two handlebody as in Theorem~\ref{thmAlexHg}.
We start the proof of this structure theorem for the Alexander form of such rational homology genus two handlebodies by proving preliminary lemmas. 
Let $\CA_{\ExtH,n}(\tua \wedge \tub)$ denote the degree $n$ part of $\CA_{\ExtH}(\tua \wedge \tub)=\sum_{n \in \NN}\CA_{\ExtH,n}(\tua \wedge \tub)$.

 \begin{lemma}
 \label{lemnormAlexH} Under the hypotheses of Theorem~\ref{thmAlexHg},
 let $\CA_{\ExtH,\CP}$ be a representative of the Alexander form of $\ExtH$. There exists a unique unit $\unit$ of $\Lambda_{\ExtH}^{\QQ}\subset \frac{\QQ[[\ua,\ub,\uc]]}{\ua+\ub+\uc=0}$ 
 such that $\CA_{\ExtH}=\unit\CA_{\ExtH,\CP}$ satisfies 
 $$\CA_{\ExtH}(\tua \wedge \tub)=|H_1(\rats)|
  +\sum_{n \in \NN, n\geq 2}\CA_{\ExtH,n}(\ua \wedge \ub)$$
  where 
  $$\CA_{\ExtH,2}(\tua \wedge \tub)= |H_1(\rats)|\ddel_{\Delta,\ExtH}(\alpha,\beta)$$
  and $\CA_{\ExtH,3}(\tua \wedge \tub)=\frac{|H_1(\rats)|}{2} \lambda^{\prime}(\ExtH;\ua,\ub,\uc)\ua\ub\uc$ for some $\lambda^{\prime}(\ExtH;\ua,\ub,\uc) \in \QQ$.
 
Furthermore, when $\granda$ and $\grandb$ are null-homologous in $\rats$, $\unit$ is a unit of $\Lambda_{\ExtH}$. In particular, when $\granda$ and $\grandb$ are null-homologous,
 $\CA_{\ExtH}(v \wedge w)\in \Lambda_{\ExtH}$ for any $v \wedge w \in \CH_{\ExtH} \wedge \CH_{\ExtH}$.
 In general, $\unit$ is a unit of 
 \begin{multline*}\Lambda(\ExtH;\ua,\ub,\uc)=\ZZ\left[\exp\Bigl( \frac{\ua}{ 2O(\grandb)}\Bigr), \exp\Bigl(-  \frac{\ua}{ 2O(\grandb)}\Bigr),\exp\Bigl( \frac{\ub}{
 2O(\granda)}\Bigr), \exp\Bigl(-  \frac{\ub}{ 2O(\granda)}\Bigr) \right] \\ \cap 
 \ZZ\left[e^{ \frac{\ub}{2O(\grandc)}}, e^{-  \frac{\ub}{2O(\grandc)}},e^{ \frac{\uc}{2O(\grandb)}}, e^{-  \frac{\uc}{2O(\grandb)}} \right] 
 \cap\ZZ\left[e^{ \frac{\uc}{2O(\granda)}}, e^{-  \frac{\uc}{2O(\granda)}},e^{ \frac{\ua}{2O(\grandc)}}, e^{-  \frac{\ua}{2O(\grandc)}} \right].
 \end{multline*}

 If $|H_1(\rats)|$ is odd (or equivalently if $\rats$ is a $\ZZ/2\ZZ$-sphere), then $\unit$
 is a unit of \begin{multline*}\Lambda_o(\ExtH;\ua,\ub,\uc)= \ZZ\left[\exp\Bigl( \frac{\ua}{ O(\grandb)}\Bigr), \exp\Bigl(-  \frac{\ua}{ O(\grandb)}\Bigr),\exp\Bigl( \frac{\ub}{O(\granda)}\Bigr), \exp\Bigl(-  \frac{\ub}{ O(\granda)}\Bigr)\right] \\ \cap 
 \ZZ\left[e^{ \frac{\ub}{O(\grandc)}}, e^{-  \frac{\ub}{O(\grandc)}},e^{ \frac{\uc}{ O(\grandb)}}, e^{-  \frac{\uc}{ O(\grandb)}} \right] 
 \cap\ZZ\left[e^{ \frac{\uc}{O(\granda)}}, e^{-  \frac{\uc}{O(\granda)}},e^{ \frac{\ua}{O(\grandc)}}, e^{-  \frac{\ua}{O(\grandc)}} \right].\end{multline*}

 \end{lemma}
 \bp Recall from Example~\ref{exampleDeltaGrandaB} that $\ub^{-1}$ is a meridian of $\granda$ and that $\ua$ is a meridian of $\grandb$, when $\granda$ and $\grandb$ are slightly pushed inside $\phi(\mhbz)$.
 There is no loss of generality in assuming that $\varepsilon\bigl(\CA_{\ExtH,\CP}(\tua \wedge \tub)\bigr) >0$, and we do.
According to Lemma~\ref{lemlambdaprimeJknot}, there exists $(k_a,k_b)\in \ZZ^2$ such that 
$$|H_1(\rats)|\left(1 + \lambda^{\prime}(\granda)\ub^2 \right) = \exp\left(\frac{k_a \ub}{2O(\granda)}\right)\mbox{ev}(\ua=0)\left( \CA_{\ExtH,\CP}(\tua \wedge \tub)\right)  +O(4)$$
and
$$|H_1(\rats)|\left(1 + \lambda^{\prime}(\grandb)\ua^2 \right) = \exp\left(\frac{k_b \ua}{2O(\grandb)}\right)\mbox{ev}(\ub=0)\left( \CA_{\ExtH,\CP}(\tua \wedge \tub)\right)  +O(4)$$
as in Example~\ref{exampleDeltaGranda}.
Set $\unit=\exp(\frac{k_b}{2O(\grandb)}\ua +\frac{k_a}{2O(\granda)}\ub)$ and $\CA_{\ExtH}=\unit\CA_{\ExtH,\CP}$. Then
$$\CA_{\ExtH}(\tua \wedge \tub)=|H_1(\rats)|+\sum_{n \in \NN, n\geq 2}\CA_{\ExtH,n}(\tua \wedge \tub)$$ 
where $$\CA_{\ExtH,2}(\tua \wedge \tub)=|H_1(\rats)|\left(\lambda^{\prime}(\granda)\ub^2 + \lambda^{\prime}(\grandb)\ua^2 + c(a,b)\ua\ub\right).$$
The coefficients of $\ua^n$ and $\ub^n$ in the degree $n$ part $\CA_{\ExtH,n}(\tua \wedge \tub)$ of $\bigl(\CA_{\ExtH}(\tua \wedge \tub) \in \QQ[[\ua,\ub]]\bigr)$ are determined by the Alexander polynomials of $\granda$ and $\grandb$. They vanish for odd $n$.
Since $\tua\wedge\tub=\exp(\ub)\tuc \wedge \tua$ in $\wedge^2\CH_{\ExtH}$ according to Lemma~\ref{lemuaubuc}, and because of the order $3$ symmetry of the setting, $\tilde{\CA_{\ExtH}}=\exp(\ub)\CA_{\ExtH}$ satisfies
$$\tilde{\CA_{\ExtH}}(\tuc \wedge \tua)=\CA_{\ExtH}(\tua \wedge \tub)
=|H_1(\rats)|\left(1 + \lambda^{\prime}(\grandc)\ua^2 + \lambda^{\prime}(\granda)\uc^2 + c(c,a)\uc\ua\right)$$
in $\QQ[[\ua,\uc]]=\frac{\QQ[[\ua,\ub,\uc]]}{\ua+\ub+\uc=0}$.
So, we get $$\begin{array}{ll}\frac{\CA_{\ExtH,2}(\tua \wedge \tub)}{|H_1(\rats)|}
&=\lambda^{\prime}(\granda)\ub(-\ua-\uc) + \lambda^{\prime}(\grandb)\ua(-\ub-\uc) +c(a,b)\ua\ub \\
&= \lambda^{\prime}(\grandc)\ua(-\ub-\uc) +
\lambda^{\prime}(\granda)\uc(-\ua-\ub) + c(c,a)\uc\ua\\
&=-\lambda^{\prime}(\granda) \ub \uc -\lambda^{\prime}(\grandb) \uc \ua  -\lambda^{\prime}(\grandc) \ua \ub\\
&=\ddel_{\Delta,\ExtH}(\alpha,\beta).
\end{array}$$
Furthermore, $\unit \in \ZZ\left[e^{ \frac{\ua}{2O(\grandb)}}, e^{-  \frac{\ua}{2O(\grandb)}},e^{ \frac{\ub}{2O(\granda)}}, e^{-  \frac{\ub}{2O(\granda)}} \right]$,
and $$\tilde{\unit}=\exp(\ub)\unit \in \ZZ\left[e^{ \frac{\ua}{2O(\grandc)}}, e^{-  \frac{\ua}{2O(\grandc)}},e^{ \frac{\uc}{2O(\granda)}}, e^{-  \frac{\uc}{2O(\granda)}} \right].$$
So $\unit \in \ZZ\left[e^{ \frac{\ua}{2O(\grandc)}}, e^{-  \frac{\ua}{2O(\grandc)}},e^{ \frac{\uc}{2O(\granda)}}, e^{-  \frac{\uc}{2O(\granda)}} \right]$. 
Similarly $\unit \in \ZZ\left[e^{ \frac{\ub}{2O(\grandc)}}, e^{-  \frac{\ub}{2O(\grandc)}},e^{ \frac{\uc}{2O(\grandb)}}, e^{-  \frac{\uc}{2O(\grandb)}} \right]$. 
We conclude that $\unit \in \Lambda(\ExtH;\ua,\ub,\uc)$.
When $\Sigma$ is null-homologous, or when $|H_1(\rats)|$ is odd,  we can remove the $2$ from the above denominators, and $\unit \in \Lambda_o(\ExtH;\ua,\ub,\uc)$, thanks to Lemma~\ref{lemlambdaprimeJknot}.

Since $\CA_{\ExtH,3}$ vanishes when at least one of the variables $\ua$, $\ub$ or $\uc$ is zero, the lemma is proved.
\eop

\begin{remark}
We fix the normalization of $\CA_{\ExtH}$ of Lemma~\ref{lemnormAlexH}, but we keep in mind that the normalization breaks the order $3$ symmetry. Indeed, according to Lemma~\ref{lemuaubuc},
$$\CA_{\ExtH}(\tub \wedge \tuc)=\exp(\uc) \CA_{\ExtH}(\tua \wedge \tub)$$
and 
$\CA_{\ExtH}(\tuc \wedge \tua)=\exp(-\ub) \CA_{\ExtH}(\tua \wedge \tub).$
\end{remark}

 \begin{figure}[h]
\begin{center}
\begin{tikzpicture} 
\draw [thick, rounded corners, fill=blue!20] (3.4,2.2) -- (3.4,3.8) -- (0,3.8) -- (0,2.2)
(0,2.2) -- (.6,2.2) 
(.6,2.2) -- (.6,3) -- (1.4,2.4) -- (1.4,2.2)
(1.4,2.2) -- (2,2.2)
(2,2.2) -- (2,2.4) -- (2.8,3) -- (2.8,2.2)
(2.8,2.2) -- (3.4,2.2);
\draw [rounded corners, fill=yellow]
(0,1.6) -- (0,0) -- (3.4,0) -- (3.4,1.6) (3.4,1.6) -- (2.8,1.6)
(2.8,1.6) -- (2.8,.8) -- (2,1.4) -- (2,1.6) (2,1.6) -- (1.4,1.6)
(1.4,1.6) -- (1.4,1.4) -- (.6,.8) -- (.6,1.6) (.6,1.6) -- (0,1.6);
\draw [rounded corners, thick]
(0,1.6) -- (0,0) -- (3.4,0) -- (3.4,1.6)
(2.8,1.6) -- (2.8,.8) -- (2,1.4) -- (2,1.6)
(1.4,1.6) -- (1.4,1.4) -- (.6,.8) -- (.6,1.6);
\draw [rounded corners,->] (3.1,3.5) node{\scriptsize $\gamma$} (2.92,3.5) -- (2.92,3.68) -- (.48,3.68) -- (.48,2.2) (2.92,2.2) -- (2.92,3.5);
\draw [dashed, rounded corners,->] (2.6,.3) node{\scriptsize $\gamma$} (2.6,.12) --  (3.28,.12) -- (3.28,1.6) (.12,1.6) -- (.12,.12) -- (2.6,.12);
\draw [rounded corners,->] (1.1,3.3) node{\scriptsize $\alpha$} (1.2,3.12) -- (1.52,3.12) -- (1.88,2.84) -- (1.88,2.2) (.12,2.2) -- (.12,3.12) -- (1.2,3.12) ;
\draw [dashed, rounded corners,-<] (1.12,.78)  node{\scriptsize $\alpha$}  (1.12,.98) -- (1.52,1.28) -- (1.52,1.6) (.48,1.6) -- (.48,.68) -- (.72,.68) -- (1.12,.98) ;
\draw [rounded corners,->] (2.3,3.3) node{\scriptsize $\beta$} (2.2,3.12) -- (3.28,3.12) -- (3.28,2.2) (1.52,2.2) -- (1.52,2.84) -- (1.88,3.12) -- (2.2,3.12);
\draw [dashed, rounded corners,->] (2.15,.75)  node{\scriptsize $\beta$}  (2.28,.98) -- (1.88,1.28) -- (1.88,1.6) (2.92,1.6) -- (2.92,.68) -- (2.68,.68) -- (2.28,.98);
\draw [thick,->] (.8,2.85) -- (1,2.7);
\draw [thick,->] (2.2,2.55) -- (2.4,2.7);
\draw [thick,->] (2,0) -- (1.5,0);
\draw [thick,->] (2,3.8) -- (1.5,3.8);
\draw [thick,->] (2.2, 1.25) -- (2.4,1.1);
\draw [thick,->] (.93,2.75) node[below]{\scriptsize $K$} (.8,.95) -- (1,1.1);
\draw (.3,1.9) node{\scriptsize $\petitb$} (-.1,1.6)  rectangle (.7,2.2);
\draw (1.7,1.9) node{\scriptsize $\petitc$} (2.1,1.6) rectangle (1.3,2.2);
\draw (3.1,1.9) node{\scriptsize $\petita$} (3.5,1.6) rectangle (2.7,2.2);
\fill[color=white] (.1 ,3.25) rectangle (.4 ,3.65) ;
       \draw (.25 ,3.45) node{\scriptsize $\Sigma$}; 
\fill[color=white] (1.5 ,.3) rectangle (1.9 ,.7) ;
       \draw (1.7 ,.5) node{\scriptsize $-\Sigma$}; 
       \draw (1.4,3.12) -- (1.7,3.4) -- (2,3.12) (1.7,3.4) -- (1.7,3.68);
 \draw (1.85,3.5) node{\scriptsize $\ast$};      
\fill (1.7,3.4) circle (.05);
\begin{scope}[xshift=5cm]
 \draw [thick, rounded corners, fill=blue!20] (3.4,2.2) -- (3.4,3.8) -- (0,3.8) -- (0,2.2)
(0,2.2) -- (.6,2.2) 
(.6,2.2) -- (.6,3) -- (1.4,2.4) -- (1.4,2.2)
(1.4,2.2) -- (2,2.2)
(2,2.2) -- (2,2.4) -- (2.8,3) -- (2.8,2.2)
(2.8,2.2) -- (3.4,2.2);
\draw [rounded corners, fill=yellow]
(0,1.6) -- (0,0) -- (3.4,0) -- (3.4,1.6) (3.4,1.6) -- (2.8,1.6)
(2.8,1.6) -- (2.8,.8) -- (2,1.4) -- (2,1.6) (2,1.6) -- (1.4,1.6)
(1.4,1.6) -- (1.4,1.4) -- (.6,.8) -- (.6,1.6) (.6,1.6) -- (0,1.6);
\draw [rounded corners, thick]
(0,1.6) -- (0,0) -- (3.4,0) -- (3.4,1.6)
(2.8,1.6) -- (2.8,.8) -- (2,1.4) -- (2,1.6)
(1.4,1.6) -- (1.4,1.4) -- (.6,.8) -- (.6,1.6);
\draw [rounded corners,->] (2.55,3.55) node{\scriptsize $\gamma$} (2.62,3.4) -- (.48,3.4) -- (.48,2.2) (2.92,2.2) -- (2.92,3.4) -- (2.62,3.4)  ; 
\draw [dashed, rounded corners,->] (2.6,.3) node{\scriptsize $\gamma$} (2.6,.12) --  (3.28,.12) -- (3.28,1.6) (.12,1.6) -- (.12,.12) -- (2.6,.12);
\draw [rounded corners,->] (1,3.25) node{\scriptsize $\alpha$} (1.2,3.12) -- (1.52,3.12) -- (1.88,2.84) -- (1.88,2.2) (.12,2.2) -- (.12,3.12) -- (1.2,3.12) ;
\draw [dashed, rounded corners,-<] (1.12,.78)  node{\scriptsize $\alpha$}  (1.12,.98) -- (1.52,1.28) -- (1.52,1.6) (.48,1.6) -- (.48,.68) -- (.72,.68) -- (1.12,.98) ;
\draw [rounded corners,->] (2.25,2.9) node{\scriptsize $\beta$} (2.2,3.12) -- (3.28,3.12) -- (3.28,2.2) (1.52,2.2) -- (1.52,2.84) -- (1.88,3.12) -- (2.2,3.12);
\draw [dashed, rounded corners,->] (2.15,.75)  node{\scriptsize $\beta$}  (2.28,.98) -- (1.88,1.28) -- (1.88,1.6) (2.92,1.6) -- (2.92,.68) -- (2.68,.68) -- (2.28,.98);
\draw [thick,->] (.8,2.85) -- (1,2.7);
\draw [thick,->] (2.2,2.55) -- (2.4,2.7);
\draw [thick,->] (2,0) -- (1.5,0);
\draw [thick,->] (2,3.8) -- (1.5,3.8);
\draw [thick,->] (2.2, 1.25) -- (2.4,1.1);
\draw [thick,->] (.93,2.75) node[below]{\scriptsize $K$} (.8,.95) -- (1,1.1);
\draw (.3,1.9) node{\scriptsize $\petitb$} (-.1,1.6)  rectangle (.7,2.2);
\draw (1.7,1.9) node{\scriptsize $\petitc$} (2.1,1.6) rectangle (1.3,2.2);
\draw (3.1,1.9) node{\scriptsize $\petita$} (3.5,1.6) rectangle (2.7,2.2);
\fill[color=white] (.1 ,3.25) rectangle (.4 ,3.65) ;
       \draw (.25 ,3.45) node{\scriptsize $\Sigma$}; 
\fill[color=white] (1.5 ,.3) rectangle (1.9 ,.7) ;
       \draw (1.7 ,.5) node{\scriptsize $-\Sigma$}; 
       \draw (1.4,3.12) -- (1.7,3.4) -- (2,3.12) (2,3.8) -- (1.7,3.4);
 \draw (2.15,3.9) node{\scriptsize $\ast$};      
\fill (2,3.8) circle (.05);
\end{scope}
\begin{scope}[xshift=10cm]
\draw [fill=blue!20,draw opacity=0] (3.4,2.2) -- (3.4,3.8) -- (0,3.8) -- (0,2.2)
(0,2.2) -- (.6,2.2) 
(.6,2.2) -- (.6,3) -- (1.4,2.4) -- (1.4,2.2)
(1.4,2.2) -- (2,2.2)
(2,2.2) -- (2,2.4) -- (2.8,3) -- (2.8,2.2)
(2.8,2.2) -- (3.4,2.2)
(1.4,1.9) rectangle (2,2.2);
\draw [thick, rounded corners] (3.4,2.2) -- (3.4,3.8) -- (0,3.8) -- (0,2.2)
(0,2.2) -- (.6,2.2) 
(.6,2.2) -- (.6,3) -- (1.4,2.4) -- (1.4,1.9)
(1.4,1.9) -- (2,1.9)
(2,1.9) -- (2,2.4) -- (2.8,3) -- (2.8,2.2)
(2.8,2.2) -- (3.4,2.2);
\draw [rounded corners, fill=yellow]
(0,1.6) -- (0,0) -- (3.4,0) -- (3.4,1.6) (3.4,1.6) -- (2.8,1.6)
(2.8,1.6) -- (2.8,.8) -- (2,1.4) -- (2,1.6) (2,1.6) -- (1.4,1.6)
(1.4,1.6) -- (1.4,1.4) -- (.6,.8) -- (.6,1.6) (.6,1.6) -- (0,1.6);
\draw [rounded corners, thick]
(0,1.6) -- (0,0) -- (3.4,0) -- (3.4,1.6)
(2.8,1.6) -- (2.8,.8) -- (2,1.4) -- (2,1.6)
(1.4,1.6) -- (1.4,1.4) -- (.6,.8) -- (.6,1.6);
\draw [rounded corners,->] (2.62,3.4) -- (.48,3.4) -- (.48,2.2) (2.92,2.2) -- (2.92,3.4) -- (2.62,3.4)  ;
\draw [dashed, rounded corners,->] (2.6,.3) node{\scriptsize $\gamma$} (2.6,.12) --  (3.29,.12) -- (3.29,1.6) (.1,1.6) -- (.1,.12) -- (2.6,.12);
\draw [rounded corners,->] (1,3.25) node{\scriptsize $\alpha$} (1.2,3.12) -- (1.52,3.12) -- (1.88,2.84) -- (1.88,1.9) (.1,2.2) -- (.1,3.12) -- (1.2,3.12) ;
\draw [dashed, rounded corners,-<] (1.12,.78)  node{\scriptsize $\alpha$}  (1.12,.98) -- (1.52,1.28) -- (1.52,1.6) (.48,1.6) -- (.48,.68) -- (.72,.68) -- (1.12,.98) ;
\draw [rounded corners,->] (2.25,2.9) node{\scriptsize $\beta$} (2.2,3.12) -- (3.29,3.12) -- (3.29,2.2) (1.52,1.9) -- (1.52,2.84) -- (1.88,3.12) -- (2.2,3.12);
\draw [dashed, rounded corners,->] (2.15,.75)  node{\scriptsize $\beta$}  (2.28,.98) -- (1.88,1.28) -- (1.88,1.6) (2.92,1.6) -- (2.92,.68) -- (2.68,.68) -- (2.28,.98);
\draw [thick,->] (.8,2.85) -- (1,2.7);
\draw [thick,->] (2.2,2.55) -- (2.4,2.7);
\draw [thick,->] (2,0) -- (1.5,0);
\draw [thick,->] (2,3.8) -- (1.5,3.8);
\draw [thick,->] (2.2, 1.25) -- (2.4,1.1);
\draw [thick,->] (.93,2.75) node[below]{\scriptsize $K$} (.8,.95) -- (1,1.1);
\draw (.3,1.9) node{\scriptsize $\petitb$} (-.1,1.6)  rectangle (.7,2.2);
\draw (1.7,1.75) node{\scriptsize $\petitc$} (2.1,1.6) rectangle (1.3,1.9);
\draw (3.1,1.9) node{\scriptsize $\petita$} (3.5,1.6) rectangle (2.7,2.2);
\fill[color=white] (1.5 ,.3) rectangle (1.9 ,.7) ;
       \draw (1.7 ,.5) node{\scriptsize $-\Sigma$}; 
       \draw (1.4,3.12) -- (2,3.8) -- (2,3.12) (2,3.8) -- (1.85,3.4);   
 \draw (2.15,3.9) node{\scriptsize $\ast$};     
 \draw [rounded corners] (2,3.8) -- (.34,3.5) -- (.34,2.52)  (2,3.8) -- (3.14,3.5) -- (3.14,2.52) (2,3.8) --(1.74,3.4) -- (1.74,2.02); 
\fill (2,3.8) circle (.05);
\begin{scope}[yshift=.53cm]
 \draw [-<] (0,2.2) .. controls (0,2.12) and (.15,2) ..  (.3,2);
 \draw (.28,1.54) node[above]{\scriptsize $\ub$} (.3,2)  .. controls (.45,2)  and  (.6,2.12)  .. (.6,2.2);
\draw [dashed,->] (0,2.2) .. controls (0,2.28) and (.15,2.4) ..  (.3,2.4);
\draw [dashed] (.3,2.4)  .. controls (.45,2.4)  and  (.6,2.28)  .. (.6,2.2);
\end{scope}
\begin{scope}[xshift=1.4cm]
 \draw [-<] (0,2.2) .. controls (0,2.12) and (.15,2) ..  (.3,2);
 \draw (-.18,1.92) node[above]{\scriptsize $\uc$} (.3,2)  .. controls (.45,2)  and  (.6,2.12)  .. (.6,2.2);
\draw [dashed,->] (0,2.2) .. controls (0,2.28) and (.15,2.4) ..  (.3,2.4);
\draw [dashed] (.3,2.4)  .. controls (.45,2.4)  and  (.6,2.28)  .. (.6,2.2);
\end{scope}
\begin{scope}[xshift=2.8cm,yshift=.53cm]
 \draw [-<] (0,2.2) .. controls (0,2.12) and (.15,2) ..  (.3,2);
 \draw (.34,1.58) node[above]{\scriptsize $\ua$} (.3,2)  .. controls (.45,2)  and  (.6,2.12)  .. (.6,2.2);
\draw [dashed,->] (0,2.2) .. controls (0,2.28) and (.15,2.4) ..  (.3,2.4);
\draw [dashed] (.3,2.4)  .. controls (.45,2.4)  and  (.6,2.28)  .. (.6,2.2);
\end{scope}
\end{scope}

\end{tikzpicture}
\caption{The based curves $\alpha=\alpha(\petita,\petitb,\petitc)$ and $\beta=\beta(\petita,\petitb,\petitc)$ on the surface $\Sigma(\petita,\petitb,\petitc)$, with their basepoint moved to $K=K(\petita,\petitb,\petitc)$}
\label{figabstractSigmaalphabet}
\end{center}
\end{figure}
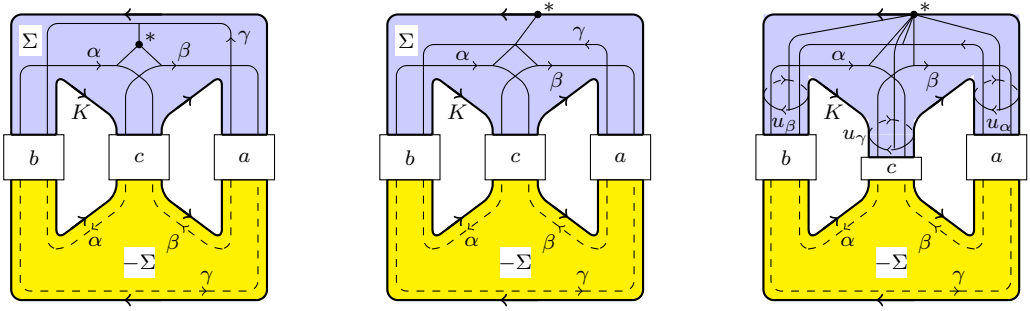

\begin{lemma} \label{lemexpralpha} Let $(\petita,\petitb,\petitc)$ be a triple of odd relative integers. Set $\amin=\frac{\petita-1}2$, $\apl=\frac{\petita+1}2$, $\bmin=\frac{\petitb-1}2$, $\bpl=\frac{\petitb+1}2$, $\cmin=\frac{\petitc-1}2$, and $\cpl=\frac{\petitb+1}2$. The homotopy classes of the based curves $\alpha(\petita,\petitb,\petitc)$ and $\beta(\petita,\petitb,\petitc)$ of Figure~\ref{figabstractSigmaalphabet} may be written as
 $$\alpha(\petita,\petitb,\petitc) =\uc^{(\petitc-1)/2}\granda \ub^{-(\petitb+1)/2} \;\;\; \mbox{and} \;\;\; \beta(\petita,\petitb,\petitc)=\ua^{(\petita-1)/2}\grandb \uc^{-(\petitc+1)/2}$$ in $\pi_1(\ExtH)$ with respect to the based curves of Figure~\ref{fighandlebobasedu}, and we have
$$\talpha(\petita,\petitb,\petitc)=\frac{\exp(\cmin\uc) -1}{\partial \tuc}\tuc
+\exp\left(\cmin\uc\right)\left(\tilde{\granda} + \exp(\granda)\frac{\exp(-\bpl\ub)-1}{\partial \tub}\tub\right)$$ and 
$$\tbeta(\petita,\petitb,\petitc)=\frac{\exp (\amin\ua) -1}{\partial \tua}\tua
+\exp\left(\amin\ua\right)\left(\tilde{\grandb} + \exp(\grandb)\frac{\exp(-\cpl\uc) -1}{\partial \tuc}\tuc\right)$$ in $\CH_{\ExtH}$.
\end{lemma}
\eop

\subsection{Proof of \texorpdfstring{Theorem~\ref{thmAlexHg}}{the structure theorem}}
\label{subsecpfstructA}

In this section, we prove Theorem~\ref{thmAlexHg} about the structure of the Alexander forms of rational homology genus two handlebodies, and we state two corollaries.

Theorem~\ref{thmAlexHg} is a direct consequence of the following lemmas~\ref{lemstructone} to \ref{lemstructfive}, in which we implicitly assume the hypotheses of Theorem~\ref{thmAlexHg} and use the normalization of $\CA_{\ExtH}$ from Lemma~\ref{lemnormAlexH}.

\begin{lemma} \label{lemstructone}
 We have $$\CA_{\ExtH}(\tub\wedge\tilde{\granda}) =\exp\left(-\ua - \grandb\right)  \CA_{\ExtH}(\tua \wedge \tilde{\grandb} ).$$
\end{lemma}
\bp 
Define $\cdel \in \CH_{\ExtH}$ to be $\cdel=\partial  \tilde{\granda} \tub - \partial \tub \tilde{\granda}$.
Then $\cdel=\exp\left(-\ua - \grandb\right) \bigl(\partial  \tilde{\grandb} \tua - \partial \tua \tilde{\grandb} \bigr)$ in $\CH_{\ExtH}$, according to Lemma~\ref{lemeqABCH}. Since $\partial \cdel=0$, Lemma~\ref{lemthreedel} implies  $\partial \tua \CA_{\ExtH}(\tub \wedge \cdel)= \partial \tub\CA_{\ExtH}(\tua\wedge \cdel)$.
So, we obtain $\partial \tua  \CA_{\ExtH}(\tub\wedge \partial \tub \tilde{\granda}) =\exp\left(-\ua - \grandb\right) \partial \tub \CA_{\ExtH}(\tua \wedge \partial \tua \tilde{\grandb}) $. \eop

\begin{lemma} \label{lemstructtwo} There exists  $v \in H_1(\ExtH;\QQ)$ 
such that $$\Delta(\granda,\grandb)=\exp\left(\frac{\ell-1}{2}\uc +v \right)\CA_{\ExtH}(\tub \wedge \tilde{\granda}).$$
If $\ell \neq 0$, then $v=0$. 

Let $(\petita,\petitb,\petitc)$ be a triple of odd integers. Set $\cmin =\frac{\petitc -1}{2}$.
Recall the curves $\alpha(\petita,\petitb,\petitc)$ and  $\beta(\petita,\petitb,\petitc)$ of Figure~\ref{figabstractSigmaalphabet} and Lemma~\ref{lemexpralpha}. If $\cmin \neq \ell$, then we have
$$\Delta(\alpha(\petita,\petitb,\petitc),\betam(\petita,\petitb,\petitc))=\exp\left(\frac{\ell}{2} \uc \right)\frac{\sih(\cmin\uc)}{\sih(\uc)}\CA_{\ExtH}(\tua \wedge\tub ) 
+ \exp\left(\frac{\cmin}{2}\uc -v \right)\Delta(\granda,\grandb).$$
\end{lemma}
\bp The first assertion follows from Example~\ref{exampleDeltaGrandaB} and Lemma~\ref{lemstructone}.
Let us compute the Alexander polynomial $\Delta(\alpha(\petita,\petitb,\petitc),\betam(\petita,\petitb,\petitc))$, for infinitely many triples $(\petita,\petitb,\petitc)$ of odd integers. 
Again, Example~\ref{exampleDeltaGrandaB} implies 
$$\Delta(\alpha(\petita,\petitb,\petitc),\betam(\petita,\petitb,\petitc)) \doteq  \CA_{\ExtH}(\tub \wedge\talpha(\petita,\petitb,\petitc) ).$$
The expression of $\talpha(\petita,\petitb,\petitc)$ in Lemma~\ref{lemexpralpha} implies 
$$\CA_{\ExtH}(\tub \wedge\talpha(\petita,\petitb,\petitc) )=\frac{\exp(\cmin\uc) -1}{\partial \tuc}\CA_{\ExtH}(\tub \wedge\tuc ) +\exp\left(\cmin\uc\right)\CA_{\ExtH}(\tub \wedge \tilde{\granda} ).
 $$
 In particular,  $\Delta(\alpha(\petita,\petitb,\petitc),\betam(\petita,\petitb,\petitc))$
depends only on the arbitrary integer $\cmin$. We denote it by $\Delta(\cmin)$. 
There exists $w(\cmin) \in H_1(\ExtH;\QQ)$ such that
$$\Delta(\cmin) 
=\exp\left(\frac{2\ell-c-1}{4}\uc + w(\cmin) \right)  \CA_{\ExtH}\bigl(\tub \wedge\talpha(\petita,\petitb,\petitc) \bigr).
$$
Then Lemma~\ref{lemuaubuc} implies
$$\Delta(\cmin)=\exp\bigl(w(\cmin)\bigr)\left(\exp\left(\frac{\ell}{2} \uc \right)\frac{\sih(\cmin\uc)}{\sih(\uc)}\CA_{\ExtH}(\tua \wedge\tub ) 
+ \exp\left(\frac{\cmin}{2}\uc -v \right)\Delta(\granda,\grandb)\right).$$
Since $\granda =\ell \ua +q\ub$ in $H_1(\ExtH;\QQ)$, for some $q \in \QQ$, the degree $0$ part $\CA_{\ExtH,0}(\tub \wedge \tilde{\granda})$ of $\CA_{\ExtH}(\tub \wedge \tilde{\granda})$ (and of $\Delta(\granda,\grandb)$) is $-|H_1(\rats)|\ell$. 
The degree $1$ part of $\Delta(\cmin)$, which must be zero, is
$$|H_1(\rats)| \left(\cmin -\ell\right)w(\cmin)+|H_1(\rats)|
\left(\cmin\frac{\ell}{2} \uc -\ell\bigl(\frac{\cmin}{2}\uc - v\bigr)\right)=|H_1(\rats)|\Bigl(
\left(\cmin -\ell\right)w(\cmin)+\ell v\Bigr).$$
In particular $\left(\cmin -\ell\right)w(\cmin)=-\ell v$. 

For any integer $\cmin$,
the degree $3$ part $\delta_3(\cmin)$ of $\frac{\Delta(\cmin)}{|H_1(\rats)|}$, which must be zero, may be expressed as follows.

$$\begin{array}{lll}\delta_3(\cmin)&=&\frac{1}{6}w(\cmin)^3\left(\cmin -\ell\right)
+ \frac{1}{2}w(\cmin)^2 \ell v\\
&&+ w(\cmin)\left(
\cmin\left(\frac{\ell^2 \uc^2 }{8}+\frac{\cmin^2 -1 }{24}\uc^2+\frac{\CA_{\ExtH,2}(\tua \wedge\tub )}{|H_1(\rats)|} \right) -
\frac{\ell}2 (\frac{\cmin}{2}\uc - v)^2+\frac{\Delta_2(\granda,\grandb)}{|H_1(\rats)|} \right) \\
&&+\cmin\left( \frac{\CA_{\ExtH,3}(\tua \wedge\tub )}{|H_1(\rats)|}
+\frac{\ell}{2}  \uc \left(\frac{\cmin^2 -1 }{24}\uc^2+\frac{\CA_{\ExtH,2}(\tua \wedge\tub )}{|H_1(\rats)|} \right)
+ \frac{\ell^3 \uc^3}{48}\right) \\ &&
- \frac{\ell}{6} (\frac{\cmin}{2}\uc - v)^3 + (\frac{\cmin}{2}\uc - v)\frac{\Delta_2(\granda,\grandb)}{|H_1(\rats)|}.
\end{array}$$
Multiply $\delta_3(\cmin)$ by $\left(\cmin -\ell\right)^2$ and
replace $\left(\cmin -\ell\right)w(\cmin)$ by $(-\ell v)$ in this equation to get
a polynomial in the variable $\cmin$. This polynomial is identically zero.
We may write it as
$\left(\cmin -\ell\right)^2\delta_3(\cmin)=\frac13 (\ell v)^3 + (\cmin -\ell)Q(\cmin)$ for some polynomial $Q$.
This implies that $\ell v=0$. Therefore, we have $v=0$ as soon as $\ell \neq 0$.

Furthermore, for any integer $\cmin$ such that $\cmin \neq \ell$, we have $w(\cmin)=0$. This proves the last assertion of the lemma.
\eop

\begin{lemma}\label{lemstructthree} If $\ell=0$ and $\Delta(\granda,\grandb)\neq 0$, then the element $v \in H_1(\ExtH;\QQ)$ of Lemma~\ref{lemstructtwo} is zero. 
\end{lemma}
\bp Let $\cmin$ be a nonzero integer.
Let $k \in \NN$ be the smallest integer such that $\Delta_{2k}(\granda,\grandb)\neq 0$.
For any $j<k$, we have $\CA_{\ExtH,2j+1}(\tua \wedge\tub )=0$. Indeed, if this holds for any $j^{\prime}<j$, then the null degree $(2j+1)$ part of $\Delta(\cmin)=\Delta(\alpha(\petita,\petitb,\petitc),\betam(\petita,\petitb,\petitc))$
is $\cmin \CA_{\ExtH,2j+1}(\tua \wedge\tub )$
according to Lemma~\ref{lemstructtwo}.
So, the null degree $(2k+1)$ part of $\Delta(\cmin)$ is
$$\cmin \left( \CA_{\ExtH,2k+1}(\tua \wedge\tub ) + \frac{1}{2}\uc \Delta_{2k}(\granda,\grandb)\right) - v \Delta_{2k}(\granda,\grandb)=0$$ for any nonzero integer $\cmin$. Thus $\CA_{\ExtH,2k+1}(\tua \wedge\tub ) =-\frac{1}{2}\uc \Delta_{2k}(\granda,\grandb)$ and $v=0$.
\eop

\begin{lemma}\label{lemstructfour} We have
 $$\sih(\uc)\Delta(\granda,\grandb) = \exp\left(-\frac{\ell}{2} \uc \right)\overline{\CA_{\ExtH}(\tua \wedge\tub )} - \exp\left(\frac{\ell}{2} \uc \right)\CA_{\ExtH}(\tua \wedge\tub ).$$
\end{lemma}
\bp According to Lemmas~\ref{lemstructtwo} and \ref{lemstructthree}, for any integer $\cmin \neq \ell$, we have
$$\Delta(\cmin)=\Delta\bigl(\alpha(\petita,\petitb,\petitc),\betam(\petita,\petitb,\petitc)\bigr)
=\exp\left(\frac{\ell}{2} \uc \right)\frac{\sih(\cmin\uc)}{\sih(\uc)}\CA_{\ExtH}(\tua \wedge\tub ) 
+ \exp\left(\frac{\cmin}{2}\uc\right)\Delta(\granda,\grandb).$$
Since $\Delta(\cmin)$ is even, we also have
$$\Delta(\cmin)
=\exp\left(-\frac{\ell}{2} \uc \right)\frac{\sih(\cmin\uc)}{\sih(\uc)}\overline{\CA_{\ExtH}(\tua \wedge\tub )}
+ \exp\left(-\frac{\cmin}{2}\uc\right)\Delta(\granda,\grandb).$$
\eop

\begin{lemma}\label{lemstructfive}
 If $\ell=0$, then the coefficient $\lambda^{\prime}(\ExtH;\ua,\ub,\uc)$ of 
Lemma~\ref{lemnormAlexH} and Theorem~\ref{thmAlexHg} is $\lambda^{\prime}(\granda,\grandb)$.
\end{lemma}
\bp 
In this case, the equations $$\begin{array}{ll}\sih(\uc)\Delta(\granda,\grandb)&=\overline{\CA_{\ExtH}(\tua \wedge\tub )} -\CA_{\ExtH}(\tua \wedge\tub )=-2\CA_{\ExtH,3}(\tua \wedge\tub )+O(4)\\&=-|H_1(\rats)|\lambda^{\prime}(\ExtH;\ua,\ub,\uc)\ua\ub\uc +O(4)\end{array}$$ imply $\Delta(\granda,\grandb)=-|H_1(\rats)|\lambda^{\prime}(\ExtH;\ua,\ub,\uc)\ua\ub +O(3)$,
where $\ua$ is a meridian of $\grandb$ and $(-\ub)$ is a meridian of $\granda$ (when $\granda$ and $\grandb$ are slightly pushed inside $\phi(\mhbz)$). So, Lemma~\ref{lemstructfive} follows from Definition~\ref{deflambdaprimetwocomp}, which is consistent with Definition~\ref{defzetatwo}, thanks to Lemma~\ref{lemzetatwocomp}.
\eop

The proof of Theorem~\ref{thmAlexHg} is finished. \eop

\begin{proposition}
 \label{propGamma}
 Let $\CA_{\ExtH,e}(\tua \wedge \tub)$ and $\CA_{\ExtH,o}(\tua \wedge \tub )$ respectively denote the even and the odd part of $\CA_{\ExtH,e}(\tua \wedge \tub )$.
Then there exists an even element $\Gamma=\Gamma(\ExtH;\ua,\ub,\uc)$ of 
$\frac{\QQ[[\ua,\ub,\uc]]}{\ua+\ub+\uc=0}$ such that 
$$\sih(\ua)\sih(\ub)\sih(\uc)\Gamma=\overline{\CA_{\ExtH}(\tua \wedge \tub )} - \CA_{\ExtH}(\tua \wedge \tub )=-2\CA_{\ExtH,o}(\tua \wedge \tub ),$$
whose degree $0$ part is
$\Gamma_0=-|H_1(\rats)|\lambda^{\prime}(\ExtH;\ua,\ub,\uc)$.
If $\granda$ and $\grandb$ are null-homologous, then $\Gamma \in \Lambda_{\ExtH}$. 
Set
$$\cosh(u)=\frac12\left(\exp \left(\frac{u}{2}\right) + \exp \left(-\frac{u}{2}\right)\right).$$
The multivariable Alexander polynomial $\Delta(\granda,\grandb)$ of $(\granda,\grandb)$ satisfies the following equation
$$\Delta(\granda,\grandb)=-\frac{\sih(\ell\uc)}{\sih(\uc)}\CA_{\ExtH,e}(\tua \wedge \tub ) + \sih(\ua)\sih(\ub)\cosh(\ell\uc)\Gamma.$$

\end{proposition}
\bp Lemma~\ref{lemnormAlexH} implies that $\left(\overline{\CA_{\ExtH}(\tua \wedge \tub )} - \CA_{\ExtH}(\tua \wedge \tub)\right)$ is a polynomial in $$\ZZ\left[\exp\Bigl( \frac{\ua}{2O(\grandb)}\Bigr), \exp\Bigl(-  \frac{\ua}{2O(\grandb)}\Bigr),\exp\Bigl( \frac{\ub}{2O(\granda)}\Bigr), \exp\Bigl(-  \frac{\ub}{2O(\granda)}\Bigr) \right].$$
According to Example~\ref{exampleDeltaGranda}, it is a multiple of $\bigl( \exp({\scriptstyle \frac{\ua}{2O(\grandb)}})-1\bigr)$.
It is similarly a multiple of $\bigl(\exp({\scriptstyle  \frac{\ub}{2O(\granda)}})-1\bigr)$ and $\bigl(\exp( {\scriptstyle \frac{\uc}{2O(\granda)}})-1\bigr)$. 

If $\granda$ and $\grandb$ are null-homologous, then $\bigl(\overline{\CA_{\ExtH}(\tua \wedge \tub )} - \CA_{\ExtH}(\tua \wedge \tub)\bigr)$ is a polynomial of  $\ZZ\left[e^{\ua}, e^{-\ua}, e^{\ub}, e^{-\ub} \right]$, and it is a multiple of
$\bigl(\left(e^{\ua} -1 \right)\left(e^{\ub} -1 \right)\left(e^{\uc} -1 \right)=\sih(\ua)\sih(\ub)\sih(\uc) \bigr)$. Use Lemma~\ref{lemstructfour} to conclude the proof.
\eop 

We can now state the following corollary of Theorem~\ref{thmAlexHg}.

\begin{corollary} \label{coreasytwo} Under the hypotheses of Theorem~\ref{thmAlexHg}, if $lk(\granda,\grandb)=0$, then $$\Delta(\granda,\grandb)=\sih(\ua)\sih(\ub)\Gamma,$$
with the series $\Gamma$ of Proposition~\ref{propGamma}.
Similarly, if $lk(\grandb,\grandc)=0$, then $\Delta(\grandb,\grandc)=\sih(\ub)\sih(\uc)\Gamma$,
and, if $lk(\grandc,\granda)=0$, then $\Delta(\grandc,\granda)=\sih(\uc)\sih(\ua)\Gamma$.

Assume that $\granda$ and $\grandb$ are rationally null-homologous in $\ExtH$.
Then we have $lk(\granda,\grandb)=lk(\grandb,\grandc)=lk(\grandc,\granda)=0$,
$ \partial \tub \tilde{\granda}=\exp\left(-\ua\right) \partial \tua \tilde{\grandb}$ in $\CH_{\ExtH}$, and, for any $x \in \CH_{\ExtH}$,
$$\begin{array}{rcl}\CA_{\ExtH}\bigl(x \wedge \tilde{\granda}\bigr)&=&\exp(\uc)\partial x \partial \tua \Gamma\\
\CA_{\ExtH}\bigl(x \wedge \tilde{\grandb}\bigr)&=&\exp(-\ub) \partial x \partial \tub  \Gamma\\  
\CA_{\ExtH}\bigl(x \wedge \tilde{\grandc}\bigr)&=&\partial x\partial \tuc \Gamma,\\
  \end{array}
$$ with the normalization of $\CA_{\ExtH}$ of Theorem~\ref{thmAlexHg}.
\end{corollary}
\noindent{\sc Proof of Corollary~\ref{coreasytwo}:}
When $\granda$ is zero in $H_1(\ExtH;\QQ)$, we have $\partial \tilde{\granda}=0$.
Then Lemma~\ref{lemthreedel} implies that $\CA_{\ExtH}(. \wedge \tilde{\granda} )$ is proportional to $\partial$. So, the expression of $\CA_{\ExtH}\bigl(\tub \wedge \tilde{\granda}\bigr)$ of Corollary~\ref{coreasyone} shows the first equality. 
According to Lemma~\ref{lemeqABCH}, we have $ \partial \tub \tilde{\granda}=\exp\left(-\ua\right) \partial \tua \tilde{\grandb}$. Therefore, the first equality implies the second one.
Since $\grandb\granda\grandc=1$ in $\pi_1(\ExtH)$, the third equality follows from Lemma~\ref{lemuaubuc}.
\eop

\begin{remark}
\label{rkDeltaAtilde} For any pair $(\granda^{\prime},\grandb^{\prime})$ of disjoint curves on $\partial \ExtH$ that are obtained from $\granda$ and $\grandb$ by Dehn twists about the curves 
$\ua$, $\ub$ and $\uc$ of Figure~\ref{fighandlebobasedu} (deprived of the paths to the basepoint), we have
$$\sih(\uc)\Delta(\granda^{\prime},\grandb^{\prime}) = \exp\left(-\frac{lk(\granda^{\prime},\grandb^{\prime})}{2} \uc \right)\overline{\CA_{\ExtH}(\tua \wedge \tub )} - \exp\left(\frac{lk(\granda^{\prime},\grandb^{\prime})}{2} \uc \right)\CA_{\ExtH}(\tua \wedge \tub )$$ as in Theorem~\ref{thmAlexHg}.
The expressions for $\CA_{\ExtH}(\tub \wedge \granda^{\prime})$ and $\CA_{\ExtH}(\tua \wedge \grandb^{\prime})$ are also valid for such curves. (Note that these curves 
coincide with $\granda$ and $\grandb$ on the top half of Figure~\ref{fighandlebobasedu}, which contains all the basepoints.) 
Let us move the basepoint on $\partial \Sigma(\petita,\petitb,\petitc)$ as in Figure~\ref{figabstractSigmaalphabet} without loss. So, the curves $\alpham$ and $\betam$  are based curves of $\Sigma_-$ obtained from $\alpha$ and $\beta$ by the natural diffeomorphism from $\Sigma$ to  $\Sigma_-$.
For any triple $(\petita,\petitb,\petitc)$, $$\bigl(\alpha(\petita,\petitb,\petitc),\cbetam(\petita,\petitb,\petitc)=\ua \betam \ua^{-1}\bigr)\;\;\mbox{and}\;\; \bigl(\calpham(\petita,\petitb,\petitc)=\ub^{-1}\alpham\ub,\beta(\petita,\petitb,\petitc)\bigr)$$ are examples of pairs $(\granda^{\prime},\grandb^{\prime})$ as above.
\end{remark}

The \emph{one-variable Alexander polynomial} of the link $(\granda,\grandb)$ whose component meridians are $m(\granda)$ and $m(\grandb)$, is obtained from the two-variable  Alexander polynomial, by setting $$t^{\frac{1}{2O(\granda)O(\grandb)}}= \exp\left(\frac{1}{2O(\granda)O(\grandb)}m(\granda)\right)= \exp\left(\frac{1}{2O(\granda)O(\grandb)}m(\grandb)\right).$$
The following corollary of Theorem~\ref{thmAlexHg} and Proposition~\ref{propGamma} allows us to compute
$$\CA_{\ExtH,3}(\tua \wedge\tub )=\frac{|H_1(\rats)|}{2}\lambda^{\prime}(\ExtH;\ua,\ub,\uc)\ua\ub\uc$$ from the one-variable Alexander polynomials of 
$\granda$, $\grandb$, $\grandc$, and $(\granda,\grandb)$. (Here, $(\granda,\grandb)$ can be replaced with any $(\granda^{\prime},\grandb^{\prime})$ as in Remark~\ref{rkDeltaAtilde}.)

\begin{corollary} \label{corlambdaprimethree} Recall $\ell=lk(\granda,\grandb)$. The degree two part $\Delta_2(\granda,\grandb)$ of $\Delta(\granda,\grandb)$ is
 $$\Delta_2(\granda,\grandb)=|H_1(\rats)|\ell\frac{1-\ell^2}{24}\uc^2 - \ell\CA_{\ExtH,2}(\tua \wedge\tub) - |H_1(\rats)|\lambda^{\prime}(\ExtH;\ua,\ub,\uc)\ua\ub,$$
and we have
$$2\CA_{\ExtH,3}(\tua \wedge\tub )=|H_1(\rats)|\ell\frac{1-\ell^2}{24}\uc^3 -\uc \bigl(\ell\CA_{\ExtH,2}(\tua \wedge\tub)+\Delta_2(\granda,\grandb)\bigr).$$
\end{corollary}
\bp 
Proposition~\ref{propGamma} implies
$\Delta(\granda,\grandb)=-\frac{\ssih(\ell\uc)}{\ssih(\uc)}\CA_{\ExtH,e}(\tua \wedge \tub ) + \sih(\ua)\sih(\ub)\cosh(\ell\uc)\Gamma.$
The first equation follows easily. The second equation follows from the first one.
\eop

\section{The Alexander series of a genus one Seifert surface}
\setcounter{equation}{0}
In Section~\ref{subintrocd}, we associated a Reidemeister torsion $\tau(\ExtH[K])$ to a genus one Seifert surface $\Sigma$. 
In this section, we study this torsion. We view it as a series $\CD(\Sigma)$, which
is described in Propositions~\ref{propgenDSigma} and \ref{propcompDEK}.

\subsection{Formulas for the Alexander series of a genus one Seifert surface}
\label{subdescDsigma}

In this subsection, we first prove the following proposition.

\begin{proposition}
\label{propgenDSigma}
Let $\phi \colon \hbz \hookrightarrow \rats$ be an embedding of the handlebody $\hbz$ in a $\QQ$-sphere $\rats$. Let $(\check{a}, \check{b}, \check{c})$  be a triple of odd integers. Set
$\Sigma=\phi( \Sigma(\check{a}, \check{b}, \check{c}) )$.
Let $\alpha=\phi(\alpha(\check{a}, \check{b}, \check{c}))$ and $\beta=\phi(\beta(\check{a}, \check{b}, \check{c}))$ be the corresponding curves of $\Sigma$, as in Figure~\ref{figSigmaabc}.
 Define $\amin=-lk(\beta,\gammam)$, $\bmin=-lk(\gamma,\alpham)$,  $\cmin=-lk(\alpha,\betam)$, $\apl=-lk(\betam,\gamma)$, $\bpl=-lk(\gammam,\alpha)$,  $\cpl=-lk(\alpham,\beta)$,
 $\petita=\apl + \amin$, $\petitb=\bpl + \bmin$, and $\petitc=\cpl + \cmin$.
 We have $\cmin=\frac{\petitc-1}2$ and $\cpl=\frac{\petitc+1}2$.
Set $$\begin{array}{llll}\CE(\petita,\petitb,\petitc)&=&&\exp\bigl(\frac{\amin}{2}\ua\bigr)
\sih(\ub)\sih\left(\alpha= \cmin\uc-\bpl \ub\right)\\
&&-&\exp\bigl(\frac{\bpl}{2}\ub\bigr)
\sih(\ua)\sih\left(\beta= \amin\ua-\cpl\uc\right).\end{array}$$
Then we have $$\begin{array}{llll}\CE(\petita,\petitb,\petitc)
&=&&\exp \left(\frac{\amin \ua +\bmin \ub -\cmin \uc}{2}\right)-\exp\left(\frac{\apl \ua +\bpl \ub - \cpl \uc}{2} \right)\\
&&+&\exp \left(\frac{\bmin \ub +\cmin \uc -\amin \ua}{2}\right)-\exp\left(\frac{\bpl \ub +\cpl \uc- \apl \ua}{2}\right)\\
&&+&\exp \left(\frac{\cmin \uc +\amin \ua -\bmin \ub}{2}\right)-\exp\left(\frac{\cpl \uc +\apl \ua- \bpl \ub}{2}\right).\end{array}$$
Furthermore, with the normalization of $\CA_{\ExtH}$ of Theorem~\ref{thmAlexHg}, we have $$\sih(\ua)\sih(\ub)\sih(\uc)\CD(\Sigma) = \CE(\petita,\petitb,\petitc)\overline{\CA_{\ExtH}(\tua \wedge \tub)}- \overline{\CE(\petita,\petitb,\petitc)}\CA_{\ExtH}(\tua \wedge \tub).
  $$
and
$$\CD(\Sigma) =\exp\left(\frac12 \bigl(-\apl \ua +\bpl \ub - \cpl \uc\bigr)\right) \CA_{\ExtH}(\talpha \wedge \tbeta).$$
\end{proposition}

\begin{remark}
When $\rats$ is a $\ZZ$-sphere, we can deduce that $\CD(\Sigma)$ is in $\Lambda_{\ExtH}$ if and only if $(\petita-\petitb)$ and $(\petita-\petitc)$ are in $4\ZZ$.
\end{remark}

In this subsection, we work under the hypotheses of Proposition~\ref{propgenDSigma}. Note that any genus one Seifert surface $\Sigma$ may be written as $\phi(\Sigma(\check{a}, \check{b}, \check{c}))$
for some (and actually any) triple $(\check{a}, \check{b}, \check{c})$  of odd integers. The exterior of $\Sigma$ is $\ExtH= \rats \setminus \phi\bigl(\mhbz\bigr)$.

We use the normalization of the Alexander form $\CA_{\ExtH}$ of Theorem~\ref{thmAlexHg}. 
The following lemma gives a first expression of the Reidemeister torsion $\tau(\ExtH[K])$ of $\ExtH[K]$ introduced in Section~\ref{subintrocd}.

\begin{lemma} \label{lemcompapalone} We have
 $$\tau(\ExtH[K]) \doteq \CA_{\ExtH}\left(\talpha \wedge\tbeta\right).$$
\end{lemma}
\bp According to our choice of sign in Section~\ref{subintrocd}, we have $\tau(\ExtH[K])\doteq \frac{\CA_{\ExtH}(\tilde{K} \wedge x)}{\partial x}$ for any $x \in \CH_{\ExtH}$
such that $\partial x \neq 0$.
Since $K=\alpha^{-1}\beta\alpha\beta^{-1}=\alpha^{-1}\left(\beta\alpha\beta^{-1}\alpha^{-1}\right) \alpha$,
we have $\tilde{K}= \exp(-\alpha)\bigl((\partial \tbeta)\talpha -(\partial \talpha) \tbeta\bigr)$ in $\CH_{\ExtH}$.
Therefore, we get the announced result since $\alpha$ and $\beta$ cannot be both rationally null-homologous in $\ExtH$ (because they intersect homologically).
\eop
 
\begin{lemma}
\label{lemAEalphab} We have
 $$\partial \tua \partial \tub \CA_{\ExtH}(\talpha \wedge \tbeta)
 =\partial \talpha \partial \tub \CA_{\ExtH}(\tua \wedge \tbeta)
 -\partial \tua \partial \tbeta \CA_{\ExtH}(\tub \wedge \talpha)
-\partial \talpha \partial \tbeta  \CA_{\ExtH}(\tua \wedge \tub).$$
\end{lemma}
\bp 
Lemma~\ref{lemthreedel} implies
$$ \partial \tub \CA_{\ExtH}(\talpha \wedge \tbeta) = \partial \talpha \CA_{\ExtH}(\tub \wedge \tbeta) - \partial \tbeta  \CA_{\ExtH}(\tub \wedge \talpha)$$
and
$$ \partial \tua \CA_{\ExtH}(\tub \wedge \tbeta) = \partial \tub \CA_{\ExtH}(\tua \wedge \tbeta) - \partial \tbeta  \CA_{\ExtH}(\tua \wedge \tub).$$
\eop

\noindent{\sc Proof of Proposition~\ref{propgenDSigma}:}
View $(-\ub)$ as a meridian of $\alpham$ (pushed inside $\phi(\mhbz)$) and $\ua$ as a meridian of $\betam$. Thus
$$\alpha=-lk(\alpha,\alpham)\ub +lk(\alpha,\betam)\ua= lk(\alpha,\betam)(\ua +\ub) + lk(\alpha,\gammam)\ub= \cmin\uc-\bpl \ub$$
in $H_1(\ExtH;\QQ)$ (as in Lemma~\ref{lemexpralpha}, when $\granda$ and $\grandb$ vanish in $H_1(\ExtH;\QQ)$). To prove that $\CE(\petita,\petitb,\petitc)$ has the given symmetric six-term expression, we expand its defining expression and group the terms that contribute to $$\exp\left(\frac12\left(\varepsilon_a\amin \ua +\varepsilon_b\bpl \ub +\varepsilon_c (\cmin \,\mbox{or}\, \cpl) \uc\right) + \mbox{monomial independent of $\petita$, $\petitb$, $\petitc$} \right)$$ for every $(\varepsilon_a,\varepsilon_b,\varepsilon_c)$ in 
$\{-1,1\}^3$.

Let $\CR=\CE(\petita,\petitb,\petitc)\overline{\CA_{\ExtH}(\tua \wedge \tub)}- \overline{\CE(\petita,\petitb,\petitc)}\CA_{\ExtH}(\tua \wedge \tub)$ be the right-hand side of the third equation of the statement.
Set $$d=\exp\left(\frac12 \left(-\apl \ua +\bpl \ub - \cpl \uc \right)\right) \sih(\ua)\sih(\ub)\sih(\uc)\CA_{\ExtH}(\talpha \wedge \tbeta).$$
Recall Lemma~\ref{lemcompapalone}. Since $\CR$ is odd, it suffices 
to prove that  
$d=\CR$ to finish the proof of Proposition~\ref{propgenDSigma}. 
According to Remark~\ref{rkDeltaAtilde} and Corollary~\ref{coreasyone}, since $lk(\alpha,\betam)=-\cmin$, we have
$$\sih(\uc)\CA_{\ExtH}(\tub \wedge \talpha)=\exp\left(\frac{1 +2\cmin}2{\uc}\right) \overline{\CA_{\ExtH}(\tua \wedge \tub)} -\exp\left(\frac{\uc}2\right)\CA_{\ExtH}(\tua \wedge \tub).$$
Similarly, $lk(\alpham,\beta)=-\cpl$. So, Remark~\ref{rkDeltaAtilde} and Theorem~\ref{thmAlexHg} imply
$$\sih(\uc)\CA_{\ExtH}(\tua \wedge \tbeta)= \exp\left(\frac{1 +2\cpl}2{\uc} +\ua+\beta\right) \overline{\CA_{\ExtH}(\tua \wedge \tub)} -\exp\left(\frac{\uc}2+\ua+\beta\right)\CA_{\ExtH}(\tua \wedge \tub).$$
We have $d=\exp\left(\frac12 \left(-\cmin \uc +\ub -\alpha -\beta\right)\right) \partial \tua\partial \tub\sih(\uc)\CA_{\ExtH}(\talpha \wedge \tbeta)$, and Lemma~\ref{lemAEalphab} implies
$d
=e \overline{\CA_{\ExtH}(\tua \wedge \tub)} +e^{\prime} \CA_{\ExtH}(\tua \wedge \tub)$
with 
$$\begin{array}{lll}e&=&\sih(\alpha)\sih(\ub)
\exp\left(\ua+\ub+ (1+2\cpl-\cmin)\frac{\uc}2 +\frac12\beta\right)\\
&&-\sih(\ua) \sih(\beta)\exp\left(\frac12 \left(\ua+\ub +(1 +\cmin)\uc -\alpha\right)\right)\\&=&\sih(\alpha)\sih(\ub)\exp\left(\frac{\amin}2\ua\right)-\sih(\ua) \sih(\beta)\exp\left(\frac{\bpl}2\ub\right)
\\&=&
\CE(\petita,\petitb,\petitc) \end{array}
$$
 and 
 $$\begin{array}{lll}e^{\prime}
&=&
 -\sih(\alpha= \cmin\uc-\bpl \ub)\sih(\ub)
\exp\left(\frac12 \left(\amin \ua - 2\cpl \uc\right)\right)\\
&&+\sih(\ua) \sih(\beta= \amin\ua-\cpl\uc)
\exp\left(\frac12 \left(\bpl \ub -2\cmin \uc\right)\right)
\\
&&-\sih(\uc)\sih( \alpha) \sih( \beta)\exp\left(\frac12 \left(-\cmin \uc +\ub\right)\right)\\
&=&-\exp\left(\frac12 \left( \apl \ua +(1-\bpl) \ub - \cpl \uc \right)\right)\sih(\ub)\\
&&+\exp\left(\frac12 \left( \amin \ua + \bpl \ub -(2\cpl +\cmin) \uc  \right)\right)\sih(\ub)\\
&&+\exp\left(\frac12 \left(  \amin \ua + \bpl \ub -(2\cmin +\cpl) \uc \right)\right)\sih(\ua)\\
&&-\exp\left(\frac12 \left( - \amin \ua +\bpl \ub +(1 -\cmin) \uc \right)\right)\sih(\ua)\\
&&-\exp\left(\frac12 \left(  \amin \ua +(1-\bpl) \ub -\cpl \uc \right)\right)\sih(\uc)\\
&&-\exp\left(\frac12 \left( -\apl \ua + \bpl \ub -\cmin \uc  \right)\right)\sih(\uc)\\
&&+\exp\left(\frac12 \left( -\amin \ua +(1-\bpl) \ub +\cpl \uc  \right)\right)\sih(\uc)\\
&&+\exp\left(\frac12 \left(  \amin \ua +(\bpl +1) \ub -(2\cmin +\cpl) \uc \right)\right)\sih(\uc).\end{array}$$
The sum of the second, third, and eighth lines is
$$\exp\left(\frac12 \left(  \amin \ua +\bpl \ub -(2\cmin +\cpl) \uc \right)\right)\left(\exp\left(-\frac{\uc}{2}\right)\sih(\ub) + \sih(\ua) + \exp\left(\frac{\ub}{2}\right)\sih(\uc)\right).$$
It vanishes thanks to Lemma~\ref{lemuaubuc}.
The sum of the other lines is 
$$\begin{array}{llll}e^{\prime}&=&&\exp\left(\frac{\apl \ua -\bpl \ub - \cpl \uc}{2} \right)-\exp \left(\frac{\amin \ua -\bmin \ub -\cmin \uc}{2}\right)\\
&&+&\exp\left(\frac{\bpl \ub -\cpl \uc- \apl \ua}{2}\right)-\exp \left(\frac{\bmin \ub -\cmin \uc -\amin \ua}{2}\right)\\
&&+&\exp\left(\frac{\cpl \uc -\apl \ua- \bpl \ub}{2}\right)-\exp \left(\frac{\cmin \uc -\amin \ua -\bmin \ub}{2}\right)
\\
&=&-&\overline{\CE(\petita,\petitb,\petitc)}. \end{array}$$
\eop
\begin{proposition}
 \label{propDabcbis}
 
Recall the surface $\Sigma(\petita,\petitb,\petitc)$ of Figure~\ref{figSigmaabc} embedded in $S^3$.
The Alexander series $\CD(\petita,\petitb,\petitc)=\CD(\Sigma(\petita,\petitb,\petitc))$
is
\begin{multline*}\CD(\petita,\petitb,\petitc)=\exp\left(\frac{\amin}2\ua\right)\frac{\sih (\bpl\ub)}{\sih(\ub)} \frac{ \sih(\cpl \uc)}{\sih(\uc)} +  
\exp\left(\frac{\bpl}2\ub\right)\frac{\sih (\cmin\uc)}{\sih( \uc)} \frac{\sih (\amin\ua)}{\sih(\ua)}  \\
+\exp\left(\frac{\cpl}2\uc + \frac12 \ub\right)\frac{\sih (\amin\ua)}{\sih( \ua)}\frac{\sih (\bpl\ub)}{\sih(\ub)}. \end{multline*}
It satisfies the following equalities:
$$\begin{array}{llll}\sih (\ua)\sih (\ub)\sih (\uc)\CD(\petita,\petitb,\petitc)&=&&\CE(\petita,\petitb,\petitc)-\overline{\CE(\petita,\petitb,\petitc)}\\
&=&&\sih \left(\amin \ua +\bmin \ub -\cmin \uc\right)-\sih\left(\apl \ua +\bpl \ub - \cpl \uc \right)\\
&&+&\sih \left(\bmin \ub +\cmin \uc -\amin \ua\right)-\sih\left(\bpl \ub +\cpl \uc- \apl \ua\right)\\
&&+&\sih \left(\cmin \uc +\amin \ua -\bmin \ub\right)-\sih\left(\cpl \uc +\apl \ua- \bpl \ub\right).\end{array}$$
  \end{proposition}
  \bp Recall that $\CA_{\ExtH_0}(\tua \wedge \tub)=1$ from Example~\ref{exaAlextriv}. Then Proposition~\ref{propgenDSigma} implies that 
$\sih (\ua)\sih (\ub)\sih (\uc)\CD(\Sigma(\petita,\petitb,\petitc))=\CE(\petita,\petitb,\petitc)-\overline{\CE(\petita,\petitb,\petitc)}$.
So, it suffices to check that the first given expression $\tilde{\CD}$ of $\CD(\petita,\petitb,\petitc)$ satisfies this equation. To check it, we compute
$$\begin{array}{lll}\sih(\ua)\sih(\ub)\sih(\uc)\tilde{\CD}=&& \exp\left(\frac12\amin\ua\right)\sih(\ua){\sih (\bpl\ub)}{ \sih(\cpl \uc)}\\&&  
+
\exp\left(\frac12\bpl\ub\right) {\sih (\amin\ua)}\sih(\ub){\sih (\cmin\uc)}\\&&  
+ \exp\left(\frac12 \cpl\uc + \frac12 \ub\right){\sih (\amin\ua)}{\sih (\bpl\ub)}\sih(\uc).\end{array}$$
Again, we expand this expression 
and group the terms
that contribute to $$\exp\left(\frac12\left(\varepsilon_a\amin \ua +\varepsilon_b\bpl \ub +\varepsilon_c (\cmin \,\mbox{or}\, \cpl) \uc\right) + \mbox{monomial independent of $\petita$, $\petitb$, $\petitc$} \right)$$ for every $(\varepsilon_a,\varepsilon_b,\varepsilon_c)$ in 
$\{-1,1\}^3$. There are no terms for which $\varepsilon_a=\varepsilon_b=\varepsilon_c=-1$. The six terms for which $\varepsilon_a=\varepsilon_b=\varepsilon_c=1$ cancel. The six terms for which exactly one $\varepsilon_.$ is $1$ correspond directly to terms in the right-hand side of the last equality of the statement.
\eop

The following proposition is an immediate corollary of Propositions~\ref{propgenDSigma} and \ref{propDabcbis}.

\begin{proposition} \label{propcompDEK}
Under the hypotheses of Proposition~\ref{propgenDSigma}, we have
 $$\CD(\Sigma)=\CD(\petita,\petitb,\petitc)\CA_{\ExtH}\left(\tua\wedge\tub\right) + \CE(\petita,\petitb,\petitc) \Gamma$$
where 
$\CD(\petita,\petitb,\petitc)$
is the even polynomial of Proposition~\ref{propDabcbis},
 $\CE(\petita,\petitb,\petitc)$ is defined in Proposition~\ref{propgenDSigma}, and $\Gamma= \Gamma(\ExtH;\ua,\ub,\uc)$ is defined in Proposition~\ref{propGamma}.
\end{proposition}
\bp Proposition~\ref{propgenDSigma} implies
$$\begin{array}{ll}\sih(\ua)\sih(\ub)\sih(\uc)\CD(\Sigma) =& \left(\CE(\petita,\petitb,\petitc)- \overline{\CE(\petita,\petitb,\petitc)}\right)\CA_{\ExtH}(\tua \wedge \tub)\\&
+ \CE(\petita,\petitb,\petitc) \left(\overline{\CA_{\ExtH}(\tua \wedge \tub)}- \CA_{\ExtH}(\tua \wedge \tub)\right).\end{array}$$ \eop

 \subsection{The low degree terms of \texorpdfstring{$\CD(\Sigma)$}{the Alexander series of a genus one Seifert surface}}
 \label{secAlexEKlow}
 
 In this section, we compute $\CD(\Sigma)$ up to order two under the hypotheses of Proposition~\ref{propgenDSigma} to prove Theorem~\ref{thmalexcanphi}.
The vector space $$H_1(\Sigma;\QQ) \otimes_s H_1(\Sigma;\QQ)=\QQ (\alpha \otimes_s \alpha) \oplus \QQ( \beta \otimes_s \beta) \oplus \QQ (\gamma \otimes_s \gamma)$$
is equipped with the canonical symmetric non-degenerate bilinear form
$$\begin{array}{llll}\CB \colon &\left(H_1(\Sigma;\QQ) \otimes_s H_1
(\Sigma;\QQ)\right)^2&\rightarrow &\QQ\\
&(u \otimes_s v, x \otimes_s y) &\mapsto& \langle u,x \rangle_{\Sigma}\langle 
v,y \rangle_{\Sigma} +\langle u,y \rangle_{\Sigma}\langle v,x 
\rangle_{\Sigma},\end{array}$$
which maps $(\alpha \otimes_s \alpha, \beta \otimes_s \beta)$ to $2$.
The tensor $$V_s^{\ast}=\frac{a}{2}  \alpha \otimes_s \alpha + \frac{b}{2} \beta \otimes_s \beta +\frac{c}
{2} \gamma \otimes_s \gamma$$ is the dual of the symmetrized Seifert form $V_s$ introduced in Section~\ref{subintrocd}, with respect to $\CB$.
$V_s(u \otimes_s v)=\CB(u \otimes_s v,V_s^{\ast})$.
So $V_s^{\ast}$ is a canonical element of $H_1(\Sigma;\QQ) \otimes_s H_1(\Sigma;\QQ)$.

\begin{lemma}\label{lemCEdegtwo} Recall the isomorphism $\fHone$ of Proposition~\ref{propfHone}.
Set $$\Pcanp=(\fHone \otimes_s \fHone)(V_s^{\ast})=\frac{a}{2}  \ua^2 +\frac{b}{2} 
\ub^2 +\frac{c}{2}\uc^2.$$
Then $\Pcanp$ is a canonical element 
of $H_1(\ExtH;\QQ) \otimes_s H_1(\ExtH;\QQ)$. Furthermore, we have $$\Pcanp = \ua\otimes_s\beta -\ub\otimes_s\alpha\;\;\;\;\mbox{and}\;\;\;\;\CE(\petita,\petitb,\petitc)= -\Pcanp  +O(3).$$
\end{lemma}
\bp We have $$\ua\otimes_s(\beta= \amin\ua-\cpl\uc) - \ub\otimes_s(\alpha=\cmin\uc-\bpl \ub)=\Pcanp +\frac12\left(-\ua^2-\uc\ua+\ub\uc + \ub^2\right)=\Pcanp.$$
According to Proposition~\ref{propgenDSigma}, we have $\CE(\petita,\petitb,\petitc)= \ub\otimes_s\alpha -\ua\otimes_s\beta  +O(3).$ \eop

\begin{lemma} \label{lemCDdegzero} 
We have $$\CD(\petita,\petitb,\petitc)=\lambda^{\prime}(\petita,\petitb,\petitc) +O(2)$$
and the degree $0$ part of $\CD(\Sigma)$ is $|H_1(\rats)|\lambda^{\prime}(K)$.
\end{lemma}
\bp The first expression of $\CD(\petita,\petitb,\petitc)$ in Proposition~\ref{propDabcbis} shows that
$$\varepsilon(\CD(\petita,\petitb,\petitc))=\bpl\cpl+\cmin\amin+\amin\bpl=\lambda^{\prime}(\petita,\petitb,\petitc)$$
while the second one shows 
that $\CD(\petita,\petitb,\petitc)$ is even.
Then Proposition~\ref{propcompDEK}, Theorem~\ref{thmAlexHg}, and Lemma~\ref{lemCEdegtwo} 
imply that the degree $0$ part of $\CD(\Sigma)$ is $|H_1(\rats)|\lambda^{\prime}(K)$.
\eop

Note the following easy lemma.
\begin{lemma} \label{lemdegtwoCDprel} Recall the polynomial invariant $\Pcanp$ from Lemma~\ref{lemCEdegtwo}, and the polynomial $\ddel_{2,\Sigma}(\alpha,\beta)$ of Theorem~\ref{thmsecondinvariantbis}.
 Set $\ddel_{2,\ExtH}(\alpha,\beta)=\fHone_{\ast}(\ddel_{2,\Sigma}(\alpha,\beta))$.
Then we have $$\begin{array}{lll}\ddel_{2,\ExtH}(\alpha,\beta)
&=&-\frac{1}{24} \left((\petita^2+2\petita\petitb+2\petita\petitc+3)\ua^2 +(\petitb^2+2\petitb\petitc+2\petitb\petita+3)\ub^2 +(\petitc^2+2\petitc\petita+2\petitc\petitb+3)\uc^2 \right)\\
 &=&\frac{a^2\ua^2 +b^2\ub^2 +c^2\uc^2}{24} - 
\frac{a+b+c}{6}\Pcanp
 -\frac{\ua^2 +\ub^2 +\uc^2}{8}\\
 &=&\frac{1}{24}\left((\petita+\petitb)^2 +8\lambda^{\prime}(\petita,\petitb,\petitc)+4\right)\ua\ub +\frac{1}{24}\left((\petitb+\petitc)^2 +8\lambda^{\prime}(\petita,\petitb,\petitc)+4\right)\ub\uc \\
 &&+\frac{1}{24} \left((\petitc+\petita)^2 +8\lambda^{\prime}(\petita,\petitb,\petitc)+4\right)\uc\ua
 .\end{array}$$
\end{lemma}
\eop

\begin{lemma} \label{lemdegtwoCD} With the notation of Lemma~\ref{lemdegtwoCDprel},
the degree two part $\CD_2(\petita,\petitb,\petitc)$ of the polynomial $\CD(\petita,\petitb,\petitc)$ of Proposition~\ref{propDabcbis} is 
$$\begin{array}{lll}\CD_2(\petita,\petitb,\petitc) 
   &=&\frac{(\petita^2-1)}{16 \times 24}(3+3\petitb\petitc+\petita(\petitb+\petitc))\ua^2
   + \frac{(\petitb^2-1)}{16 \times 24}(3+3\petitc\petita+\petitb(\petitc+\petita))\ub^2\\
   &&+\frac{(\petitc^2-1)}{16 \times 24}(3+3\petita\petitb+\petitc(\petita+\petitb))\uc^2\\
  &=& - \frac1{12} \wdel(\petita,\petitb,\petitc)\Pcanp + \frac{\lambda^{\prime}(\petita,\petitb,\petitc)}{4}\ddel_{2,\ExtH}(\alpha,\beta) 
   ,
  \end{array}
$$
where  $\wdel(\petita,\petitb,\petitc)$ is defined in Theorem~\ref{thmsimpleinvt}.
\end{lemma}
\bp Recall
$$\begin{array}{lll}\CD(\petita,\petitb,\petitc)&=&\exp\left(\frac12\amin\ua\right)\frac{\sih (\bpl\ub)}{\sih(\ub)} \frac{ \sih(\cpl \uc)}{\sih(\uc)} +  
\exp\left(\frac12\bpl\ub\right)\frac{\sih (\cmin\uc)}{\sih( \uc)} \frac{\sih (\amin\ua)}{\sih(\ua)}  \\
&&+\exp\left(\frac12 \cpl\uc + \frac12 \ub\right)\frac{\sih (\amin\ua)}{\sih( \ua)}\frac{\sih (\bpl\ub)}{\sih(\ub)}, \end{array}$$
 $\sih(u)=\exp\left(\frac{u}2\right) - \exp\left(-\frac{u}2\right)=u\left(1+\frac{u^2}{24}+O(4) \right)$
and
$\frac{\ssih \left(\amin\ua\right)}{\ssih (\ua)}=
\amin\left(1+\frac{\amin^2-1}{24}\ua^2\right)+O(4).$
Since the degree $2$ part of $\exp\left(\frac12\amin\ua\right)\frac{\ssih (\bpl\ub)}{\ssih(\ub)} \frac{ \ssih(\cpl \uc)}{\ssih(\uc)}$ is
$$\bpl\cpl\left(\frac18 \amin^2\ua^2 + \frac{\bpl^2-1}{24}\ub^2+ \frac{\cpl^2-1}{24}\uc^2\right),$$
we get
$$\begin{array}{lll}\CD_2(\petita,\petitb,\petitc)&=&\frac18 \left( \amin^2\bpl\cpl \ua^2 +\bpl^2\cmin\amin \ub^2 + \cpl^2\amin\bpl \uc^2 + \amin\bpl \ub^2 + 2\cpl\amin\bpl \ub\uc  \right)\\
&&+\frac{\amin(\bpl+\cmin)}{24}(\amin^2-1)\ua^2 
+\frac{\bpl(\cpl+\amin)}{24}(\bpl^2-1)\ub^2 +   
\frac{\bpl\cpl(\cpl^2-1)+\amin\cmin(\cmin^2-1)}{24} \uc^2\\
&=&\frac18 \left( \amin^2\bpl\cpl \ua^2 +\bpl^2\cmin\amin \ub^2 + \cpl^2\amin\bpl \uc^2 + \amin\bpl \ub^2 + \cpl\amin\bpl (\ua^2-\ub^2-\uc^2)  \right)\\
&&+\frac{\amin\apl(\petitb+\petitc)}{48}(\amin-1)\ua^2 +\frac{\bpl\bmin(\petitc+\petita)}{48}(\bpl+1)\ub^2 
+\frac{\cmin\cpl (\bpl(\cpl+1)+\amin(\cmin-1))}{24} \uc^2.\end{array}$$
Thanks to the symmetry of $\CD_2(\petita,\petitb,\petitc)=d_a \ua^2 +d_b \ub^2 +d_c \uc^2$, it suffices to compute 
$$\begin{array}{ll}d_a&=\frac{\amin\apl}{48}\left(6\bpl\cpl +(\petitb+\petitc)(\amin-1)\right) =\frac{\petita^2-1}{8\times 48}\left(3 (\petitb+1)(\petitc+1) +(\petitb+\petitc)(\petita -3)\right)\\&=\frac{\petita^2-1}{8\times 48}\left(3 + 3 \petitb\petitc +\petita(\petitb+\petitc)\right) \end{array}$$
to find the first expression of $\CD_2(\petita,\petitb,\petitc)$.

For the second one, recall $16\wdel(\petita,\petitb,\petitc)=2(\petita+\petitb)(\petitb+\petitc)(\petitc+\petita)-4(\petita \petitb + \petitb \petitc + \petitc \petita +1)(\petita+\petitb+\petitc)$. We have
$$\frac{\lambda^{\prime}(\petita,\petitb,\petitc)}{4}\ddel_{2,\ExtH}(\alpha,\beta)- \frac1{12} \wdel(\petita,\petitb,\petitc)\Pcanp =c_a \ua^2 +c_b \ub^2 +c_c \uc^2$$
where
$$\begin{array}{lll} c_a&=&-\frac{1}{16 \times 24}(\petita\petitb+\petita\petitc+\petitb\petitc+1)(\petita^2+2\petita\petitb+2\petita\petitc+3)\\
&&+ \frac{1}{16 \times 24}\left(2\petita (\petita^2(\petitb+\petitc) + \petitb^2(\petita+\petitc) +\petitc^2(\petita+\petitb)) +8\petita^2\petitb\petitc +4\petita(\petita+\petitb+\petitc)\right)\\
 &=&\frac{(\petita^2-1)}{16 \times 24} (3+3\petitb\petitc+\petita(\petitb+\petitc)).
  \end{array}$$
\eop

\begin{proposition}
\label{propCDtwo} Let $\Sigma$ be a surface as in Proposition~\ref{propgenDSigma}.
Set $\ddel_{\ExtH}(\alpha,\beta)=\ddel_{2,\ExtH}(\alpha,\beta)+4\ddel_{\Delta,\ExtH}(\alpha,\beta)$, where $\ddel_{2,\ExtH}(\alpha,\beta)$ and $\ddel_{\Delta,\ExtH}(\alpha,\beta)$ are respectively defined in Lemma~\ref{lemdegtwoCDprel} and in Theorem~\ref{thmAlexHg}. The degree two part $\CD_2(\Sigma)$ of $\CD(\Sigma)$
satisfies
 $$\frac{1}{|H_1(\rats)|}\CD_2(\Sigma)=
w_{SL}(\Sigma) \Pcanp+ \frac{\lambda^{\prime}(\partial \Sigma)}{4}\ddel_{\ExtH}(\alpha,\beta),$$
where $w_{SL}$ and $\Pcanp$ are respectively defined in Theorem~\ref{thmwsl} and Lemma~\ref{lemCEdegtwo}.
\end{proposition}
\bp 
According to Proposition~\ref{propcompDEK},
$\CD(\Sigma)=\CD(\petita,\petitb,\petitc)\CA_{\ExtH}\left(\tua\wedge\tub\right) + \CE(\petita,\petitb,\petitc) \Gamma.$
Lemmas~\ref{lemCDdegzero} and \ref{lemdegtwoCD} imply
$$\CD(\petita,\petitb,\petitc)=\lambda^{\prime}(\petita,\petitb,\petitc) - \frac1{12} \wdel(\petita,\petitb,\petitc)\Pcanp + \frac{\lambda^{\prime}(\petita,\petitb,\petitc)}{4}\ddel_{2,\ExtH}(\alpha,\beta) +O(4)$$ where $\lambda^{\prime}(\petita,\petitb,\petitc)=\lambda^{\prime}(\phi(K))$.
Lemma~\ref{lemCEdegtwo} implies 
$\CE(\petita,\petitb,\petitc)= -\Pcanp+O(3)$.
Theorem~\ref{thmAlexHg} implies $\CA_{\ExtH}\left(\tua\wedge\tub\right)=|H_1(\rats)|\left(1 + \ddel_{\Delta,\ExtH}(\alpha,\beta) \right) +O(3)$,
and Proposition~\ref{propGamma} implies $\Gamma=-|H_1(\rats)|\lambda^{\prime}(\ExtH;\ua,\ub,\uc) +O(2)$.
Therefore, we have
$$\frac{1}{|H_1(\rats)|}\CD_2(\Sigma)= \frac{\lambda^{\prime}(\petita,\petitb,\petitc)}{4}\bigl(\ddel_{2,\ExtH}(\alpha,\beta)+4\ddel_{\Delta,\ExtH}(\alpha,\beta)\bigr)- \frac1{12} \wdel(\petita,\petitb,\petitc)\Pcanp  +\lambda^{\prime}(\ExtH;\ua,\ub,\uc)\Pcanp,$$
where 
$w_{SL}(\Sigma)$ is defined to be $w_{SL}(\Sigma)=\lambda^{\prime}(\ExtH;\ua,\ub,\uc)- \frac1{12} \wdel(\petita,\petitb,\petitc)  $ in Theorem~\ref{thmwsl}.
\eop

Recall 
$$\begin{array}{llll}W_s=V_s \circ (\fHone^{-1} \otimes_s \fHone^{-1}) \colon &\otimes_s^2 H_1(\ExtH=\rats\setminus (\mathring{H}))/\mbox{Torsion}
&\rightarrow &\QQ\\
&\ua^2 &\mapsto& \petitb+\petitc\\
&\ua\ub &\mapsto& -\petitc
\end{array}$$ from the statement of Theorem~\ref{thmalexcanphi}.
Note the following easy lemma.

\begin{lemma} \label{lemreleasy} The polynomial $\wdel(\petita,\petitb,\petitc)$ of Theorem~\ref{thmsimpleinvt} equals
$$\wdel(\petita,\petitb,\petitc)=- \frac{\petita^2(\petitb+\petitc) +\petitb^2(\petitc+\petita) +\petitc^2(\petita+\petitb)}{ 
8}-\frac{\petita \petitb \petitc}{2}-\frac{\petita+\petitb+\petitc}{4}$$
and we have $W_s(\ddel_{2,\ExtH}(\alpha,\beta))= \wdel(a,b,c)$,
$W_s(\ddel_{\Delta,\ExtH}(\alpha,\beta))=a\lambda^{\prime}(\phi(\granda)) +  
b\lambda^{\prime}(\phi(\grandb))  + c\lambda^{\prime}(\phi(\grandc))$,
$W_s(\ddel_{\ExtH}(\alpha,\beta))=\wdel(\Sigma)$,
and 
$W_s(\Pcanp)=4\lambda^{\prime}(\phi(K))-1$.
\end{lemma}
 \bp Thanks to Lemma~\ref{lemdegtwoCDprel}, we have
 $$W_s(\ddel_{2,\ExtH}(\alpha,\beta))= \frac{\petita^2(\petitb+\petitc) +\petitb^2(\petitc+\petita) +\petitc^2(\petita+\petitb)}{ 
24} - \frac{\petita+\petitb+\petitc}{6}(4\lambda^{\prime}(\phi(K)) -1)
 -\frac{\petita+\petitb+\petitc}{4}.$$
 \eop

 Theorem~\ref{thmalexcanphi} is a direct consequence of Proposition~\ref{propCDtwo} and Lemma~\ref{lemreleasy}. \eop
 
 \subsection{The invariant \texorpdfstring{$\ddel$}{delta} of a genus one Seifert surface}
 \label{subdelta}
 
 In this section, we prove that $\left(\ddel_{\Sigma}(\alpha,\beta)=\fHone_{\ast}^{-1}\bigl(\ddel_{\ExtH}(\alpha,\beta)\bigr)\right)$
is an invariant of $\Sigma$ such that $\wdel(\Sigma)=V_s(\ddel_{\Sigma})$, as stated in Theorem~\ref{thmsecondinvariantbis}.
 
 \begin{remark}
Corollary~\ref{corwww} implies that $w_{SL}(\Sigma)$ is an invariant of $\Sigma$ when $\Sigma$ is null-homologous.
Since $\Pcanp$ is also an invariant of $\Sigma$, Proposition~\ref{propCDtwo} implies Theorem~\ref{thmsecondinvariantbis} when $\Sigma$ is null-homologous and $\lambda^{\prime}(K) \neq 0$. The proof below shows that Theorem~\ref{thmsecondinvariantbis} holds without these restrictive assumptions.
\end{remark}

 First note the following lemma.
 
 \begin{lemma}
The polynomial $\ddel_{\ExtH}(\alpha,\beta)$ of Proposition~\ref{propCDtwo} depends only on the homology classes of $\alpha$ and $\beta$ on $\Sigma$.
\end{lemma}
\bp Since $a$, $b$, and $c$ are determined by the Seifert form, they depend only on the homology classes of $\alpha$ and $\beta$ on $\Sigma$. So do $\lambda^{\prime}(\phi(A))$, $\lambda^{\prime}(\phi(B))$, and $\lambda^{\prime}(\phi(C))$ according to Lemma~\ref{lemhomcurvesisot}.
\eop

The group of diffeomorphisms of $\Sigma$ up to isotopy acts on $H_1(\Sigma)=\ZZ[\alpha] \oplus \ZZ[\beta]$ as $SL_2(\ZZ)$. It is generated by the Dehn twists about the curves $\alpha$ and $\beta$ of Figure~\ref{figSigmaabc}. Thus it is sufficient to prove that $\ddel_{\Sigma}(\alpha,\beta)$ is invariant under a Dehn twist about $\alpha$ and a Dehn twist about $\beta$ in order to prove Theorem~\ref{thmsecondinvariantbis}.
Let $\tau_{\alpha}$ be the Dehn twist about a curve parallel to $\alpha$ that changes $\beta$ to $\beta^{\prime}= \beta \alpha$.
Thanks to the order $3$ symmetry of the setting, it suffices to prove that $\ddel_{\Sigma}(\alpha,\beta)$ is invariant under $\tau_{\alpha}$.

The twist $\tau_{\alpha}$ changes $\betam$ to $\beta^{\prime}_-= \beta_- \alpham$. It leaves $\alpha$ and $\alpham$ unchanged. Therefore, it leaves the curve 
$\ua=\alpha\alpham^{-1}$ of Figures~\ref{fighandlebobasedu} and \ref{figabstractSigmaalphabet} unchanged ($\ua^{\prime}=\ua$), and it changes $\ub =\betam^{-1}\beta$ to $\ub^{\prime}=\alpham^{-1}\ub\alpha$. Recall $\calpham=\ub^{-1}\alpham\ub$ from Remark~\ref{rkDeltaAtilde}. Then we have $$\ub^{\prime}=\ub\calpham^{-1}\alpha.$$

\begin{lemma}\label{diffdeltaDel} We have
$$\ddel_{\Delta,\ExtH}(\alpha^{\prime},\beta^{\prime}) - \ddel_{\Delta,\ExtH}(\alpha,\beta)=
 \frac1{8}(b^2-1)\ua\uc.$$
\end{lemma}
\bp Since $\tub^{\prime}=\tub + \exp(\ub - \alpham)(\talpha-\tcalpham)$, we have
$$\CA_{\ExtH}(\tua^{\prime}\wedge\tub^{\prime})=\CA_{\ExtH}(\tua\wedge\tub)+\exp(\ub-\alpham)\CA_{\ExtH}(\tua\wedge\talpha-\tcalpham)$$
where 
$\partial \tub \CA_{\ExtH}(\tua\wedge\talpha)= \partial \tua \CA_{\ExtH}(\tub\wedge\talpha) + \partial \talpha \CA_{\ExtH}(\tua\wedge\tub)$
according to Lemma~\ref{lemthreedel}.
Using the similar equality for $\calpham$ and the expressions
of $\CA_{\ExtH}(\tub\wedge\talpha)$ and $\CA_{\ExtH}(\tub\wedge\tcalpham)$ given by Corollary~\ref{coreasyone} and Remark~\ref{rkDeltaAtilde}, we get
$$\begin{array}{lll}\sih(\uc)\partial \tub\CA_{\ExtH}(\tua^{\prime}\wedge\tub^{\prime})&=&\sih(\uc)\partial \tub\CA_{\ExtH}(\tua\wedge\tub)\\
&&+\exp(\ub-\alpham)\sih(\uc)\exp(\alpham)\left(\exp(\ua)-1\right)\CA_{\ExtH}(\tua\wedge\tub)\\
&&+
\partial \tua\exp(\ub-\alpham)
\left(\exp\left( \frac{1+2\cmin}{2}\uc\right) -\exp\left( \frac{1+2\cpl}{2}\uc\right)\right) 
\overline{\CA_{\ExtH}(\tua\wedge\tub)}.\end{array}$$
 So, we obtain
 $$\begin{array}{lll}\partial \tub \CA_{\ExtH}(\tua^{\prime}\wedge\tub^{\prime})&=&\left(\partial\tub+\exp(\ub)\partial\tua\right)\CA_{\ExtH}(\tua\wedge\tub)\\
 &&-\partial \tua \exp(\ub-\alpham) \exp(\cpl \uc) \overline{\CA_{\ExtH}(\tua\wedge\tub)} \\
 &=&(\exp(-\uc)-1)\CA_{\ExtH}(\tua\wedge\tub) -\partial \tua \exp((1+\bmin)\ub)\overline{\CA_{\ExtH}(\tua\wedge\tub)}.\end{array}$$
 Set $\CA_{\ExtH}^{\prime}=\exp(\bmin \ua)\CA_{\ExtH}$. Then we have
 $$\begin{array}{lll}\sih(\ub) \CA_{\ExtH}^{\prime}(\tua^{\prime}\wedge\tub^{\prime})&=&-
\sih(\uc)\exp\left(\left(\frac12+\bmin\right) \ua\right)\CA_{\ExtH}(\tua\wedge\tub) \\ &&
-\sih( \ua )\exp\left(\left(-\frac12 - \bmin\right)\uc\right)\overline{\CA_{\ExtH}(\tua\wedge\tub)}\\
&=&-\sih(\uc)\exp\left(\frac{\petitb}2 \ua\right)\CA_{\ExtH}(\tua\wedge\tub)-\sih( \ua )\exp\left(-\frac{\petitb}2\uc\right)\overline{\CA_{\ExtH}(\tua\wedge\tub)},\end{array}
 $$
 and the degree two part of $\sih(\ub)\CA_{\ExtH}^{\prime}(\tua^{\prime}\wedge\tub^{\prime})$ is zero. Therefore, the degree one part of $\CA_{\ExtH}^{\prime}(\tua^{\prime}\wedge\tub^{\prime})$ is zero, and 
 $\ddel^{\prime}=\ddel_{\Delta,\ExtH}(\alpha^{\prime},\beta^{\prime})$ is the degree two part of $\frac{1}{|H_1(\rats)|}\CA_{\ExtH}^{\prime}(\tua^{\prime}\wedge\tub^{\prime})$.
 Set $\ddel=\ddel_{\Delta,\ExtH}(\alpha,\beta)$. We have
 $$\begin{array}{lll}\ub \left(1 +\frac{\ub^2}{24}\right)(1+\ddel^{\prime})&=&-\uc
  \left(1 +\frac{\uc^2}{24}\right)\left(1 +\frac{\petitb^2}8\ua^2\right)(1+{\ddel})\\
  &&-\ua
  \left(1 +\frac{\ua^2}{24}\right)\left(1 +\frac{\petitb^2}8\uc^2\right)(1+{\ddel}) + O(4)\\
  &=&\ub\left(1 + \ddel + \frac{\uc^2-\ua\uc +\ua^2}{24} +\frac18\petitb^2  \ua\uc \right) + O(4).\end{array}$$
So
 $\ddel^{\prime}=\ddel +\frac18\petitb^2\ua\uc + \frac{\uc^2-\ua\uc +\ua^2- (\ua+\uc)^2}{24}=\ddel +\frac18(\petitb^2-1)\ua\uc$.
 \eop

 \begin{lemma}
 \label{lemvarDehnone}
 The Dehn twist $\tau_{\alpha}$ about $\alpha$ changes $(\petita,\petitb,\petitc)$ to $$(\petita^{\prime}=\petita +2 \petitb,\petitb^{\prime}= 2\petitb +\petitc,\petitc^{\prime}= -\petitb)$$
 and the element  $(\ua,\ub,\uc)$ of $H_1(\partial \ExtH)^3$ to $(\ua^{\prime}=\ua,\ub^{\prime}=\ub +\ua=-\uc,\uc^{\prime}=\uc -\ua)$.
\end{lemma}
\bp In $H_1(\partial \ExtH)$, we have $\alpham= \cpl \uc-\bmin \ub$,
$\alpha^{\prime}=\alpha $, $\beta^{\prime}=\beta + \alpha$, and $\gamma^{\prime}=\gamma -\alpha$.
  We also have 
  $lk(\betam,\gamma)=-\apl$, $lk(\gammam,\alpha)=-\bpl=-\frac{\petitb+1}{2}$, and $lk(\alpham,\beta)=-\cpl$.
  Since $lk(\gamma^{\prime}_-, \alpha) = -\bpl -\frac{\petitb+\petitc}{2}$, we obtain
  $\petitb^{\prime} = 2\petitb +\petitc.$
  Similarly $lk(\alpha_-, \beta^{\prime})=-\cpl^{\prime}=-\cpl + \frac{\petitb+\petitc}{2}$, and $\petitc^{\prime} = -\petitb$. Since $lk(\beta^{\prime}_-,\gamma^{\prime})=-\apl - \frac{\petitb+\petitc}{2} -\bmin +  \cmin$, we get
$\petita^{\prime} = \petita +2 \petitb.$
\eop 

\begin{lemma} \label{lemdeltatwovar} We have
$$\ddel_{2,\ExtH}(\alpha^{\prime},\beta^{\prime}) - \ddel_{2,\ExtH}(\alpha,\beta)=
 \frac12 (1-\petitb^2)\ua\uc$$
\end{lemma}
\bp        
Set $\ddel^{\prime\prime}_2=\ddel_{2,\ExtH}(\alpha^{\prime},\beta^{\prime})- \ddel_{2,\ExtH}(\alpha,\beta)$ and recall
$$\ddel_{2,\ExtH}(\alpha,\beta)= \frac{\petita^2\ua^2 +\petitb^2\ub^2 +\petitc^2\uc^2}{24} - 
\frac{\petita+\petitb+\petitc}{6}\Pcanp
 -\frac{\ua^2 +\ub^2 +\uc^2}{8}$$
from Lemma~\ref{lemdegtwoCDprel}. We have 
$$\begin{array}{lll}\ddel^{\prime\prime}_2&=&-\frac{\petitb}{3}\Pcanp
-\frac{1}{8}\left(\ua^2-\ub^2 +\uc^2- (\ub^2-\ua^2-\uc^2)\right)\\
&&+\frac{1}{24}\left((4\petitb^2+4\petita\petitb)\ua^2 +(2\petitb+\petitc)^2\uc^2 -\petitb^2 \ub^2 +b^2(\uc-\ua)^2-c^2 \uc^2\right)\\
&=&-\frac{1}{24}(4ab \ua^2 +4b^2\ub^2 + 4\petitb\petitc\uc^2)
-\frac{1}{24}\left(6\ua^2-6\ub^2 +6\uc^2\right)\\
&&+\frac{1}{24}\left((4b^2+4ab)\ua^2 +(4b^2+4bc+c^2)\uc^2 -\petitb^2 \ub^2 +\petitb^2(2\uc^2+2\ua^2-\ub^2)-\petitc^2 \uc^2\right)\\
&=&\frac{1}{24}(6\petitb^2-6)( \ua^2-\ub^2 +\uc^2)=\frac12 (1-\petitb^2)\ua\uc.
\end{array}$$
\eop

\begin{lemma}
 \label{lemvarddelExtH} Under the hypotheses of Proposition~\ref{propCDtwo}, the polynomial $\ddel_{\ExtH}(\alpha,\beta)$ is an invariant of $\Sigma$.
\end{lemma}
\bp This is a direct consequence of Lemmas~\ref{diffdeltaDel} and \ref{lemdeltatwovar}.
\eop

 \noindent{\sc Proof of Theorem~\ref{thmsecondinvariantbis}:}
Lemma~\ref{lemvarddelExtH} implies that 
$\ddel_{\Sigma}(\alpha,\beta)=\fHone_{\ast}^{-1}\bigl(\ddel_{\ExtH}(\alpha,\beta)\bigr)$ is an invariant of $\Sigma$.
\eop

\noindent{\sc Proof of Theorem~\ref{thmwsl}:}
Since $\CD_2(\Sigma)$, the polynomial $\Pcanp$ of Lemma~\ref{lemCEdegtwo}, and $\ddel_{\ExtH}(\alpha,\beta)$ are invariants of $\Sigma$, and since $\Pcanp \neq 0$, Proposition~\ref{propCDtwo} implies that $w_{SL}(\Sigma)$ is an invariant of $\Sigma$, too. \eop

\begin{lemma}
 \label{lemvarwdelDehn}
With the notation of Lemma~\ref{lemvarDehnone}, we have
$$\wdel(\petita^{\prime},\petitb^{\prime},\petitc^{\prime})-\wdel(\petita,\petitb,\petitc)
= \frac12\petitb (\petitb^2-1).$$
\end{lemma}
\bp This follows from Lemmas~\ref{lemdeltatwovar} and \ref{lemreleasy}. \eop

\section{Proof of invariance\texorpdfstring{ of $W_s\left(\CD_2(\Sigma)\right)$}{}}
\setcounter{equation}{0}

\subsection{Invariance \texorpdfstring{of $W_s\left(\CD_2(\Sigma)\right)$}{} under cobordism}
\label{secinvcob}

In this subsection, we prove that $W_s\left(\CD_2(.)\right)$ is invariant under $K$-cobordism of genus one Seifert surfaces, as stated in Lemma~\ref{leminvKcob}.

Let $\Sigma$ and $\Sigma^{\prime}$ be two null-homologous genus one Seifert surfaces of a knot $K$ in a rational homology $3$-sphere $\rats$. Assume 
$\Sigma \cap \Sigma^{\prime}=K$. Also assume that $\Sigma^{\prime} \cup_K (-\Sigma)$ is smooth, and consider a tubular neighborhood $[-1,1] \times \bigl(\Sigma^{\prime} \cup_K (-\Sigma)\bigr)$ of  $\Sigma^{\prime} \cup_K (-\Sigma)$ in $\rats$. Set 
$\Sigma^{\prime}_+=\{1\} \times \Sigma^{\prime}$, $\Sigma^{\prime}_-=\{-1\} \times \Sigma^{\prime}$, $(-\Sigma_+)=\{-1\} \times(-\Sigma)$, and $\Sigma_-=\{1\} \times(\Sigma)$.
The $3$-manifold
$\rats \setminus \bigl( \left]-1,1\right[ \times \bigl(\Sigma^{\prime} \cup_K (-\Sigma)\bigr)\bigr)$ has two connected components $Y$ and $Z$
with $\partial Y = \Sigma^{\prime}_- \cup (-\Sigma_+)$ and $\partial Z = \Sigma_- \cup (-\Sigma^{\prime}_+)$.
(This notation for $\Sigma^{\prime}_-$ and $\Sigma_-$ is compatible with the previous notation except in a tubular neighborhood of $K$.)

\begin{figure}[h]
\begin{center}
\begin{tikzpicture}
\fill[lightgray] (-4,.5) rectangle (4,1.5) (-4,-.5) rectangle (4,-1.5);
\fill[white] (-2.4,-.7) rectangle (-1.6,-1.3) (2.4,.7) rectangle (1.6,1.3);
\draw [thick] (0,0) node{\scriptsize $[-1,1] \times \left(\Sigma^{\prime} \cup (-\Sigma)\right)$}
(3,.75) node{\scriptsize $\Sigma^{\prime}_+$}
(3,-.75) node{\scriptsize $\Sigma^{\prime}_-$}
(-3,.75) node{\scriptsize $\Sigma_-$}
(-3,-.75) node{\scriptsize $\Sigma_+$}
(-2,-1) node{\scriptsize $Y$} (2,1) node{\scriptsize $Z$}
(0,.75) node{\scriptsize $\{1\} \times K$}
(0,-.75) node{\scriptsize $\{-1\} \times K$}
(-4,.5) -- (4,.5) (-4,-.5) -- (4,-.5);
\fill  (0,.5) circle (1.5pt) (0,-.5) circle (1.5pt);

\begin{scope}[xshift=8cm]
\fill[lightgray] (-2,-1.5) rectangle (2,1.5);
\draw[fill=gray] (0,0) circle (1cm);
\fill  (0,1) circle (1.5pt);
\draw (-.2,1.1) node {\scriptsize $\star$};
\draw [thick,->] (0,1.3)  node[right]{\scriptsize $K$} (0,-1.5) -- (0,1.5);
\fill[white] (-1.8,-1.3) rectangle (-1.1,-.7) (1.8,-1.3) rectangle (1.1,-.7) (.1,-.25) rectangle (.9,.25);
\draw (-1.45,-1) node{\scriptsize $\Sigma^{\prime}$}
(1.45,-1) node{\scriptsize $-\Sigma$} (.5,0) node {\scriptsize $D(\star)$};
\end{scope}
\end{tikzpicture}
\caption{The manifolds $Y$ and $Z$, and the disk $D(\star)$ in $\left(\Sigma^{\prime} \cup (-\Sigma)\right)$}
\label{figYZ} \end{center}
\end{figure}
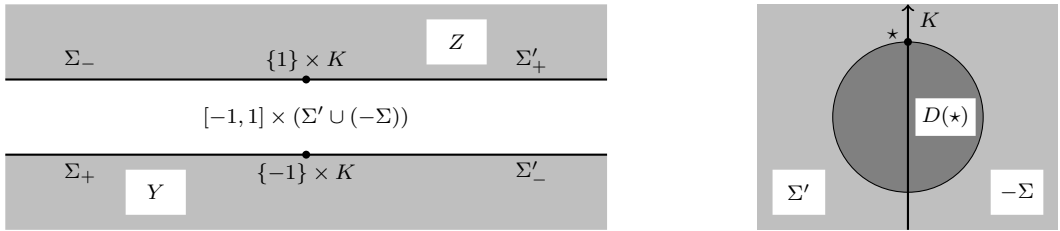

The \emph{Lagrangian} of a $3$-manifold $X$ with boundary is the kernel of the map induced by the inclusion from $H_1(\partial X;\QQ)$ to $H_1(X;\QQ)$.

Let $\star$ be a basepoint of $K$. Let $D(\star)$ be a disk of $\left(\Sigma^{\prime} \cup_K (-\Sigma)\right)$ such that $\star$ is on $\partial D(\star)$ and the intersection $D(\star) \cap K$ is an interval, as in the right part of Figure~\ref{figYZ}. Let $I_K \subset K$ be the interval $K \setminus (K\cap \mathring{D}(\star))$.

Define $Y \ast Z$ to be $Y \cup Z \cup \left([-1,1] \times D(\star)\right)$.
Then $\partial (Y \ast Z)$ separates $\rats$ into two genus $4$ rational homology handlebodies $Y \ast Z$ and $E_{YZ}\cong[-1,1] \times  \left(\Sigma^{\prime} \cup_{I_K} (-\Sigma) \right)$.
For an element $u$ of $H_1(\Sigma)$, $u_+$ denotes the class of $\{-1\} \times u$ (the curve $u$ is pushed in the direction of the positive normal to $\Sigma$, which is the negative normal to $(-\Sigma)$).
Let $\CL_{\Sigma}$ be the image of the injective map
$$\begin{array}{lll} H_1(\Sigma;\QQ) &\to &H_1(\partial (Y \ast Z);\QQ)\\
   u & \mapsto & [u_+] -[u_-],
  \end{array}
$$
 and let $\CL_{\Sigma^{\prime}}$ be the image of the similar map $(u^{\prime} \mapsto [u^{\prime}_+] -[u^{\prime}_-])$ from 
$H_1(\Sigma^{\prime};\QQ)$ to $H_1(\partial (Y \ast Z);\QQ)$.
Then the Lagrangian of $E_{YZ}$ is the direct sum $\CL_{\Sigma} \oplus \CL_{\Sigma^{\prime}}$.

The Lagrangian of $Y \ast Z$ is the direct sum of the Lagrangians $\CL_Y$ and $\CL_Z$ of the two rational homology handlebodies $Y$ and $Z$.
We have
$$\begin{array}{ll} H_1(\partial Y;\QQ)&=H_1(\Sigma^{\prime}_-;\QQ) \oplus H_1(\Sigma_+;\QQ),\\
H_1(\partial Z;\QQ)&=H_1(\Sigma_-;\QQ) \oplus H_1(\Sigma^{\prime}_+;\QQ),\\
H_1\left(\partial (Y\ast Z);\QQ\right)&= H_1(\partial Y;\QQ) \oplus H_1(\partial Z;\QQ)\\
H_1(Y;\QQ)&=H_1(\partial Y;\QQ)/\CL_Y,\\
H_1(Z;\QQ)&=H_1(\partial Z;\QQ)/\CL_Z,\\
  \end{array}
$$
 and $$H_1(Y\ast Z;\QQ)=H_1(Y;\QQ)\oplus H_1(Z;\QQ) = H_1\bigl(\partial (Y\ast Z);\QQ\bigr)/(\CL_Y \oplus \CL_Z). $$

Construct the link exterior $Y[K]$ (resp. $Z[K]$) by attaching a $2$-handle along $\{-1\} \times K$ to $Y$ (resp. $\{1\} \times K$ to $Z$). Choose a symplectic basis
$(\alpha,\beta)$ of $H_1(\Sigma)$ and a symplectic basis
$(\alpha^{\prime},\beta^{\prime})$ of $H_1(\Sigma^{\prime})$. Normalize
the Alexander forms $\CA_{Y}$ and $\CA_{Z}$ of $Y$ and $Z$, with respect to these bases,
so that the even normalized Reidemester torsions of $Y[K]$ and $Z[K]$ are
$$\CD(Y[K])=\CA_Y(\talpha_+ \wedge \tbeta_+) \in \Lambda_Y^{\QQ}=\ZZ[H_1(Y;\QQ)]$$ and
$$\CD(Z[K])=\CA_Z(\talpha^{\prime}_+ \wedge \tbeta^{\prime}_+)\in \Lambda_Z^{\QQ}=\ZZ[H_1(Z;\QQ)],$$
respectively, as in Lemma~\ref{lemcompapalone}.
As shown in Remark~\ref{rkcanon} below, there is no canonical sign for the Reidemeister torsions of $Y[K]$ and $Z[K]$.
We fix arbitrary signs so that the above equalities hold.

In the definition of $Y \ast Z$ as $Y \cup Z \cup \left([-1,1] \times D(\star)\right)$, 
shrink $[-1,1] \times  \partial D(\star)$ to $\partial D(\star)$ so that the basepoint $\star$ is in $ Y\cap Z$. 
Note that $\ExtH =\rats \setminus \left(]1,1[ \times (\Sigma \setminus K)\right)$ is obtained from $Y \ast Z$ by identifying $\Sigma^{\prime}_-$ with $\Sigma^{\prime}_+$. The manifold
$\ExtH^{\prime} =\rats \setminus \left(]1,1[ \times (\Sigma^{\prime} \setminus K)\right)$ is similarly obtained from $Y \ast Z$ by identifying $\Sigma_-$ with $\Sigma_+$. Let $i_{\ExtH} \colon Y \ast Z \hookrightarrow  \ExtH$ and $i_{\ExtH^{\prime}} \colon Y \ast Z \hookrightarrow \ExtH^{\prime}$ denote the corresponding inclusion maps (whose homotopy classes are well-defined).
Then we have the following lemma.
\begin{lemma}
We can normalize $\CA_{Y \ast Z}$ so that it has the following properties.
 \begin{itemize}
  \item The Alexander form $\CA_{Y \ast Z}$ is valued in $\Lambda_{Y\ast Z}^{\QQ}=\ZZ[H_1(Y;\QQ)\oplus H_1(Z;\QQ)]$, which contains $\Lambda_Y^{\QQ}$ and $\Lambda_Z^{\QQ}$.\\
  \item If $u_Y$, $v_Y$ are two based curves of $\partial Y$, and if $u_Z$, $v_Z$ are two based curves of $\partial Z$, then 
$$\CA_{Y \ast Z}(\tilde{u}_Y\wedge \tilde{v}_Y \wedge \tilde{u}_Z \wedge \tilde{v}_Z)=\CA_{Y}(\tilde{u}_Y\wedge \tilde{v}_Y) \CA_{Z} (\tilde{u}_Z \wedge \tilde{v}_Z).$$
\item If $u$, $v$, and $w$ are three based curves of $\partial Y$ (resp. of $\partial Z$) in $\CH_{Y \ast Z}$, then $\CA_{Y \ast Z}(\tilde{u}\wedge \tilde{v} \wedge \tilde{w} \wedge.)$ is zero.
\item There exists $\eta\in \{- 1,1\}$ such that $$\varepsilon\left(\CA_{Y \ast Z}(\talpha^{\prime}_+-\talpham^{\prime} \wedge \tbeta^{\prime}_+ -\betam^{\prime} \wedge \talpha_+ -\talpham \wedge \tbeta_+-\betam) \right)=\eta |H_1(\rats)|.$$
 \end{itemize}
Let $\mbox{ev}(\ua^{\prime}=\ub^{\prime}=0)$ (resp. $\mbox{ev}(\ua=\ub=0)$) be the map from $\Lambda_{Y\ast Z}^{\QQ}$ to $\Lambda_{\ExtH}^{\QQ}$ (resp. $\Lambda_{\ExtH^{\prime}}^{\QQ}$) induced by $i_{\ExtH}$ (resp. by $i_{\ExtH^{\prime}}$). 
Then 
$$\CD(\Sigma)=\eta \mbox{ev}(\ua^{\prime}=\ub^{\prime}=0)\left(\CA_{Y}(\talpha_+ \wedge \tbeta_+)\CA_{Z}(\talpha^{\prime}_+ \wedge \tbeta^{\prime}_+)\right)$$ 
and $$\CD(\Sigma^{\prime})=\eta\mbox{ev}(\ua=\ub=0)\left(\CA_{Y}(\talpha_+ \wedge \tbeta_+)\CA_{Z}(\talpha^{\prime}_+ \wedge \tbeta^{\prime}_+)\right).$$
\end{lemma}
\bp The first four properties are easy to prove. The Alexander form $\CA_{\ExtH^{\prime}}$ of
$\ExtH^{\prime}=\rats \setminus \left(\left]-1,1\right[ \times (\Sigma^{\prime} \setminus K)\right)$ is normalized so that
$\varepsilon\bigl(\CA_{\ExtH^{\prime}}(\talpha^{\prime}-\talpham^{\prime}\wedge\tbeta^{\prime}-\tbetam^{\prime})\bigr)=|H_1(\rats)|$
as in Section~\ref{subintrocd}.
Since $\ExtH^{\prime}$ is obtained from $Y \ast Z$ by attaching two $2$-handles along $\ua = \alpha_+\alpham^{-1}$ and $\ub =\betam^{-1}\beta_+$, we have
$$\CA_{\ExtH^{\prime}} \doteq \eta \mbox{ev}(\ua=\ub=0)\left(\CA_{Y \ast Z}(\talpha_+-\talpham \wedge \tbeta_+ -\tbetam \wedge .)\right)$$
where $\mbox{ev}(\ua=\ub=0)$ sends the homology classes of $\ua$ and $\ub$ to zero.
Thus, we get
$$\begin{array}{ll}\CD(\Sigma^{\prime})&\doteq \eta\mbox{ev}(\ua=\ub=0)\left(\CA_{Y \ast Z}(\talpha_+-\talpham \wedge \tbeta_+ -\tbetam \wedge \talpha^{\prime}_+ \wedge \tbeta^{\prime}_+)\right)\\
   &\doteq \eta\mbox{ev}(\ua=\ub=0)\left(\CA_{Y \ast Z}(\talpha_+ \wedge \tbeta_+  \wedge \talpha^{\prime}_+ \wedge \tbeta^{\prime}_+)\right)\\
&\doteq \eta\mbox{ev}(\ua=\ub=0) \left(\CA_{Y}(\talpha_+ \wedge \tbeta_+)\CA_{Z}(\talpha^{\prime}_+ \wedge \tbeta^{\prime}_+)\right).
  \end{array}$$
 Since $\CA_{Y}(\talpha_+ \wedge \tbeta_+)\CA_{Z}(\talpha^{\prime}_+ \wedge \tbeta^{\prime}_+)$ is even, we can replace $\doteq$ by $= $.
We similarly get $\CD(\Sigma)=\eta \mbox{ev}(\ua^{\prime}=\ub^{\prime}=0)\bigl(\CA_{Y}(\talpha_+ \wedge \tbeta_+)\CA_{Z}(\talpha^{\prime}_+ \wedge \tbeta^{\prime}_+)\bigr)$.
\eop

Since $\CA_{Y}(\talpha_+ \wedge \tbeta_+)$ and $\CA_{Z}(\talpha^{\prime}_+ \wedge \tbeta^{\prime}_+)$ are both even, the following lemma allows us to conclude the proof of Lemma~\ref{leminvKcob}. 

\begin{lemma} \label{lemkeyalginvcob} With the notation of Theorem~\ref{thmalexcanphi} for $W_s$,
for any element $S$ of
$$\Bigl(H_1(\partial Y;\QQ) \otimes_s H_1(\partial Y;\QQ) \Bigr)\oplus \Bigl(H_1(\partial Z;\QQ) \otimes_s H_1(\partial Z;\QQ)\Bigr),$$ we have
$$W_s(\Sigma^{\prime})\bigl(\mbox{\rm ev}(\ua=\ub=0)(S)\bigr)=W_s(\Sigma)\bigl(\mbox{\rm ev}(\ua^{\prime}=\ub^{\prime}=0)(S)\bigr).$$
\end{lemma}
\bp It suffices to prove the lemma for an element $S=y \otimes y$ where $y \in H_1(\partial Y) $ or $y \in H_1(\partial Z)$.
Let $y \in H_1(\partial Y)=H_1(\Sigma^{\prime}_-) \oplus H_1(\Sigma_+)$. 
Since $H_1(\partial (Y \ast Z);\QQ) =\CL_Y \oplus \CL_Z \oplus \CL_{\Sigma} \oplus \CL_{\Sigma^{\prime}}$,
there exists a unique $(u,u^{\prime}) \in H_1(\Sigma;\QQ) \times H_1(\Sigma^{\prime};\QQ)$ such that
$$y=u_+-u_-+u^{\prime}_+-u^{\prime}_-$$ in $H_1(\partial (Y \ast Z);\QQ)/(\CL_Y \oplus \CL_Z)$.
Since $i_{\ExtH\ast}(y)=u_+-u_-$, we have $u=\fHone^{-1}(i_{\ExtH\ast}(y))$,
and $$W_s(\Sigma)\left(\mbox{ev}(\ua^{\prime}=\ub^{\prime}=0)(S)=i_{\ExtH\ast}(y) \otimes i_{\ExtH\ast}(y)\right)=2lk(u_+,u_-).$$
We similarly get $$W_s(\Sigma^{\prime})\left(\mbox{ev}(\ua=\ub=0)(S)\right)=2lk(u^{\prime}_+,u^{\prime}_-).$$
Since 
$(u^{\prime}_+-u_-) \in \CL_Z$, we have
$lk(u^{\prime}_+,u^{\prime}_-)=lk(u_-,u^{\prime}_-)=lk(u,u^{\prime})$,
and $lk(u_+,u_-)=lk(u_+,u^{\prime}_+)=lk(u,u^{\prime})$.
 For $y \in H_1(\partial Z;\QQ)$ such that $y=u_+-u_-+u^{\prime}_+-u^{\prime}_-$ in $H_1(Y \ast Z;\QQ)$, we similarly obtain
 $$W_s(\Sigma)\Bigl(\mbox{ev}(\ua^{\prime}=\ub^{\prime}=0)(y \otimes y)\Bigr)=W_s(\Sigma^{\prime})\Bigl(\mbox{ev}(\ua=\ub=0)(y \otimes y)\Bigr)=lk(u,u^{\prime}).$$
\eop

Lemma~\ref{leminvKcob} is proved. \eop

\begin{remark}
 The hypothesis that $S$ belongs to
$\left( \otimes_s^2 H_1(\partial Y;\QQ) \right)\oplus \left(\otimes_s^2 H_1(\partial Z;\QQ)  \right)$ is necessary in Lemma~\ref{lemkeyalginvcob}. For example, for $S=(\alpha_+-\alpha_-) \otimes (\alpha_+-\alpha_-)$, $W_s(\Sigma^{\prime})\left(\mbox{ev}(\ua=\ub=0)(S)\right)=0$, while $W_s(\Sigma)\left(\mbox{ev}(\ua^{\prime}=\ub^{\prime}=0)(S)\right)=lk(\alpha_+,\alpha_-)$ is not zero, in general.
\end{remark}

\begin{remark} \label{rkcanon} 
Fix a normalization $\CA_Y$ of the Alexander form of $Y$. For two based curves $u_1$ and $u_2$ of $Y$, the integer
$\varepsilon \left(\CA_Y(\tilde{u}_1 \wedge \tilde{u}_2) \right)$ depends only of the homology classes of $u_1$ and $u_2$ in $Y$. It is equal to $\pm \bigl| H_1(Y)/\left(\ZZ[u_1] +  \ZZ[u_2]\right)\bigr|$. We denote it by $\varepsilon \CA_Y\left({u}_1 \wedge {u}_2 \right)$.
Choose a symplectic basis $(\alpha,\beta)$ of $H_1(\Sigma)$ such that $\alpha \notin \CL_Y$.

Let $\psi \colon \Sigma \to \Sigma^{\prime}$ be a diffeomorphism.
Then $\rats(Y,\psi) =Y/(\Sigma \sim_{\psi} \Sigma^{\prime})$ is a $\QQ$-sphere if and only if $\varepsilon\CA_Y\left(\alpha_+ - \psi_{\ast}(\alpha)_- \wedge \beta_+ -\psi_{\ast}(\beta)_-\right)  \neq 0$.
Below, we omit the subscripts $+$ and $-$ because they are useless here. 

We are about to prove that there exist diffeomorphisms $\psi_r \colon \Sigma \to \Sigma^{\prime}$ and $\psi_s \colon \Sigma \to \Sigma^{\prime}$ such that $\varepsilon\CA_Y\bigl(\alpha - \psi_{r \ast}(\alpha) \wedge \beta -\psi_{r\ast}(\beta)\bigr) \varepsilon\CA_Y\bigl(\alpha - \psi_{s \ast}(\alpha) \wedge \beta -\psi_{s\ast}(\beta)\bigr)<0$.
So, such diffeomorphisms do not provide a canonical sign for the Reidemeister torsions of $Y[K]$.

Let $(k_1,k_2,v_1,v_2)$ be a basis of $H_1(\partial Y,\QQ)$ such that $(k_1,k_2)$ is a basis of $\CL_Y$, and $\langle k_i, v_j\rangle_{\partial Y}=\delta_{ij}$. ($\delta_{ii}=1$, and $\delta_{ij}=0$ if $i \neq j$.) Then $\varepsilon \CA_Y({v}_1 \wedge {v}_2) \neq 0$, and, for any two curves $u_1$ and $u_2$ of $\partial Y$, we have
$$\varepsilon \CA_Y({v}_1 \wedge {v}_2)k_1 \wedge k_2 \wedge {u}_1 \wedge {u}_2 = \varepsilon \CA_Y({u}_1 \wedge {u}_2)k_1 \wedge k_2 \wedge {v}_1 \wedge {v}_2$$
in $\bigwedge^4 H_1(\partial Y,\QQ)$.

Define $\fdel_Y({u}_1 \wedge {u}_2)= \langle k_1, u_1\rangle_{\partial Y} \langle k_2, u_2\rangle_{\partial Y} -\langle k_1, u_2\rangle_{\partial Y} \langle k_2, u_1\rangle_{\partial Y}$. 
Then $\fdel_Y({v}_1 \wedge {v}_2)=1$ and $$\fdel_Y( {u}_1 \wedge {u}_2)= \frac{\varepsilon \CA_Y({u}_1 \wedge {u}_2)}{\varepsilon \CA_Y({v}_1 \wedge {v}_2)}.$$
In particular, $\rats(Y,\psi)$ is a $\QQ$-sphere if and only $\fdel_Y\left(\alpha - \psi_{\ast}(\alpha) \wedge \beta -\psi_{\ast}(\beta)\right) \neq 0$.

If $\CL_Y \cap H_1(\Sigma;\QQ) \neq \{0\}$, then there exist a primitive element $\gamma_Y$ of $H_1(\Sigma)$ such that $\CL_Y \cap H_1(\Sigma;\QQ) =\QQ \gamma_Y$ and a primitive element $\gamma^{\prime}_Y$ of $H_1(\Sigma^{\prime})$ such that $\CL_Y \cap H_1(\Sigma^{\prime};\QQ) =\QQ \gamma^{\prime}_Y$.
Then we choose the above basis $(k_1,k_2,v_1,v_2)$ so that $k_1=\gamma_Y$ and $k_2=\gamma^{\prime}_Y$.
So, we have 
$$\fdel_Y\Bigl(\alpha - \psi_{\ast}(\alpha) \wedge \beta -\psi_{\ast}(\beta)\Bigr) = 
\langle \gamma_Y, \alpha \rangle_{\Sigma} \langle \gamma^{\prime}_Y,\psi_{\ast}(\beta)  \rangle_{\Sigma^{\prime}} - \langle \gamma_Y, \beta \rangle_{\Sigma} \langle \gamma^{\prime}_Y,\psi_{\ast}(\alpha)  \rangle_{\Sigma^{\prime}}.$$

If $\CL_Y \cap H_1(\Sigma;\QQ) = \{0\}$, then there exists a symplectic isomorphism 
$\lambda_Y \colon H_1(\Sigma;\QQ) \to H_1(\Sigma^{\prime};\QQ)$ such that
$\CL_Y= \QQ (\alpha-\lambda_Y(\alpha)) \oplus \QQ (\beta-\lambda_Y(\beta))$.
Choose 
the above basis $(k_1,k_2,v_1,v_2)$ so that $k_1=\alpha-\lambda_Y(\alpha)$ and $k_2=\beta-\lambda_Y(\beta)$.
Thus, we have 
$$\begin{array}{lll}\fdel_Y\Bigl(\alpha - \psi_{\ast}(\alpha) \wedge \beta -\psi_{\ast}(\beta)\Bigr)&
=&
\langle \lambda_Y(\alpha),\psi_{\ast}(\alpha)\rangle_{\Sigma^{\prime}}\langle \lambda_Y(\beta),\psi_{\ast}(\beta)\rangle_{\Sigma^{\prime}} \\
&&+ \bigl(1 -\langle \lambda_Y(\alpha),\psi_{\ast}(\beta)\rangle_{\Sigma^{\prime}}\bigr)\bigl(1 +\langle \lambda_Y(\beta),\psi_{\ast}(\alpha)\rangle_{\Sigma^{\prime}}\bigr) \\
&=& 2- \langle \lambda_Y(\alpha),\psi_{\ast}(\beta)\rangle_{\Sigma^{\prime}} + \langle \lambda_Y(\beta),\psi_{\ast}(\alpha)\rangle_{\Sigma^{\prime}}.\end{array}$$

Choose a symplectic basis
$(\psi_{0}(\alpha),\psi_{0}(\beta))$ of $H_1(\Sigma^{\prime})$ such that 
$\langle \psi_{0}(\alpha),\gamma^{\prime}_Y\rangle_{\Sigma^{\prime}} \neq 0$ if $\CL_Y \cap H_1(\Sigma;\QQ) \neq \{0\}$, and 
$\langle\lambda_Y(\alpha), \psi_{0}(\alpha)  \rangle_{\Sigma^{\prime}} \neq 0$ if $\CL_Y \cap H_1(\Sigma;\QQ) = \{0\}.$
This condition can be rephrased as
$\psi_{0}(\alpha) \notin \CL_Y \oplus \QQ \alpha$. 

Let $\psi_k$ be the symplectic map that transforms $(\alpha,\beta)$ to $(\psi_k(\alpha)=\psi_0(\alpha),\psi_k(\beta)=\psi_0(\beta) +k\psi_0(\alpha))$. Such a map is induced by a diffeomorphism still denoted by $\psi_k$ from $\Sigma$ to $\Sigma^{\prime}$.
Thus, there exists $R \in \NN$ such that 
for any $r >R$, we have 
$$\fdel_Y\Bigl(\alpha - \psi_{r}(\alpha) \wedge \beta -\psi_{r}(\beta)\Bigr)\fdel_Y\Bigl(\alpha - \psi_{-r}(\alpha) \wedge \beta -\psi_{-r}(\beta)\Bigr) < 0.$$

We can simultaneously get the same results for $Z$, starting with a symplectic basis $(\alpha,\beta)$ such that $\alpha \notin \CL_Y \cup \CL_Z$ and a symplectic basis
$(\psi_{0}(\alpha),\psi_{0}(\beta))$ of $H_1(\Sigma^{\prime})$ such that $\psi_{0}(\alpha)
\notin (\CL_Z \oplus \QQ \alpha) \cup (\CL_Y \oplus \QQ \alpha)$.
In particular, there exists $R_2 \in \NN$ such that for any integer $r$ in $\RR \setminus ]-R_2,R_2[$, 
 $\rats(Y,\psi_r) =Y/(\Sigma \sim_{\psi_r} \Sigma^{\prime})$ and $\rats(Z,\psi_r) =Z/(\Sigma \sim_{\psi_r} \Sigma^{\prime})$
are $\QQ$-spheres. 
\end{remark}

\subsection{From a genus one Seifert surface to another one}
\label{secKcob}

Say that two genus one Seifert surfaces $\Sigma$ and $\Sigma^{\prime}$ of $K$ in a $\QQ$-sphere $\rats$ are \emph{weakly $K$-cobordant} if there exists a sequence $(\Sigma_i)_{i=1, \dots, k}$ of genus one Seifert surfaces of $K$ such that $\Sigma=\Sigma_1$, $\Sigma^{\prime}=\Sigma_k$, and, for any $i=1, \dots, k-1$, the surfaces $\Sigma_i$ and $\Sigma_{i+1}$ are $K$-cobordant.
In this section, we prove that any two genus one Seifert surfaces of the same knot $K$ in a 
$\QQ$-sphere are weakly $K$-cobordant. This proves Lemma~\ref{lemcob}.

We start with two oriented genus one surfaces $\Sigma$ and $\Sigma^{\prime}$ with boundary $K$ in $\rats$.
 Since they induce the same framing of $K$, we can perform an isotopy of $\Sigma$ so that $\Sigma \cap \Sigma^{\prime}$ is the disjoint union of $K$ and $n$ circles embedded in the interiors of $\Sigma$ and $\Sigma^{\prime}$, for some $n \in \NN$.

We proceed by induction on $n=n(\Sigma,\Sigma^{\prime})$. When $n=0$, $\Sigma$ and $\Sigma^{\prime}$ 
are $K$-cobordant. 
Otherwise, it suffices to prove that, for any pair $(\Sigma,\Sigma^{\prime})$ of surfaces as above, there exist genus one surfaces $\tilde{\Sigma}$ and $\tilde{\Sigma}^{\prime}$, respectively weakly $K$-cobordant to $\Sigma$ and $\Sigma^{\prime}$, such that $n(\tilde{\Sigma},\tilde{\Sigma}^{\prime})<n(\Sigma,\Sigma^{\prime})$, to conclude by induction.

Let $C$ be a connected component of $\Sigma \cap \Sigma^{\prime} \setminus K$. It is oriented as in the beginning of Section~\ref{subSato}.
On $\Sigma$, the curve $C$ can be of one of the three following types:
\begin{itemize}
\item $(\Sigma, D)$: when $\pm C$ bounds a disk embedded in $\Sigma$,
\item $(\Sigma, K)$: when $\pm C$ cobounds an annulus ($[0,1] \times S^1$) embedded in $\Sigma$ with $K$
\item $(\Sigma, p)$: when $C$ does not separate $\Sigma$.
\end{itemize}
Similarly, $C$ is of one of the types $(\Sigma^{\prime}, D)$ (when $\pm C$ bounds a disk embedded in $\Sigma^{\prime}$), $(\Sigma^{\prime}, K)$ or $(\Sigma^{\prime}, p)$ on $\Sigma^{\prime}$.

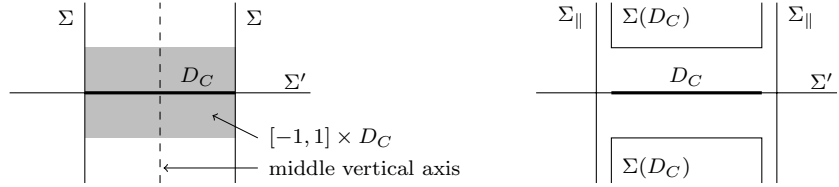
\begin{figure}[h]\begin{center}
 \begin{tikzpicture}\useasboundingbox (-5,-1.2) rectangle (5,1.2);
 \begin{scope}[xshift=-3.5cm]
\fill [lightgray] (-1,-.6) rectangle (1,.6);
\draw (-2,0) -- (2,0) (-1,-1.2) -- (-1,1.2) (1,-1.2) -- (1,1.2);
\draw [very thick] (-1,0) -- (1,0);
\draw [very thin,->] (1.3,-.6) -- (.7,-.3);
\draw [very thin,->] (1.3,-1) -- (.05,-1);
\draw [dashed] (0,-1.2) -- (0,1.2);
\draw (.5,-.05) node[above]{\scriptsize $D_C$}  (-1,1) node[left]{\scriptsize $\Sigma$} (1,1) node[right]{\scriptsize $\Sigma$} 
(1.8,-.1) node[above]{\scriptsize $\Sigma^{\prime}$}
(1.3,-.6) node[right]{\scriptsize $[-1,1] \times  D_C$}
(1.3,-1) node[right]{\scriptsize middle vertical axis};
\end{scope}
\begin{scope}[xshift=3.5cm]
\draw (-2,0) -- (2,0) (-1,-1.2) -- (-1,-.6) -- (1,-.6) -- (1,-1.2) (1,1.2) --  (1,.6) -- (-1,.6) -- (-1,1.2) (-1.2,-1.2) -- (-1.2,1.2) (1.2,-1.2) -- (1.2,1.2);
\draw [very thick] (-1,0) -- (1,0);
\draw (0,-.05) node[above]{\scriptsize $D_C$}  
(-1,1) node[right]{\scriptsize $\Sigma(D_C)$} 
(-1,-1) node[right]{\scriptsize $\Sigma(D_C)$} 
(1.2,1) node[right]{\scriptsize $\Sigma_{\parallel}$} 
(-1.2,1) node[left]{\scriptsize $\Sigma_{\parallel}$} 
(1.8,-.1) node[above]{\scriptsize $\Sigma^{\prime}$};
\end{scope}
\end{tikzpicture}
\caption{Around the disk $D_C$. (The $3$--dimensional picture is obtained by rotating the picture around the middle vertical axis.)}
\label{figsurgerDC}\end{center}
\end{figure}

\begin{lemma}
\label{lemcobone}
If $\Sigma \cap \Sigma^{\prime} \setminus K$ contains a curve of type $(\Sigma^{\prime}, D)$ or a curve of type $(\Sigma, D)$, then $\Sigma$ and $\Sigma^{\prime}$ are respectively $K$-cobordant to genus one surfaces $\tilde{\Sigma}$ and $\tilde{\Sigma}^{\prime}$ such that $n(\tilde{\Sigma},\tilde{\Sigma}^{\prime})<n(\Sigma,\Sigma^{\prime})$.
\end{lemma}
\bp
Assume that there is a curve $C$ of type $(\Sigma^{\prime}, D)$. Then there is such a curve such that the disk $D_C$ bounded by $C$ in $\pm \Sigma^{\prime}$ does not contain another component of $\Sigma \cap \Sigma^{\prime}$.
We can assume that a product $[-1,1] \times D_C$ meets $\Sigma$ along $[-1,1] \times \partial D_C$, as in Figure~\ref{figsurgerDC}, and we do.
We have $\varepsilon[-1,1] \times \partial D_C \subset \Sigma$, for some $\varepsilon \in \{-1,1\}$. Let $\Sigma(D_C)$ be the surface obtained from $\Sigma$ by replacing $\varepsilon [-1,1] \times \partial D_C$ with $\varepsilon(\partial [-1,1]) \times  D_C$. Let $\Sigma(D_C,K)$ be the connected component of $\Sigma(D_C)$ that contains $K$.

The number  $n\bigl(\Sigma(D_C,K),\Sigma^{\prime}\bigr)$ of connected components of $\Sigma(D_C,K)\cap \Sigma^{\prime}$ distinct from $K$ satisfies $n\bigl(\Sigma(D_C,K),\Sigma^{\prime}\bigr) < n(\Sigma,\Sigma^{\prime})$.
Let $\Sigma_{\parallel}$ be a surface obtained from $\Sigma$ by pushing the interior of $\Sigma$ in its normal direction that does not go towards $D_C$ near $C$. Then $\Sigma(D_C,K) \cap \Sigma_{\parallel}=K$. Thus $\Sigma(D_C,K)$ is $K$-cobordant to $\Sigma$.

When $C$ is of type $(\Sigma, D)$, $\Sigma(D_C,K)$ is a genus one surface $K$-cobordant to $\Sigma$.
Otherwise, $\Sigma(D_C,K)$ is a disk. This disk is $K$-cobordant to a genus one surface $\tilde{\Sigma}(D_C,K)$ obtained by connected sum with the boundary of a tiny solid torus near $\Sigma(D_C,K)$. We can perform such a connected sum without affecting the intersection with $\Sigma^{\prime}$ and $\Sigma_{\parallel}$.
\eop

\begin{lemma}
\label{lemcobtwo}
 Assume that $\Sigma \cap \Sigma^{\prime} \setminus K$ contains a curve
 $C$, which is of type $(\Sigma, K)$ on $\Sigma$, and of type  $(\Sigma^{\prime}, K)$ on $\Sigma^{\prime}$. Then $\Sigma$ is isotopic to a surface $\tilde{\Sigma}$, and  $\Sigma^{\prime}$ is isotopic to a surface $\tilde{\Sigma}^{\prime}$ such that $n(\tilde{\Sigma},\tilde{\Sigma}^{\prime})<n(\Sigma,\Sigma^{\prime})$.
 \end{lemma}
 \bp Consider the surfaces $\tilde{\Sigma}$ and $\tilde{\Sigma}^{\prime}$ such that $\tilde{\Sigma}$ (resp. $\tilde{\Sigma}^{\prime}$) is the closure of the surface obtained from $\Sigma$ (resp. from $\Sigma^{\prime}$) by removing the annulus between $K$ and $C$ in $\Sigma$ from $\Sigma$ (resp. in $\Sigma^{\prime}$ from $\Sigma^{\prime}$).
 \eop

\begin{lemma}
\label{lemcobthree} 
 Assume that $\Sigma \cap \Sigma^{\prime} \setminus K$ contains two curves $C_1$ and $C_2$ such that
 \begin{itemize}
 \item  there exist an annulus $A^{\prime}$ in $\Sigma^{\prime}$ and $\varepsilon_1 \in \{-1,1\}$ such that $\partial A^{\prime} = \varepsilon_1 (C_1 \cup C_2)$ and the interior of $A^{\prime}$ does not meet $\Sigma$,
 \item neither $C_1$ nor $C_2$ is of type $(\Sigma,D)$.
 \end{itemize}
 Then there exists a genus one surface $\tilde{\Sigma}$ $K$-cobordant to $\Sigma$ such that $n(\tilde{\Sigma},{\Sigma}^{\prime})<n(\Sigma,\Sigma^{\prime})$. 
\end{lemma}
 \bp Let $[-1,1] \times A^{\prime}$ be a collar of $A^{\prime}$ that intersects $\Sigma^{\prime}$ along $A^{\prime}=\{0\} \times A^{\prime}$,
 and $(\pm \Sigma)$ along $[-1,1] \times \partial A^{\prime}$.
 Our hypotheses guarantee that the coorientation of $\Sigma$ points either towards $A^{\prime}$ on $[-1,1] \times \partial A^{\prime}$ or towards $\Sigma^{\prime} \setminus A^{\prime}$.
 In other words, there exists $\varepsilon_2 \in \{-1,1\}$ such that the orientation of both components of $[-1,1] \times \partial A^{\prime}$ matches the orientation of $\varepsilon_2\Sigma$.
Set  $\check{\Sigma}= \Sigma \setminus (\varepsilon_2 \left]-1,1\right[ \times  \partial A^{\prime})$. The surface $\check{\Sigma}$ is the disjoint union of 
\begin{itemize}
 \item an annulus and a disk with two holes, when $C_1$ and $C_2$ are both of type $(\Sigma,p)$,
 \item a genus one surface with one boundary component and two annuli,
 when $C_1$ and $C_2$ are both of type $(\Sigma,K)$,
 \item an annulus and a disk with two holes otherwise, when one of the curves $C_1$ and $C_2$ is of type $(\Sigma,p)$ and the other one is of type $(\Sigma,K)$.
\end{itemize}

 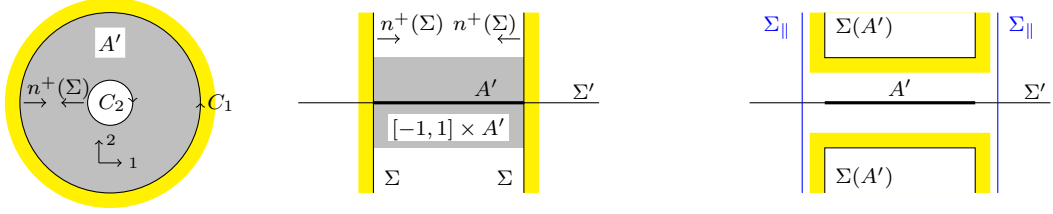
\begin{figure}[h]\begin{center}
 \begin{tikzpicture}
  \begin{scope}[xshift=-7cm]
\fill [yellow] (0,0) circle (1.4);
\fill [lightgray] (0,0) circle (1.2);
\fill [white] (0,0) circle (.3) (-.2,.6) rectangle (.2,1);
\draw (0,.8) node{\scriptsize $A^{\prime}$}
(1.15,0) node[right]{\scriptsize $C_1$}
(.35,0) node[left]{\scriptsize $C_2$}
(-.2,-.5) node[right]{\tiny $2$}
(.1,-.8) node[right]{\tiny $1$};
\draw [very thin,->] (-.15,-.8) -- (.15,-.8);
\draw [very thin,->] (-.15,-.8) -- (-.15,-.5);
\draw [->] (1.2,0) arc(0:360:1.2);
\draw [<-] (.3,0) arc(0:360:.3);
\draw [->] (-1.25,.2) node[right]{\scriptsize $n^+(\Sigma)$} 
(-.35,0) -- (-.65,0);
\draw [->]  
(-1.15,0) -- (-.85,0);
\end{scope}
 \begin{scope}[xshift=-2.5cm]
\fill [lightgray] (-1,-.6) rectangle (1,.6);
\fill [white] (-.8,-.5) rectangle (.8,-.15);
\fill [yellow] (1,-1.2) rectangle (1.2,1.2) (-1,-1.2) rectangle (-1.2,1.2);
\draw (-2,0) -- (2,0) (-1,-1.2) -- (-1,1.2) (1,-1.2) -- (1,1.2);
\draw [very thick] (-1,0) -- (1,0);
\draw (.5,-.05) node[above]{\scriptsize $A^{\prime}$}  (-1,-1) node[right]{\scriptsize $\Sigma$} (1,-1) node[left]{\scriptsize $\Sigma$} 
(1.8,-.1) node[above]{\scriptsize $\Sigma^{\prime}$}
(0,-.35) node{\scriptsize $[-1,1] \times  A^{\prime}$};
\draw [->] (-1.05,1.05) node[right]{\scriptsize $n^+(\Sigma)$} 
(-.95,.85) -- (-.65,.85);
\draw [->] (1.05,1.05) node[left]{\scriptsize $n^+(\Sigma)$} 
(.95,.85) -- (.65,.85);
\end{scope}
\begin{scope}[xshift=3.5cm]
\fill [yellow] (-1.2,.4) rectangle (1.2,1.2) (-1.2,-.4) rectangle (1.2,-1.2);
\fill [white]  (-1,.6) rectangle (1,1.2) (-1,-.6) rectangle (1,-1.2);
\draw (-2,0) -- (2,0) (-1,-1.2) -- (-1,-.6) -- (1,-.6) -- (1,-1.2) (1,1.2) --  (1,.6) -- (-1,.6) -- (-1,1.2);
\draw [blue] (-1.3,-1.2) -- (-1.3,1.2) (1.3,-1.2) -- (1.3,1.2)
(1.3,1) node[right]{\scriptsize $\Sigma_{\parallel}$} 
(-1.3,1) node[left]{\scriptsize $\Sigma_{\parallel}$};
\draw [very thick] (-1,0) -- (1,0);
\draw (0,-.05) node[above]{\scriptsize $A^{\prime}$}  
(-1,1) node[right]{\scriptsize $\Sigma(A^{\prime})$} 
(-1,-1) node[right]{\scriptsize $\Sigma(A^{\prime})$} 
(1.8,-.1) node[above]{\scriptsize $\Sigma^{\prime}$}; 
\end{scope}
\end{tikzpicture}
\caption{Around the annulus $A^{\prime}$, when the positive normals to $\Sigma$ point towards $A^{\prime}$.\\(For the second or third figure, the $3$--dimensional picture is the product of this figure by a circle. The side of the negative normal to $\Sigma$ or $\Sigma(A^{\prime})$ is colored light.)}
\label{figsurgerAprime}\end{center}
\end{figure}

Set
$\Sigma(A^{\prime})= \check{\Sigma}\cup (\varepsilon_2 (\partial [-1,1]) \times A^{\prime})$.
The surface $\Sigma(A^{\prime})$ is the union of
pieces of non-positive Euler characteristics glued along circles. 
It is an oriented surface whose boundary is $K$. Its Euler characteristic is $(-1)$.
Its connected components have non-positive Euler characteristics. Since the sum of the Euler characteristics of its components is $-1$, exactly one of its connected components has Euler characteristic $-1$, and it must be a genus one Seifert surface $\tilde{\Sigma}$ of $K$.

Furthermore, $\Sigma(A^{\prime})$ is $K$-cobordant to a surface $\Sigma_{\parallel}$ obtained from $\Sigma$ by pushing the interior of $\Sigma$ in its normal direction that does not go towards $A^{\prime}$ near $\partial A^{\prime}$. 
The component $\tilde{\Sigma}$ of $\Sigma(A^{\prime})$ is $K$-cobordant to $\Sigma_{\parallel}$, too.
 \eop

Of course, we can exchange the roles of $\Sigma$ and $\Sigma^{\prime}$
in the statement of Lemma~\ref{lemcobthree} and get the lemma.
\begin{lemma}
\label{lemcobthreebis} 
 If $\Sigma \cap \Sigma^{\prime} \setminus K$ contains two curves $C_1$ and $C_2$ such that
 \begin{itemize}
 \item  there exist an annulus $A$ in $\Sigma$ and $\varepsilon \in \{-1,1\}$ such that $\partial A = \varepsilon (C_1 \cup C_2)$ and the interior of $A$ does not meet $\Sigma^{\prime}$,
 \item neither $C_1$ nor $C_2$ is of type $(\Sigma^{\prime},D)$.
 \end{itemize}
 Then there exists a genus one surface $\tilde{\Sigma}^{\prime}$ $K$-cobordant to $\Sigma^{\prime}$ such that $n({\Sigma},\tilde{\Sigma}^{\prime})<n(\Sigma,\Sigma^{\prime})$.
\end{lemma}
\eop

The following lemma allows us to conclude the proof of Lemma~\ref{lemcob} using the previous four lemmas.

\begin{lemma} \label{lemconclucob}
Let $\Sigma$ and $\Sigma^{\prime}$ be two genus one surfaces of a knot $K$ in a $\QQ$-sphere $\rats$ such that $\Sigma \cap \Sigma^{\prime}$ is the disjoint union of $K$ and $n$ circles which are embedded in the interiors of $\Sigma$ and $\Sigma^{\prime}$, with $n \in \NN$,
and $\Sigma \cap \Sigma^{\prime} \setminus K$
contains 
\begin{itemize}
\item no curve of type $(\Sigma, D)$ on $\Sigma$, no curve of type $(\Sigma^{\prime}, D)$ on $\Sigma^{\prime}$,
 \item no curve which is both of type $(\Sigma, K)$ and $(\Sigma^{\prime}, K)$,
 \item no pair $(C_1,C_2)$ of curves (as in Lemma~\ref{lemcobthree}) such that $\pm (C_1 \cup C_2)$ is the boundary of an annulus of $\Sigma^{\prime}$ whose interior does not meet $\Sigma$,
 \item no pair $(C_1,C_2)$ of curves (as in Lemma~\ref{lemcobthreebis}) such that $\pm (C_1 \cup C_2)$ is the boundary of an annulus of $\Sigma$ whose interior does not meet $\Sigma^{\prime}$.
 \end{itemize} Then $\Sigma \cap \Sigma^{\prime}=K$.
\end{lemma}
\bp 
Since $\rats$ is a $\QQ$-sphere and since $(\Sigma^{\prime} \cup -\Sigma)$ is an oriented $2$-cycle, the algebraic intersection of any closed curve with $(\Sigma^{\prime} \cup -\Sigma)$ is zero.
Let $s$ be a point near $K$ on the negative side of $(\Sigma^{\prime} \cup -\Sigma)$ (the side of the negative normal).
For any point of $\rats \setminus (\Sigma^{\prime} \cup (-\Sigma))$,
define $W(p)$ as the algebraic intersection of a path from $s$ to $p$ transverse to 
 $(\Sigma^{\prime} \cup (-\Sigma))$ with $(\Sigma^{\prime} \cup (-\Sigma))$.
Then $\rats \setminus (\Sigma^{\prime} \cup (-\Sigma))$ has several connected components $\Omega_i$ weighted by integers $W(\Omega_i)$,
where $W(\Omega_i) =W(p \in \Omega_i)$.
The difference $W(\Omega_j)-W(\Omega_i)$ is the algebraic intersection of a path from a point of $\Omega_i$ to a point of $\Omega_j$ with $\Sigma^{\prime} \cup (-\Sigma)$. The components adjacent to $K$ are weighted by $0$ and $1$. 

\begin{figure}[h]
\begin{center}
\begin{tikzpicture}\useasboundingbox (-3.5,-1) rectangle (6.5,1);
\fill [yellow] (-3.5,0) rectangle (3.5,.2);
\fill [blue!40,fill opacity=0.5] (0,-1) rectangle (.2,1);
\draw [very thick] (-3.5,0) -- (3.5,0) (0,-1) -- (0,1);
\draw [->] (-.05,.85) node[right]{\scriptsize $n^+(\Sigma^{\prime})$} (.1,.6) -- (.6,.6);
\draw [->] (3.4,.25) node{\scriptsize $n^+(-\Sigma)$} (2.8,.1) -- (2.8,.5);
\draw (-1.75,.6) node{\scriptsize $\Omega_{+-}$} (-3.2,.25) node{\scriptsize $\Sigma$}  (-.25,-.85) node{\scriptsize $\Sigma^{\prime}$}
(.45,-.25) node{\scriptsize $\Sigma \cap \Sigma^{\prime}$} 
(1.75,-.6) node{\scriptsize $\Omega_{+-}$} (-1.75,-.6) node{\scriptsize $\Omega_{--}$} (1.75,.6) node{\scriptsize $\Omega_{++}$} 
(5,.5) node[right]{$W(\Omega_{++})=W(\Omega_{-+})+1$} (5,0) node[right]{$W(\Omega_{+-})=W(\Omega_{-+})$} (5,-.5) node[right]{$W(\Omega_{--})=W(\Omega_{-+})-1$};
\end{tikzpicture}
\caption{Weights around a component of $\Sigma \cap \Sigma^{\prime}\setminus K$}
\label{figweight}
 \end{center}
\end{figure}
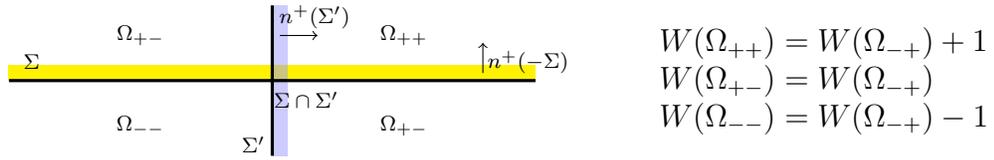

Weight the connected components of $\Sigma^{\prime} \setminus (\Sigma \cap \Sigma^{\prime})$ by the
weight of the adjacent component of $\rats \setminus (\Sigma^{\prime} \cup (-\Sigma))$ in the direction of the positive normal to $\Sigma^{\prime}$. Similarly weight the connected components of $\Sigma \setminus (\Sigma \cap \Sigma^{\prime})$ by the
weight of the adjacent component of $\rats \setminus (\Sigma^{\prime} \cup (-\Sigma))$ in the direction of the negative normal to $\Sigma$. So, the components adjacent to $K$ get the weight one, and the weight on $\Sigma \setminus (\Sigma \cap \Sigma^{\prime})$ or on $\Sigma^{\prime} \setminus (\Sigma \cap \Sigma^{\prime})$ varies by $\pm 1$ across a component of $\Sigma \cap \Sigma^{\prime}$. In particular, $\Sigma \cap \Sigma^{\prime}$ contains an even number of (necessarily parallel) components of type $(\Sigma^{\prime},p)$.
Half of them are oriented in one direction, and the other half are oriented in the opposite direction. So, if there are curves of $(\Sigma \cap \Sigma^{\prime})\setminus K$ of type $(\Sigma^{\prime},p)$, there is an annulus $A^{\prime}$ between two such consecutive transverse intersection curves $C_1$ and $C_2$ of type $(\Sigma^{\prime},p)$ with opposite direction. This is forbidden by our assumptions.
Transverse intersection curves of type $(\Sigma,p)$ are similarly forbidden. \eop

So our reduction is finished, and Lemma~\ref{lemcob} is proved. \eop

\subsection{Curves on genus one Seifert surfaces}
\label{secpfisotopcurv}

In this section, we prove that two homologous non-separating simple closed curves $\xi$ and $\xi^{\prime}$ of a genus one oriented surface $\Sigma$ with one boundary component are isotopic in $\Sigma$, as stated in Lemma~\ref{lemhomcurvesisot}.
 
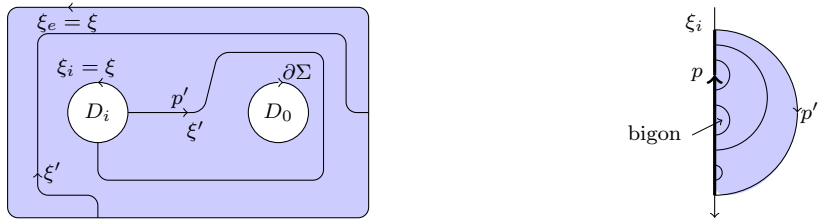
\begin{figure}[h]
\begin{center}
\begin{tikzpicture}
\begin{scope}[xshift=-4cm]
\fill [rounded corners,blue!20] (-2.4,-1.4) rectangle (2.4,1.4);
\draw [rounded corners,->] ( -1.6,1.4) -- (-2.4,1.4) -- (-2.4,-1.4) -- ( 2.4,-1.4) -- ( 2.4,1.4) -- ( -1.6,1.4);
\draw (-1.6,1.2) node{\scriptsize $\xi_e=\xi$};
\fill [white] (-1.2,0) circle (.4);
\fill [white] (1.2,0) circle (.4);
\draw [->] (-1.2,.4) arc (90:450:.4);
 \draw (-1.35,.6) node{\scriptsize $\xi_i=\xi$} (-1.2,0) node {\scriptsize $D_i$};
\draw [-<] (1.2,.4) arc (-270:90:.4);
 \draw  (1.2,0) node {\scriptsize $D_0$} (1.45,.55) node{\scriptsize $\partial \Sigma$};
\draw [rounded corners,->]  (-.1,.2) node{\scriptsize $p^{\prime}$}  (.1,-.25) node{\scriptsize $\xi^{\prime}$} (0,0) -- (.2,0) -- (.4,.8) -- (1,.8) -- (1.8,.8) --  (1.8,-.9) -- (-1.2,-.9) -- (-1.2,-.4) (-.8,0) -- (0,0);
\draw [rounded corners,->]  (-2,-.8) -- (-2,1.05)  -- (0,1.05) -- (2.1,1.05) -- (2.1,0) -- (2.4,0) (-1.8,-.8) node{\scriptsize $\xi^{\prime}$}(-1.2,-1.4) -- (-1.2,-1.1) -- (-2,-1.1) -- (-2,-.8);
\end{scope}
\begin{scope}[xshift=3cm]
\fill [rounded corners,blue!20] (0,1.1) -- (0,-1.1) arc (-90:90:1.1);
\draw[->] (0,1.2) node[left]{\scriptsize $\xi_i$} (0,1.4) -- (0,-1.4);
\draw [->] (1,0) node[right]{\scriptsize $p^{\prime}$} (0,-1.1) arc (-90:0:1.1) (0,1.1) arc (90:0:1.1);
\draw (0,-.9) arc (-90:90:.1) (0,-.5) arc (-90:90:.7) (0,-.3) arc (-90:90:.2)  (0,.3) arc (-90:90:.2);
\draw [very thick, ->] (0,.5) node[left]{\scriptsize $p$} (0,.5) -- (0,1.1) (0,-1.1) -- (0,.5);
\draw [very thin, ->] (-.3,-.3) node[left]{\scriptsize bigon} (-.3,-.3) -- (.1,-.1);
\end{scope}
\end{tikzpicture}
\caption{The surface $\Sigma$ cut along $\xi$, and a curve $\xi^{\prime}$ homologous to $\xi$, on the left.}
\label{figSigmacutxixi}
\end{center}
\end{figure}

Assume that these curves $\xi$ and $\xi^{\prime}$ are transverse without loss of generality.
Cut $\Sigma$ along $\xi$ to obtain a disk with two holes, as in Figure~\ref{figSigmacutxixi}.
We view this disk with two holes as $D_{e} \setminus (\mathring{D}_i \sqcup \mathring{D}_0)$, for a disk $D_{e}$ of $\RR^2$ and two disjoint disks $D_0$ and $D_i$ in the interior $\mathring{D}_e$ of $D_e$, where $\partial D_0 = -\partial \Sigma$, and $\partial D_e$ and  $\partial D_i$ are respectively identified to two copies $\xi_e$ and $\xi_i$ of $\xi$. According to the Jordan–Schoenflies theorem, any simple closed curve of $D_e$ bounds a disk in $D_e$.

If $\xi^{\prime}$ does not meet $\xi$, then the disk bounded by $\xi^{\prime}$ in $D_e$ contains the disk $D_i$ bounded by the internal copy $\xi_i$ of $\xi$. So, the curve $\xi^{\prime}$ is parallel to one of the copies of $\xi$.

Otherwise, we show that we can always reduce the intersection $\xi \cap \xi^{\prime}$ of $\xi$ and $\xi^{\prime}$ by isotopy. Since the algebraic intersection number of $\xi$ and $\xi^{\prime}$ is zero, if $\xi \cap \xi^{\prime}$ is not empty, then there exists an arc $p^{\prime}$ of $\xi^{\prime}$ from $\xi_i$ to itself on $D_{e} \setminus (\mathring{D}_i \sqcup \mathring{D}_0)$. Let $p$ be the arc of $\pm \xi_i$ with boundary $-\partial p^{\prime}$ such that the disk bounded by $\pm (p\cup p^{\prime})$ in $D_e$ does not contain $D_i$.  If this disk bounded by  $\pm (p\cup p^{\prime})$ does not contain $D_0$, then there is a bigon (as in the right-hand side of Figure~\ref{figSigmacutxixi}) that can be eliminated, and we can reduce $\xi \cap \xi^{\prime}$ by an isotopy. Otherwise, 
pick an arc $q$ from the external copy $\xi_e$ of $\xi$ to itself. The arc $q$ cannot meet the interior of the disk bounded by $\pm (p\cup p^{\prime})$. Therefore, it cobounds a disk with an arc of $\xi_e$ in $D_{e} \setminus (\mathring{D}_i \sqcup \mathring{D}_0)$. So, there is a least a bigon that can be suppressed by isotopy.  \eop

\def\cprime{$'$}
\providecommand{\bysame}{\leavevmode ---\ }
\providecommand{\og}{``}
\providecommand{\fg}{''}
\providecommand{\smfandname}{\&}
\providecommand{\smfedsname}{\'eds.}
\providecommand{\smfedname}{\'ed.}
\providecommand{\smfmastersthesisname}{M\'emoire}
\providecommand{\smfphdthesisname}{Th\`ese}

 \end{document}